\documentclass{article}%
\usepackage{graphicx}
\usepackage{amsmath}
\usepackage{amsfonts}
\usepackage{amssymb}%
\newtheorem{theorem}{Theorem}[section]
\newtheorem{corollary}[theorem]{Corollary}
\newtheorem{lemma}[theorem]{Lemma}
\newtheorem{proposition}[theorem]{Proposition}
\newenvironment{proof}[1][Proof]{\textbf{#1.} }{\ \rule{0.5em}{0.5em}}
\begin{document}

\title{On finitely generated profinite groups, I: strong completeness and uniform bounds}
\author{Nikolay Nikolov\thanks{Work done while the first author held a Golda-Meir
Fellowship at the Hebrew University of Jerusalem} \ and Dan Segal}
\maketitle

\begin{abstract}
We prove that in every finitely generated profinite group, every subgroup of
finite index is open; this implies that the topology on such groups is
determined by the algebraic structure. This is deduced from the main result
about finite groups: let $w$ be a `locally finite' group word and
$d\in\mathbb{N}$. Then there exists $f=f(w,d)$ such that in every
$d$-generator finite group $G$, every element of the verbal subgroup $w(G)$ is
equal to a product of $f$ $w$-values.

An analogous theorem is proved for commutators; this implies that in every
finitely generated profinite group, each term of the lower central series is closed.

The proofs rely on some properties of the finite simple groups, to be
established in Part II.

\end{abstract}

\section*{Contents}

\noindent\S 1. \ Introduction

\noindent\S 2. \ The Key Theorem

\noindent\S 3. \ Variations on a theme

\noindent\S 4. \ Proof of the Key Theorem

\noindent\S 5. \ The first inequality: lifting generators

\noindent\S 6. \ Exterior squares and quadratic maps

\noindent\S 7. \ The second inequality, soluble case

\noindent\S 8. \ Word combinatorics

\noindent\S 9. \ Equations in semisimple groups, 1: the second inequality

\noindent\S 10. Equations in semisimple groups, 2: powers

\noindent\S 11. Equations in semisimple groups, 3: twisted commutators

\section{Introduction}

A profinite group $G$ is the inverse limit of some inverse system of finite
groups. Thus it is a compact, totally disconnected topological group;
properties of the original system of finite groups are reflected in properties
of the topological group $G$. An algebraist may ask: does this remain true if
one forgets the topology? Now a base for the neighbourhoods of $1$ in $G$ is
given by the family of all open subgroups of $G$, and each such subgroup has
finite index; so if \emph{all} subgroups of finite index were open we could
reconstruct the topology by taking these as a base for the neighbourhoods of
$1$.

Following \cite{RZ} we say that $G$ is \emph{strongly complete} if it
satisfies any of the following conditions, which are easily seen to be equivalent:

\begin{description}
\item[(a)] every subgroup of finite index in $G$ is open,

\item[(b)] $G$ is equal to its own profinite completion,

\item[(c)] every group homomorphism from $G$ to any profinite group is continuous.
\end{description}

This seems \emph{a priori} an unlikely property for a profinite group, and it
is easy to find counterexamples. Indeed, any countably based but not finitely
generated pro-$p$ group will have $2^{2^{\aleph_{0}}}$ subgroups of index $p$
but only countably many open subgroups; more general examples are given in
\cite{RZ}, \S 4.2, and some examples of a different kind will be indicated
below. Around 30 years ago, however, J.-P. Serre showed that every
\emph{finitely generated} pro-$p$ group is strongly complete. We generalize
this to

\begin{theorem}
\label{strong}Every finitely generated profinite group is strongly complete.
\end{theorem}

(Here, `finitely generated' is meant in the topological sense.) \ This answers
Question 7.37 of the 1980 Kourovka Notebook \cite{K}, restated as Open
Question 4.2.14 in \cite{RZ}. It implies that the topology of a finitely
generated profinite group is completely determined by its underlying abstract
group structure, and that the category of finitely generated profinite groups
is a full subcategory of the category of (abstract) groups.

The theorem is a consequence of our major result. This concerns \emph{finite}
groups having a bounded number of generators, and the values taken by certain
\emph{group words}. Let us say that a group word $w$ is $d$\emph{-locally
finite} if every $d$-generator (abstract) group $H$ satisfying $w(H)=1$ is
finite (in other words, if $w$ defines a variety of groups all of whose
$d$-generator groups are finite).

\begin{theorem}
\label{general-w}Let $d$ be a natural number, and let $w$ be a group word.
Suppose \emph{either} that $w$ is $d$-locally finite \emph{or} that $w$ is a
simple commutator$.$ Then there exists $f=f(w,d)$ such that: in any finite
$d$-generator group $G$, every product of $w$-values in $G$ is equal to a
product of $f$ $w$-values.
\end{theorem}

\noindent Here, by `simple commutator' we mean one of the words
\begin{align*}
\lbrack x_{1},x_{2}]  &  =x_{1}^{-1}x_{2}^{-1}x_{1}x_{2},\\
\lbrack x_{1},\ldots,x_{n}]  &  =[[x_{1},\ldots,x_{n-1}],x_{n}]\qquad(n>2),
\end{align*}
and a $w$\emph{-value} means an element of the form $w(g_{1},g_{2}%
,\ldots)^{\pm1}$ with the $g_{j}\in G$.

\subsection*{Profinite results}

The proof of Serre's theorem sketched in \S 4.2 of \cite{Sr} proceeds by
showing that if $G$ is a finitely generated pro-$p$ group then the subgroup
$G^{p}[G,G]$, generated (algebraically) by all $p$th powers and commutators,
is open in $G$. To state an appropriate generalization, consider a group word
$w=w(x_{1},\ldots,x_{k})$. For any group $H$ the corresponding \emph{verbal
subgroup} is%
\[
w(H)=\left\langle w(h_{1},\ldots,h_{k})\mid h_{1},\ldots,h_{k}\in
H\right\rangle ,
\]
the subgroup generated (\emph{algebraically}, whether or not $H$ is a
topological group) by all $w$-values in $H$. We prove

\begin{theorem}
\label{w-closed}Let $w$ be a $d$-locally finite group word and let $G$ be a
$d$-generator profinite group. Then the verbal subgroup $w(G)$ is open in $G$.
\end{theorem}

To deduce Theorem 1.1, let $G$ be a $d$-generator profinite group and $K$ a
subgroup of finite index in $G$. Then $K$ contains a normal subgroup $N$ of
$G$ with $G/N$ finite. Now let $F$ be the free group on free generators
$x_{1},\ldots,x_{d}$ and let%
\[
D=\bigcap_{\theta\in\Theta}\ker\theta
\]
where $\Theta$ is the (finite) set of all homomorphisms $F\rightarrow G/N$.
Then $D$ has finite index in $F$ and is therefore finitely generated: say%
\[
D=\left\langle w_{1}(x_{1},\ldots,x_{d}),\ldots,w_{m}(x_{1},\ldots
,x_{d})\right\rangle .
\]
It follows from the definition of $D$ that $w_{i}(\mathbf{u})\in D$ for each
$i$ and any $\mathbf{u}\in F^{(d)}$; so putting%
\[
w(\mathbf{y}_{1},\ldots,\mathbf{y}_{m})=w_{1}(\mathbf{y}_{1})\ldots
w_{m}(\mathbf{y}_{m})
\]
where $\mathbf{y}_{1},\ldots,\mathbf{y}_{m}$ are disjoint $d$-tuples of
variables we have $w(F)=D$. This implies that the word $w$ is $d$-locally
finite; and as $w_{i}(\mathbf{g})\in N$ for each $i$ and any $\mathbf{g}\in
G^{(d)}$ we also have $w(G)\leq N$. Theorem 1.3 now shows that $w(G)$ is an
open subgroup of $G$, and as $K\geq N\geq w(G)$ it follows that $K$ is open.

The statement of Theorem \ref{w-closed} is really the concatenation of two
facts: the deep result that $w(G)$ is \emph{closed} in $G$, and the triviality
that this entails $w(G)$ being open. To get the latter out of the way, say $G$
is generated (topologically) by $d$ elements, and let $\mu(d,w)$ denote the
order of the finite group $F_{d}/w(F_{d})$ where $F_{d}$ is the free group of
rank $d$. Now suppose that $w(G)$ is closed. Then $w(G)=\bigcap\mathcal{N}$
\ where $\mathcal{N}$ is the set of all open normal subgroups of $G$ that
contain $w(G)$. For each $N\in\mathcal{N}$ the finite group $G/N$ is an
epimorphic image of $F_{d}/w(F_{d})$, hence has order at most $\mu(d,w)$; it
follows that $\mathcal{N}$ is finite and hence that $w(G)$ is open.

Though not necessarily relevant to Theorem \ref{strong}, the nature of other
verbal subgroups may also be of interest. Using a variation of the same
method, we shall prove

\begin{theorem}
\label{[comm]}Let $G$ be a finitely generated profinite group and $H$ a closed
normal subgroup of $G$. Then the subgroup $[H,G]$, generated (algebraically)
by all commutators $[h,g]=h^{-1}g^{-1}hg$ ($h\in H,\,g\in G$), is closed in
$G$.
\end{theorem}

\noindent This implies that the (algebraic) derived group $\gamma
_{2}(G)=[G,G]$ is closed, and then by induction that each term $\gamma
_{n}(G)=[\gamma_{n-1}(G),G]$ of the lower central series of $G$ is closed. It
is an elementary (though not trivial) fact that $\gamma_{n}(G)$ is actually
the verbal subgroup for the word $\gamma_{n}(x_{1},x_{2},\ldots,x_{n}%
)=[x_{1},x_{2},\ldots,x_{n}]$.

\begin{theorem}
\label{powers}Let $q\in\mathbb{N}$ and let $G$ be a finitely generated
non-universal profinite group. Then the subgroup $G^{q}$ is open in $G$.
\end{theorem}

\noindent Here $G^{q}$ denotes the subgroup generated (algebraically) by all
$q$th powers in $G$, and $G$ is said to be \emph{non-universal} if there
exists at least one finite group that is not isomorphic to any open section
$B/A$ of $G$ (that is, with $A\vartriangleleft B\leq G$ and $A$ open in $G$).
We do not know whether this condition is necessary for Theorem \ref{powers};
it seems to be necessary for our proof.

Although the word $w=x^{q}$ is not in general locally finite, we may still
infer that $G^{q}$ is open once we know that $G^{q}$ is \emph{closed} in $G$.
The argument is exactly the same as before; far from being a triviality,
however, it depends in this case on Zelmanov's theorem \cite{Z} which asserts
that there is a finite upper bound $\overline{\mu}(d,q)$ for the order of any
\emph{finite} $d$-generator group of exponent dividing $q$ (the solution of
the restricted Burnside problem).

The words $\gamma_{n}$ (for $n\geq2)$ are also not locally finite. Could it be
that verbal subgoups of finitely generated profinite groups are in general
closed? The answer is no: Romankov \cite{R} has constructed a finitely
generated (and soluble) pro-$p$ group $G$ in which the second derived group
$G^{\prime\prime}$ is not closed; and $G^{\prime\prime}=w(G)$ where
$w=[[x_{1},x_{2}],[x_{3},x_{4}]]$.

\subsection*{Uniform bounds for finite groups}

Qualitative statements about profinite groups may often be interpreted as
quantitative statements about (families of) finite groups. For example, a
profinite group $G$ is finitely generated if and only if there exists a
natural number $d$ such that every continuous finite quotient of $G$ can be
generated by $d$ elements.

To re-interpret the theorems stated above, consider a group word
$w=w(x_{1},x_{2},\ldots,x_{k})$. For any group $G$ we write%
\[
G^{\{w\}}=\left\{  w(g_{1},g_{2},\ldots,g_{k})^{\pm1}\mid g_{1},g_{2}%
,\ldots,g_{k}\in G\right\}  ,
\]
and call this the set of $w$-values in $G$. If the group $G$ is profinite, the
mappings $\mathbf{g}\mapsto w(\mathbf{g})$ and $\mathbf{g}\mapsto
w(\mathbf{g})^{-1}$ from $G^{(k)}$ to $G$ are continuous, so the set
$G^{\{w\}}$ is compact. For any subset $S$ of $G$ let us write%
\[
S^{\ast n}=\{s_{1}s_{2}\ldots s_{n}\mid s_{1},\ldots,s_{n}\in S\}.
\]
Then for each natural number $n$, the set $\left(  G^{\{w\}}\right)  ^{\ast
n}$ of all products of $n$ $w$-values in $G$ is compact, hence closed in $G$.

Now $w(G)$ is the ascending union of compact sets%
\[
w(G)=\bigcup_{n=1}^{\infty}\left(  G^{\{w\}}\right)  ^{\ast n}.
\]
If $w(G)$ is closed in $G$, a straightforward application of the Baire
category theorem (see \cite{Hr}) shows that for some finite $n$ one has%
\begin{equation}
w(G)=\left(  G^{\{w\}}\right)  ^{\ast n}. \label{w(G)bounded}%
\end{equation}
The \emph{converse} (which is more important here) is obvious. Thus $w(G)$ is
closed if and only if (\ref{w(G)bounded}) holds for some natural number $n$.
Now this is a property that can be detected in the finite quotients of $G$.
That is,

\begin{itemize}
\item $w(G)=\left(  G^{\{w\}}\right)  ^{\ast n}$ if and only if
$w(G/N)=\left(  (G/N)^{\{w\}}\right)  ^{\ast n}$ for every open normal
subgroup $N$ of $G$.
\end{itemize}

\noindent The ``only if'' is obvious; to see the other implication, write
$\mathcal{N}$ for the set of all open normal subgroups of $G$ and observe that
if $w(G/N)=\left(  (G/N)^{\{w\}}\right)  ^{\ast n}$ for each $N\in\mathcal{N}$
then%
\[
w(G)\subseteq\bigcap_{N\in\mathcal{N}}w(G)N=\bigcap_{N\in\mathcal{N}}\left(
G^{\{w\}}\right)  ^{\ast n}N=\left(  G^{\{w\}}\right)  ^{\ast n}%
\]
because $\left(  G^{\{w\}}\right)  ^{\ast n}$ is closed.

It follows that Theorem \ref{w-closed} is equivalent to

\begin{theorem}
\label{w-finite} Let $d$ be a natural number and let $w$ be a $d$-locally
finite word. Then there exists $f=f(w,d)$ such that in every finite
$d$-generator group $G$, every element of the verbal subgroup $w(G)$ is a
product of $f$ $w$-values.
\end{theorem}

A similar argument shows that Theorem \ref{[comm]} is a consequence of

\begin{theorem}
\label{[H,G]}Let $G$ be a finite $d$-generator group and $H$ a normal subgroup
of $G$. Then every element of $[H,G]$ is equal to a product of $g(d)$
commutators $[h,y]$ with $h\in H$ and $y\in G$, where $g(d)=12d^{3}+O(d^{2})$
depends only on $d$.
\end{theorem}

In particular, this shows that in any finite $d$-generator group $G$, each
element of the derived group $\gamma_{2}(G)=[G,G]$ is equal to a product of
$g(d)$ commutators. Now let $n>2$. It is easy to establish identities of the
following type: (a) $[y_{1},\ldots,y_{n}]^{-1}=[y_{2},y_{1},y_{3}^{\prime
},\ldots,y_{n}^{\prime}]$ where $y_{j}^{\prime}$ is a certain conjugate of
$y_{j}$ for $j\geq3$, and (b) for $k\geq2$, $[c_{1}\ldots c_{k},x]=[c_{1}%
^{\prime},x_{1}^{\prime}]\ldots\lbrack c_{k}^{\prime},x_{k}^{\prime}]$ where
$c_{j}^{\prime}$ is conjugate to $c_{j}$ and $x_{j}^{\prime}$ is conjugate to
$x$ for each $j$. Using these and arguing by induction on $n$ we infer that
each element of $\gamma_{n}(G)=[\gamma_{n-1}(G),G]$ is a product of
$g(d)^{n-1}$ terms of the form $[y_{1},\ldots,y_{n}]$.
Thus Theorems \ref{w-finite} and \ref{[H,G]} together imply Theorem
\ref{general-w}.

For a finite group $H$ let us denote by $\alpha(H)$ the largest integer $k$
such $H$ involves the alternating group $\mathrm{Alt}(k)$ (i.e. such that
$\mathrm{Alt}(k)\cong M/N$ for some $N\vartriangleleft M\leq H$). Evidently, a
profinite group $G$ is non-universal if and only if the numbers $\alpha
(\widetilde{G})$ are bounded as $\widetilde{G}$ ranges over all the finite
continuous quotients of $G$, and we see that Theorem \ref{powers} is
equivalent to

\begin{theorem}
\label{finitepowers}Let $q$, $d$ and $c$ be natural numbers. Then there exists
$h=h(c,d,q)$ such that in every finite $d$-generator group $G$ with
$\alpha(G)\leq c$, every element of $G^{q}$ is a product of $h$ $q$th powers.
\end{theorem}

It is worth remarking (though not surprising) that the functions $f,\,g$ and
$h$ necessarily depend on the number of generators $d$ (i.e. they must be
unbounded as $d\rightarrow\infty$). This can be seen e.g. from the examples
constructed by Holt in \cite{Ho}, Lemma 2.2: among these are finite groups $K$
(with $\alpha(K)=5$) such that $K=[K,K]=K^{2}$ but with $\log\left|  K\right|
/\log\left|  K^{\{w\}}\right|  $ unbounded, where $w(x)=x^{2}$ (for the
application to $g$ note that every commutator is a product of three squares).
The Cartesian product $G$ of infinitely many such groups is then a
topologically perfect profinite group (i.e. $G$ has no proper open normal
subgroup with abelian quotient), but the subgroup $G^{2}$ is not closed; in
particular $G>G^{2}$ so $G$ contains a (non-open) subgroup of index $2$.

The proofs depend ultimately on two theorems about finite simple groups. We
state these here, but postpone their proofs, which rely on the Classification
and use quite different methods, to Part II \cite{NS}.

Let $\alpha,\,\beta$ be automorphisms of a group $S$. For $x,\,y\in S,$ we
define the \textquotedblleft twisted commutator\textquotedblright%
\[
T_{\alpha,\beta}(x,y)=x^{-1}y^{-1}x^{\alpha}y^{\beta},
\]
and write $T_{\alpha,\beta}(S,S)$ for the \emph{set} $\{T_{\alpha,\beta
}(x,y)\mid x,\,y\in S\}$ (in contrast to our convention that $[S,S]$ \ denotes
the \emph{group} generated by all $[x,y]$). Recall that a group $S$ is said to
be \emph{quasisimple} if $S=[S,S]$ and $S/\mathrm{Z}(S)$ is simple (here
$\mathrm{Z}(S)$ denotes the centre of $S$).

\begin{theorem}
\label{twisted}There is an absolute constant $D\in\mathbb{N}$ such that if $S$
is a finite quasisimple group and $\alpha_{1},\,\beta_{1},\ldots,\alpha
_{D},\,\beta_{D}$ \ are any automorphisms of $S$ then%
\[
S=T_{\alpha_{1},\beta_{1}}(S,S)\cdot\ldots\cdot T_{\alpha_{D},\beta_{D}%
}(S,S).
\]

\end{theorem}

\begin{theorem}
\label{autos}Let $q$ be a natural number. There exist natural numbers $C=C(q)$
and $M=M(q)$ such that if $S$ is a finite quasisimple group with $\left|
S/\mathrm{Z}(S)\right|  >C$, $\beta_{1},\ldots,\beta_{M}$ \ are any
automorphisms of $S$, and $q_{1},\ldots,q_{M}$ are any divisors of $q$, then
there exist inner automorphisms $\alpha_{1},\ldots,\alpha_{M}$ of $S$ such
that%
\[
S=[S,(\alpha_{1}\beta_{1})^{q_{1}}]\cdot\ldots\cdot\lbrack S,(\alpha_{M}%
\beta_{M})^{q_{M}}].
\]

\end{theorem}

(Here the notation $[S,\gamma]$ stands for the \emph{set} of all
$[x,\gamma],\,x\in S$, not the group they generate.)

\subsection*{Arrangement of the paper}

The rest of the paper is devoted to the proofs of Theorems \ref{w-finite},
\ref{[H,G]} and \ref{finitepowers}. All groups henceforth will be assumed
finite (apart from the occasional appearance of free groups).

In \S 2 we state what we call the \emph{Key Theorem}, a slightly more
elaborate version of Theorem \ref{[H,G]}, and show that it implies Theorem
\ref{w-finite}. Once this is done, we can forget all about the mysterious word
$w$. Section 3 presents two variants of the Key Theorem, and the deduction of
Theorems \ref{[H,G]} and \ref{finitepowers}.

The proof of the Key Theorem is explained in \S 4. The argument is by
induction on the group order, and the inductive step requires a number of
subsidiary results. These are established in \S \S 5, 7, 9, 10 and 11, while
Sections 6 and 8 contain necessary preliminaries. (To see just the complete
proof of Theorem \ref{strong}, the reader may skip \S 3, the last subsection
of \S 4 and \S 11.)

\subsection*{Historical remarks}

The special cases of Theorems \ref{strong}, \ref{[comm]} and \ref{powers}
relating to prosoluble groups were established in \cite{Sg}, and the global
strategy of our proofs follows the same model.

The special case of Theorem \ref{finitepowers} where $q$ is odd was the main
result of \cite{N1}. Theorem \ref{finitepowers} for simple groups $G$ (the
result in this case being independent of $\alpha(G)$) was obtained by
\cite{MZ} and \cite{SW}; a common generalization of this result and of Theorem
\ref{w-finite} for simple groups is given in \cite{LS2}, and is the starting
point of our proof. Theorem \ref{twisted} generalizes a result from \cite{W}.

The material of Sections 6 and 7 generalizes (and partly simplifies) methods
from \cite{Sg}, while that of Sections 8-11 extends techniques introduced in
\cite{N1} and \cite{N2}.

We are indebted to J. S. Wilson for usefully drawing our attention to the
verbal subgroup $w(G)$ where $w$ defines the variety generated by a finite group.

\subsection*{Notation}

Here $G$ denotes a group, $x\in G$, $\,y\in G$ or $y\in\mathrm{Aut}%
(G),\,S,\,T\subseteq G$, $q\in\mathbb{N}$.
\begin{align*}
x^{y}  &  =y^{-1}xy,\qquad\,[x,y]=x^{-1}x^{y}\\
\lbrack S,y]  &  =\left\{  [s,y]\mid s\in S\right\} \\
\mathfrak{c}(S,T)  &  =\left\{  [s,t]\mid s\in S,\,t\in T\right\} \\
S^{\{q\}}  &  =\left\{  s^{q}\mid s\in S\right\} \\
ST  &  =\left\{  st\mid s\in S,\,t\in T\right\} \\
S^{\ast q}  &  =\{s_{1}s_{2}\ldots s_{q}\mid s_{1},\ldots,s_{q}\in S\}\\
&  =SS\ldots S\text{ (}q\text{ factors),}%
\end{align*}
and $\left\langle S\right\rangle $ denotes the subgroup generated by $S$. If
$H,\,K\leq G$ (meaning that $H$ and $K$ are subgroups of $G$),%
\begin{align*}
\lbrack H,K]  &  =[H,_{1}K]=\left\langle \mathfrak{c}(H,K)\right\rangle ,\\
\lbrack H,_{n}K]  &  =[[H,_{n-1}K],K]\qquad(n>1),\\
\lbrack H,_{\omega}K]  &  =\bigcap_{n\geq1}[H,_{n}K],\\
H^{\prime}  &  =[H,H].
\end{align*}
The $n$th Cartesian power of a set $S$ is generally denoted $S^{(n)}$, and
$n$-tuples are conventionally denoted by boldface type: $(s_{1},\ldots
,s_{n})=\mathbf{s}$.

$\alpha(G)$ denotes the largest integer $k$ such that $G$ involves the
alternating group $\mathrm{Alt}(k).$

The term `simple group' will mean `non-abelian finite simple group'.

\section{The Key Theorem}

The following theorem is the key to the main results. We make an \emph{ad hoc}

\medskip

\noindent\textbf{Definition} Let $H$ be a normal subgroup of a finite group
$G$. Then $H$ is \emph{acceptable} if

\begin{description}
\item[(i)] $H=[H,G],$

\item[(ii)] if $Z<N$ are normal subgroups of $G$ contained in $H$ then $N/Z$
is neither a (non-abelian) simple group nor the direct product of two
isomorphic (non-abelian) simple groups.
\end{description}

\medskip

\noindent\textbf{Key Theorem }\ \emph{Let} $G=\left\langle g_{1},\ldots
,g_{d}\right\rangle $ \emph{be a finite group and }$H$ \emph{an acceptable
normal subgroup of }$G$.\emph{ Let }$q$\emph{ be a natural number. Then}%
\[
H=\left(  [H,g_{1}]\cdot\ldots\cdot\lbrack H,g_{d}]\right)  ^{\ast h_{1}%
(d,q)}\cdot(H^{\{q\}})^{\ast z(q)}%
\]
\emph{where }$h_{1}(d,q)$\emph{ and }$z(q)$\emph{ depend only on the indicated
arguments.}

\medskip

Assuming this result, let us prove Theorem \ref{w-finite}. Fix an integer
$d\geq2$ and a group word $w=w(x_{1},\ldots,x_{k})$; we assume that%
\[
\mu:=\mu(d,w)=\left|  F_{d}/w(F_{d})\right|
\]
is finite, where $F_{d}$ denotes the free group of rank $d$. Let $q$ denote
the order of $C/w(C)$ where $C$ is the infinite cyclic group. Evidently
$q\mid\mu$, and it is easy to see that $C^{q}=w(C)=C^{\{w\}}$; hence
\begin{equation}
h^{q}\in H^{\{w\}} \label{xpower}%
\end{equation}
for any group $H$ and $h\in H$.

Let $\mathcal{S}$ denote the set of simple groups $S$ that satisfy $w(S)=1$.
It follows from the Classification that every simple group can be generated by
two elements
; therefore $\left|  S\right|  \mid\mu(2,w)$ for each $S\in\mathcal{S}$, so
the set $\mathcal{S}$ is \emph{finite}. We shall denote the complementary set
of simple groups by $\mathcal{T}$.

An important special case of our theorem was established by Liebeck and Shalev
(it is valid for arbitrary words $w$; in the present case, it may also be
deduced, via (2), from the main result of [MZ] and [SW], together with the
fact that there are only finitely many simple groups of exponent dividing $q$):

\begin{proposition}
\label{LiSha}\emph{(\cite{LS2}, Theorem 1.6.) }There exists a constant $c(w)$
such that%
\[
S=(S^{\{w\}})^{\ast c(w)}%
\]
for every $S\in\mathcal{T}.$
\end{proposition}

The next result is due to Hamidoune:

\begin{lemma}
\label{hami}\emph{\cite{Hm}} Let $X$ be a generating set of a group $G$ such
that $1\in X$ and $\left|  G\right|  \leq r\left|  X\right|  $. Then
$G=X^{\ast2r}$.
\end{lemma}

We call a group $Q$ \emph{semisimple} if $Q$ is a direct product of simple
groups, and \emph{quasi-semisimple }if $Q=Q^{\prime}$ and $Q/\mathrm{Z}(Q)$ is
semisimple. In this case, $Q$ is a central quotient of its universal covering
group $\widetilde{Q}$, and $\widetilde{Q}$ is a direct product of quasisimple groups.

\begin{corollary}
\label{qss}Let $Q$ be a quasi-semisimple group having no composition factors
in $\mathcal{S}$. Then%
\[
Q=(Q^{\{w\}})^{\ast n_{1}}%
\]
where $n_{1}=2q^{2}c(w)+q$.
\end{corollary}

\begin{proof}
In view of the preceding remark, we may assume that $Q$ is in fact
quasisimple. Write $Z=\mathrm{Z}(Q)$ and put $X=Q^{\{w\}}$. It is evident that
$X$ generates $Q$ modulo $Z$; since $\left\langle X\right\rangle
\vartriangleleft Q=Q^{\prime}$ it follows that $X$ generates $Q$. According to
Proposition \ref{LiSha} we have $Q=ZX^{\ast c}$ where $c=c(w)$.

Now it follows from the Classification (see \cite{G}, Table 4.1 or \cite{GLS},
\S 6.1) that $Z$ has rank at most $2$. If we assume for the moment that
$Z^{q}=1$, we may infer that $\left|  Z\right|  \leq q^{2},$ so $\left|
Q\right|  \leq q^{2}\left|  X^{\ast c}\right|  $. In this case, Hamidoune's
lemma yields $Q=X^{\ast2q^{2}c}$. In general we may conclude that%
\[
Q=Z^{q}\cdot X^{\ast2q^{2}c},
\]
and the result follows since $Z$ is abelian and every $q$th power is a $w$-value.
\end{proof}

\begin{lemma}
\label{nilp-comm}Let $G$ be a group, $H$ a normal subgroup and suppose that
$G=G^{\prime}\left\langle x_{1},\ldots,x_{m}\right\rangle $. Then%
\[
\lbrack H,G]=[H,x_{1}]\ldots\lbrack H,x_{m}][H,_{n}G]
\]
for every $n\geq1$.
\end{lemma}

\begin{proof}
Suppose this holds for a certain value of $n\geq1$. To deduce that it holds
with $n+1$ in place of $n$ we may as well assume that $[H,_{n+1}G]=1$. This
implies that $[[H,_{n-1}G],G^{\prime}]=1$. Now $[H,_{n}G]$ is generated by
elements of the form $[w,g]$ with $w\in\lbrack H,_{n-1}G]$ and $g\in G$. As
$[H,_{n}G]$ is central in $G$ it follows that every element of $[H,_{n}G]$
takes the form
\[
z=[w_{1},x_{1}]\ldots\lbrack w_{m},x_{m}]
\]
with $w_{i}\in\lbrack H,_{n-1}G]$ for each $i$. For any $h_{1},\ldots,h_{m}\in
H$ we then have%
\begin{align*}
\lbrack h_{1},x_{1}]\ldots\lbrack h_{m},x_{m}]\cdot z  &  =[w_{1},x_{1}%
][h_{1},x_{1}]\ldots\lbrack w_{m},x_{m}][h_{m},x_{m}]\\
&  =[w_{1}h_{1},x_{1}]\ldots\lbrack w_{m}h_{m},x_{m}],
\end{align*}
again because each $[w_{i},x_{i}]$ is central. Thus%
\[
\lbrack H,G]=[H,x_{1}]\ldots\lbrack H,x_{m}][H,_{n}G]=[H,x_{1}]\ldots\lbrack
H,x_{m}]
\]
as required.
\end{proof}

\bigskip

\textbf{Proof of Theorem \ref{w-finite}}.\ \ Let $G$ be a $d$-generator finite
group and put $\mathfrak{X}=G^{\{w\}}$. We shall show that%
\begin{equation}
w(G)=\mathfrak{X}^{\ast f}, \label{X*f}%
\end{equation}
where $f=f(w,d)$ is a number that will be specified in due course.

We begin by setting up a configuration to which the Key Theorem may be
applied. Set%
\begin{align*}
G_{1}  &  =w(G),\\
H_{1}  &  =\bigcap_{\theta\in\Theta}\ker\theta
\end{align*}
where $\Theta$ is the set of all homomorphisms from $G_{1}$ to $\mathrm{Aut}%
(S\times S)$ with $S\in\mathcal{S}.$ Set%
\[
H_{2}=[H_{1},_{\omega}G_{1}].
\]
Then $H_{1}/H_{2}$ is nilpotent and $H_{2}=[H_{2},G_{1}]$. Define $H_{3}$ to
be the smallest normal subgroup of $H_{1}$ such that $H_{1}/H_{3}$ is soluble;
then $H_{3}\leq H_{2}$ and $H_{3}=H_{3}^{\prime}$. Set%
\[
H_{4}=\bigcap\mathcal{N}%
\]
where $\mathcal{N}$ is the set of all normal subgroups $K$ of $H_{3}$ such
that $H_{3}/K\in\mathcal{T}$. Finally, put%
\[
H_{5}=[H_{4},H_{3}].
\]
Note that $H_{3}/H_{4}$ is a semisimple group; it follows that $H_{3}/H_{5}$
is quasi-semisimple.

Next, we choose a nice generating set for $G_{1}$. Since $F_{d}/w(F_{d})$ is
finite, the group $w(F_{d})$ is generated by finitely many $w$-values in
$F_{d}$:%
\[
w(F_{d})=\left\langle w(\mathbf{u}_{1}),\ldots,w(\mathbf{u}_{d^{\prime}%
})\right\rangle .
\]
Choose an epimorphism $\pi:F_{d}\rightarrow G$ and put $g_{i}=\pi
(w(\mathbf{u}_{i}))$ for $i=1,\ldots,d^{\prime}$. Then%
\[
G_{1}=w(G)=\left\langle g_{1},\ldots,g_{d^{\prime}}\right\rangle
\]
and for each $i$ we have $g_{i}=w(\pi(\mathbf{u}_{i}))\in\mathfrak{X}$. Note
that $d^{\prime}$ depends only on $w$ and $d$, and that%
\begin{equation}
\lbrack h,g_{i}]=g_{i}^{-h}g_{i}\in\mathfrak{X}^{\ast2} \label{[h,g]in X*2}%
\end{equation}
for each $i$ and any $h\in G$.

Now we build up to the proof of (\ref{X*f}) in steps.

\medskip\textbf{Step 1}\emph{.} $H_{5}\subseteq\mathfrak{X}^{\ast n_{2}}$
\emph{where} $n_{2}=2d^{\prime}h_{1}(d^{\prime},q)+z(q)$. We show first that
$H_{5}$ \emph{is an acceptable subgroup of} $G_{1}$. To verify condition (i),
observe that $H_{5}=[H_{5},H_{3}]$ because $H_{3}=H_{3}^{\prime},$ so
$H_{5}=[H_{5},G_{1}]$. For condition (ii), suppose that $Z<N$ are normal
subgroups of $G_{1}$ contained in $H_{5}$ and that $N/Z=S_{1}\times
\cdots\times S_{n}$ where $n\leq2$ and the $S_{j}$ are isomorphic simple
groups. If $S_{1}\in\mathcal{S}$ then $H_{1}$ must act trivially by
conjugation on $N/Z$, which is impossible since $N\leq H_{1}$ and $N/Z$ is
non-abelian. Therefore $S_{1}\in\mathcal{T}$. Now $H_{3}$ permutes the factors
$S_{j}$ by conjugation, and as $H_{3}=H_{3}^{\prime}$ and $n\leq2$ it follows
that $S_{1}\vartriangleleft H_{3}/Z$. Since the outer automorphism group of
$S_{1}$ is soluble (Schreier's conjecture, \cite{G},
the action of $H_{3}$ on $S_{1}$ induces precisely the group of inner
automorphisms of $S_{1}$; consequently $H_{3}/\mathrm{C}_{H_{3}}(S_{1})\cong
S_{1}$. Hence $\mathrm{C}_{H_{3}}(S_{1})\geq H_{4}\geq N$, a contradiction
since $S_{1}$ is non-abelian.

We may now apply the Key Theorem to the pair $(G_{1},H_{5})$. This shows that
each element of $H_{5}$ is equal to one of the form%
\[
\prod_{j=1}^{h_{1}(d^{\prime},q)}\prod_{i=1}^{d^{\prime}}[a_{ij},g_{i}%
]\cdot\prod_{j=1}^{z(q)}b_{j}^{q},
\]
and the claim follows by (\ref{[h,g]in X*2}) and (\ref{xpower}).

\medskip\textbf{Step 2.} $H_{3}\subseteq\mathfrak{X}^{\ast n_{1}}H_{5}$
\emph{where} $n_{1}=2q^{2}c(w)+q$. This follows from Corollary \ref{qss}
applied to the quasi-semisimple group $H_{3}/H_{5}$.

\medskip\textbf{Step 3.} $H_{2}\subseteq\mathfrak{X}^{\ast n_{2}}H_{3}$. It is
clear that $H_{2}/H_{3}$ is an acceptable subgroup of $G_{1}/H_{3}$. The claim
now follows just as in Step 1, on applying the Key Theorem to the pair
$(G_{1}/H_{3},H_{2}/H_{3})$.

\medskip\textbf{Step 4. }$[H_{1},G_{1}]H_{1}^{q}\subseteq\mathfrak{X}%
^{\ast(2d^{\prime}+1)}H_{2}$. Note that $H_{2}=[H_{1},_{n}G_{1}]$ for some
$n$; now Lemma \ref{nilp-comm}, with (\ref{[h,g]in X*2}), shows that
$[H_{1},G_{1}]\subseteq\mathfrak{X}^{\ast2d^{\prime}}H_{2}$, and the claim
follows by (\ref{xpower}) since $H_{1}^{q}\subseteq\lbrack H_{1},G_{1}]\cdot
H_{1}^{\{q\}}$.

\medskip\textbf{Step 5.} $G_{1}\subseteq\mathfrak{X}^{\ast n_{3}}[H_{1}%
,G_{1}]H_{1}^{q}$ \emph{where }$n_{3}$ \emph{depends only on }$d^{\prime}%
$\emph{ and }$w$.\emph{ }Let $\nu$ denote the maximal order of $\mathrm{Aut}%
(S\times S)$ as\emph{ }$S$ ranges over\emph{ }$\mathcal{S}$ (it is easy to see
that $\nu\leq2\mu(2,w)^{4}$.) For each such $S$ the number of homomorphisms
$G_{1}\rightarrow\mathrm{Aut}(S\times S)$ is at most $\nu^{d^{\prime}}$, so
$\left|  G_{1}:H_{1}\right|  \leq\nu^{\nu^{d^{\prime}}}=\rho$, say. It follows
that $H_{1}$ can be generated by $\rho d^{\prime}$ elements, and hence that
$\left|  H_{1}:[H_{1},G_{1}]H_{1}^{q}\right|  \leq q^{\rho d^{\prime}}$. Thus
$\left|  G_{1}/[H_{1},G_{1}]H_{1}^{q}\right|  \leq n_{3}$ where $n_{3}=q^{\rho
d^{\prime}}\rho$; consequently each of its elements can be written as a word
of length at most $n_{3}$ in the images of the generators $g_{i}$.

\medskip\textbf{Conclusion.} Putting Steps 1 -- 5 together we obtain
(\ref{X*f}) with%
\[
f=n_{1}+2n_{2}+2d^{\prime}+1+n_{3}.
\]

\section{Variations on a theme}

In this section we present two variants of the Key Theorem, and use them to
deduce Theorems \ref{[H,G]} and \ref{finitepowers}. The variants will be
proved at the end of \S 4.

The first variant of the Key Theorem has the same hypotheses, but a new
conclusion (its proof will not need Theorem 1.10 or the material of \S 10).

\medskip

\noindent\textbf{Key Theorem }(B)\textbf{ } \emph{Let} $G=\left\langle
g_{1},\ldots,g_{d}\right\rangle $ \emph{be a finite group and }$H$ \emph{an
acceptable normal subgroup of }$G$.\emph{ Then}%
\[
H=\left(  [H,g_{1}]\cdot\ldots\cdot\lbrack H,g_{d}]\right)  ^{\ast h_{2}%
(d)}\cdot\mathfrak{c}(H,H)^{\ast D}%
\]
\emph{where }$h_{2}(d)=6d^{2}+O(d)$\emph{ depends only on }$d$\emph{ and }%
$D$\emph{ is an absolute constant (given in Theorem 1.9).}

\medskip

\textbf{Proof of Theorem \ref{[H,G]}.} Let $G=\left\langle g_{1},\ldots
,g_{d}\right\rangle $ be a finite group and $H$ a normal subgroup of $G$.
Putting%
\[
\mathfrak{X}=\mathfrak{c}(H,G)
\]
we shall show that%
\begin{equation}
\lbrack H,G]=\mathfrak{X}^{\ast g(d)} \label{[H,G]=Xstar}%
\end{equation}
where $g(d)\leq2dh_{2}(d)+O(d)$ is a number that depends only on $d$.
Obviously, if $H$ is acceptable this follows at once from Key Theorem (B),
with $g(d)=dh_{2}(d)+D$. For the general case, we take a step by step approach
as in the preceding section.

Put%
\[
H_{1}=[H,_{\omega}G];
\]
let $H_{2}$ be the smallest normal subgroup of $H$ such that $H/H_{2}$ is
soluble; let
\[
H_{3}=\bigcap\mathcal{N}%
\]
where $\mathcal{N}$ denotes the set of all normal subgroups $K$ of $H_{2}$
such that $H_{2}/K$ is (non-abelian) simple; and put%
\[
H_{4}=[H_{3},H_{2}].
\]
As in the preceding section, we see that $H_{4}=[H_{4},H_{2}]=[H_{4},G]$ and
that $H_{2}/H_{4}$ is a quasi-semisimple group. We shall need

\begin{lemma}
\label{c(Q,Q)}If $Q$ is a quasi-semisimple group then $Q=\mathfrak{c}%
(Q,Q)^{\ast D}$.
\end{lemma}

\noindent Replacing $Q$ by its universal cover, we may suppose that $Q$ is a
direct product of quasisimple groups; in that case, the result follows from
the special case of Theorem 1.9 where all the automorphisms $\alpha_{j}$ and
$\beta_{j}$ are equal to the identity (this special case may be quickly
deduced, using Lemma \ref{hami}, from Wilson's theorem \cite{W}, Prop. 2.1).
\medskip

\textbf{Step 1}$_{B}$. $H_{4}\subseteq\mathfrak{X}^{\ast(dh_{2}(d)+D)}$. As
remarked above, this holds provided $H_{4}$ is an acceptable normal subgroup
of $G$. That this is the case follows, just as in Step 1 of the preceding
section, from the fact that $H_{3}$ is contained in the kernel of every
homomorphism $H_{2}\rightarrow\mathrm{Aut}(S\times S)$, $S$ any simple group;
the argument is now much simpler since we may ignore the distinction made
there between different kinds of simple group.\medskip

\textbf{Step 2}$_{B}$. $H_{2}\subseteq\mathfrak{X}^{\ast D}H_{4}$. This
follows from Lemma \ref{c(Q,Q)} applied to the quasi-semisimple group
$H_{2}/H_{4}$.\medskip

\textbf{Step 3}$_{B}$. $H_{1}\subseteq\mathfrak{X}^{\ast(dh_{2}(d)+D)}H_{2}$.
This follows from Key Theorem (B) applied to the pair $(G/H_{2},H_{1}/H_{2})$;
it is clear that $H_{1}/H_{2}$ is an acceptable normal subgroup of $G/H_{2}%
$.\medskip

\textbf{Step 4}$_{B}$. $[H,G]\subseteq\mathfrak{X}^{\ast d}H_{1}$. This is
immediate from Lemma \ref{nilp-comm}.\medskip

\textbf{Conclusion.} Putting the steps together we obtain (\ref{[H,G]=Xstar})
with%
\[
g(d)=2dh_{2}(d)+3D+d=12d^{3}+O(d^{2}).
\]
$\blacksquare$

\bigskip

The second variant of the Key Theorem has a weaker hypothesis: as we shall
see, this is necessary because the failure of the word $w(x)=x^{q}$ to be
locally finite means that we have less control over the generators of the
verbal subgroup $G^{q}$. (The proof of this variant will not need Theorem 1.10
or the material of \S \S 10, 11.)

\medskip

\noindent\textbf{Key Theorem }(C) \ \emph{Let} $G$\emph{ be a }$d$%
\emph{-generator finite group and }$H$ \emph{an acceptable normal subgroup of
}$G$. \emph{Suppose that }$G=H\left\langle g_{1},\ldots,g_{r}\right\rangle
$\emph{. Then}%
\[
H=\left(  [H,g_{1}]\cdot\ldots\cdot\lbrack H,g_{r}]\right)  ^{\ast h_{3}(d,c)}%
\]
\emph{where }$h_{3}(d,c)$\emph{ depends only on }$d$\emph{ and }$c=\alpha(G)$.

\bigskip

\textbf{Proof of Theorem \ref{finitepowers}.} Let $G$ be a $d$-generator group
with $\alpha(G)\leq c$, let $q$ be a natural number, and put $\mathfrak{X}%
=G^{\{q\}}$. We will prove that%
\begin{equation}
G^{q}=\mathfrak{X}^{\ast h}, \label{Gq=Xh}%
\end{equation}
where $h=h(c,d,q)$ will be determined below.

To this end, we take $w(x)=x^{q}$ and then define $G_{1}=G^{q}$ and normal
subgroups $H_{1}\geq\ldots\geq H_{5}$ exactly as in the proof of Theorem 1.6
in \S 2. The argument now follows that proof step by step, but we have to
carry out the steps in reverse order: this is necessary in order to obtain
substitutes for the `global generators' $g_{i}$ used in \S 2.

As in the preceding section, we will repeatedly use the fact that if $h\in G$
and $g\in\mathfrak{X}$ then $[h,g]\in\mathfrak{X}^{\ast2}$.

Set%
\[
\mu=\overline{\mu}(d,q),
\]
the maximal order of a finite $d$-generator group of exponent dividing $q$;
this is finite by the positive solution of the restricted Burnside problem
\cite{Z}. Then $\left|  G:G_{1}\right|  \leq\mu$, and it follows that $G_{1}$
can be generated by $d^{\prime}=d\mu$ elements.

Since $G_{1}$ is generated by $\mathfrak{X}$, the argument of \S 2, Step 5 now gives

\medskip

\textbf{Step 5}$_{C}$. $G_{1}\subseteq\mathfrak{X}^{\ast n_{3}}[H_{1}%
,G_{1}]H_{1}^{q}$ \emph{where }$n_{3}$ \emph{depends only on }$d^{\prime}%
$\emph{ and }$q$.

\medskip

The next step depends on the following simple observation, where $\sigma(q)$
will denote the number of distinct prime divisors of $q$.

\begin{lemma}
If $H=\left\langle X\right\rangle $ is an $r$-generator abelian group then%
\[
H=\left\langle y_{1}^{q},\ldots,y_{r}^{q},\,x_{1},\ldots,x_{r\sigma
(q)}\right\rangle
\]
for some $y_{1},\ldots,y_{r}\in H$ and some $x_{1},\ldots,x_{r\sigma(q)}\in X$.
\end{lemma}

\begin{proof}
Let $P$ be a Sylow $p$-subgroup of $H$. Write $\pi:H\rightarrow P$ for the
projection. $P$ is an $r$-generator $p$-group generated by $\pi(X)$ so
$P=\left\langle \pi(X_{p})\right\rangle $ for some subset $X_{p}$ of $X$ of
size $r$ (because $P/\mathrm{Frat}(P)$ is an $r$-dimensional $\mathbb{F}_{p}$-vector
space). Thus if $p_{1},\ldots,p_{\sigma}$ are the primes dividing $q$ and
$P_{1},\ldots,P_{\sigma}$ the corresponding Sylow subgroups, then the subgroup
$R=\left\langle X_{p_{1}}\cup\ldots\cup X_{p_{\sigma}}\right\rangle $ projects
onto each $P_{i}$. It follows that $\left|  H:R\right|  $ is coprime to $q$
and hence that $H=QR$ where $Q$ is a direct factor of $H$ of order coprime to
$q$. Thus $Q$ is an $r$-generator group and each element of $Q$ is a $q$th
power, so $Q=\left\langle y_{1}^{q},\ldots,y_{r}^{q}\right\rangle $ for some
$y_{1},\ldots,y_{r}$. The lemma follows.
\end{proof}

\bigskip

Applying this lemma to $G_{1}/G_{1}^{\prime}$, we deduce that%
\[
G_{1}=G_{1}^{\prime}\left\langle h_{1},\ldots,h_{d^{\prime\prime}%
}\right\rangle
\]
where each $h_{i}\in\mathfrak{X}$ and $d^{\prime\prime}=d^{\prime}%
(1+\sigma(q))$. Now Lemma \ref{nilp-comm} gives%
\[
\lbrack H_{1},G_{1}]=\prod_{i=1}^{d^{\prime\prime}}[H_{1},h_{i}]\cdot
H_{2}\subseteq\mathfrak{X}^{\ast2d^{\prime\prime}}H_{2}.
\]
As $H_{1}^{q}\subseteq\lbrack H_{1},G_{1}]H_{1}^{\{q\}}$ we have established

\medskip

\textbf{Step 4}$_{C}$. $[H_{1},G_{1}]H_{1}^{q}\subseteq\mathfrak{X}%
^{\ast(2d^{\prime\prime}+1)}H_{2}.$

\medskip

Putting the last two steps together gives $G_{1}=\mathfrak{X}^{\ast n_{4}%
}H_{2}$ where $n_{4}$ depends only on $d$ and $q$. As $G_{1}$ is generated by
$d^{\prime}$ elements, it follows that there exist $g_{1},\ldots,g_{r}%
\in\mathfrak{X}$, where $r=n_{4}d^{\prime},$ such that%
\[
G_{1}=H_{2}\left\langle g_{1},\ldots,g_{r}\right\rangle .
\]
Since $H_{2}/H_{3}$ is an acceptable normal subgroup of $G_{1}/H_{3}$, Key
Theorem (C) may be applied to give

\medskip

\textbf{Step 3}$_{C}$. $H_{2}\subseteq\mathfrak{X}^{\ast n_{5}}H_{3}$ where
$n_{5}=2rh_{3}(d^{\prime},c)$.

\medskip

\textbf{Step 2}$_{C}$. $H_{3}\subseteq\mathfrak{X}^{\ast n_{1}}H_{5}$ where
$n_{1}$ depends only on $q$. This is identical to Step 2 in \S 2.

\medskip

\textbf{Step 1}$_{C}$. $H_{5}\subseteq\mathfrak{X}^{\ast n_{2}}$ where $n_{2}$
depends only on $q,\,d$ and $c$. We proved in Step 1 of \S 2 that $H_{5}$ is
an acceptable normal subgroup of $G_{1}$. So the claim will follow by Key
Theorem (C) if we can show that $G_{1}=H_{5}\left\langle g_{1}^{\prime}%
,\ldots,g_{s}^{\prime}\right\rangle $ where each $g_{i}^{\prime}%
\in\mathfrak{X}$ and $s$ depends only on $q,\,d$ and $c$. But this follows
from the preceding four steps: for $G_{1}$ is generated by $d^{\prime}$
elements, each of which lies in $\mathfrak{X}^{\ast n_{6}}H_{5}$ where
$n_{6}=n_{1}+n_{5}+n_{4}$, so we may take $s=d^{\prime}n_{6}$.

\medskip

\textbf{Conclusion.} Altogether we obtain (\ref{Gq=Xh}) with $h=n_{6}+n_{2}.$

\section{Proof of the Key Theorem}

\subsection*{The general idea}

Before getting down to specifics, let us outline the general plan of attack.
The Key Theorem asserts that, under suitable hypotheses on the finite group
$G$ and its normal subgroup $H$, every element of $H$ is equal to a product of
a specific form. Thus what has to be established is the solvability of
equations like%
\begin{equation}
h=\Phi(u_{1},\ldots,u_{m}) \label{modeleq}%
\end{equation}
where%
\[
\Phi(u_{1},\ldots,u_{m})=U(g_{1},\ldots,g_{r},\,u_{1},\ldots,u_{m});
\]
here the `constant' $h$ is an arbitrary element of $H$, $U$ is a specific
group word, $g_{1},\ldots,g_{r}$ are some fixed parameters from $G$, and the
`unknowns' $u_{1},\ldots,u_{m}$ are to be found in $H$. The idea of the proof
is modelled on that of Hensel's Lemma: one shows that an approximate solution
of (\ref{modeleq}) can be successively refined to an exact solution.

What makes Hensel's Lemma work is a hypothesis that ensures the surjectivity
of a certain linear map: the relevant derivative must be non-singular modulo
$p$. This translates in a straightforward way to our context.

\medskip

\noindent\textbf{Definition} Let $\mathbf{v}\in H^{(m)}$. The mapping
$\Phi_{\mathbf{v}}^{\prime}:H^{(m)}\rightarrow H$ is defined by%
\[
\Phi(\mathbf{x}\cdot\mathbf{v})=\Phi_{\mathbf{v}}^{\prime}(\mathbf{x}%
)\cdot\Phi(\mathbf{v})\qquad(\mathbf{x}\in H^{(m)})
\]
where $\mathbf{x}\cdot\mathbf{v}$ denotes the $m$-tuple $(x_{1}v_{1}%
,\ldots,x_{m}v_{m})$.

\medskip

Suppose now that $K$ is a normal subgroup of $G$ contained in $H$, and that we
have found a solution of (\ref{modeleq}) modulo $K$; that is, we have
$\mathbf{v}\in H^{(m)}$ such that%
\[
h=\kappa\cdot\Phi(\mathbf{v})
\]
for some $\kappa\in K$. Then $\mathbf{u}=\mathbf{x}\cdot\mathbf{v}$ is a
solution of (\ref{modeleq}) if and only if%
\begin{equation}
\Phi_{\mathbf{v}}^{\prime}(\mathbf{x})=\kappa. \label{firstsoln}%
\end{equation}
Thus our `approximate solution' $\mathbf{v}$ can be lifted to an exact
solution provided the image of the map $\Phi_{\mathbf{v}}^{\prime}$ contains
$K$. Let us call $\mathbf{v}$ `liftable' in this case. To ensure that the
process can be iterated, however, we require that the `new' solution
$\mathbf{x}\cdot\mathbf{v}$ is again liftable in the appropriate sense. This
will be achieved by a `probabilistic' argument: we establish independently (a)
that a relatively large proportion of the elements $\mathbf{x}$ in a suitable
domain are solutions of (\ref{firstsoln}), and (b) that a relatively large
proportion of the $\mathbf{x}$ in the same domain have the property that
$\mathbf{x}\cdot\mathbf{v}$ is liftable. It will follow that at least some of
these elements $\mathbf{x}$ will have both properties.

Here is a final remark. All our main results about finite groups concern
functions that are uniformly bounded in terms of $d$, the number of
generators. Why is this the dominant parameter? There are two reasons. The
first is evident in the statement of the Key Theorem: each of the $d$
generators appears explicitly in the statement. The second, hidden in the
proof, is to do with the way the generators have to act on chief factors of
the group; it comes down to the following obvious but crucial observation:

\begin{lemma}
\label{fixedpts}Let $G=\left\langle g_{1},\ldots,g_{d}\right\rangle $ be a group.

\emph{(i)} If $G$ acts without fixed points on a set of size $n$ then at least
one of the $g_{i}$ moves at least $n/d$ points.

\emph{(ii)} If $G$ acts linearly on a vector space $V$ of dimension $n$, and
fixes only $0$, then at least one of the $g_{i}$ satisfies $\dim\mathrm{C}%
_{V}(g_{i})\leq(1-\frac{1}{d})n$.
\end{lemma}

\noindent(Here $\mathrm{C}_{V}(g)$ denotes the fixed-point set of $g$.)

\subsection*{Solvability of equations}

Let $G$ be a finite group. A normal subgroup $N$ of $G$ will be called
\emph{quasi-minimal} if $N=[N,G]>1$ and $N$ is minimal with this property. It
is easy to see that in this case, there is a uniquely determined normal
subgroup $Z=Z_{N}$ of $G$ maximal subject to $Z<N$; indeed, if $Z_{1}$ and
$Z_{2}$ were two distinct such subgroups then $N=Z_{1}Z_{2}$ would imply
$[N,G]=[N,Z_{1}][N,Z_{2}]=1$.

We write `QMN' for `quasi-minimal normal subgroup', and recall the definition
of `acceptable' from \S 3. The Frattini subgroup of $G$ is denoted
$\mathrm{Frat}(G)$.

\begin{lemma}
\label{qmn}Let $N$ be a QMN of $G$ and put $Z=Z_{N}$. Suppose that $N\leq H$
where $H$ is an acceptable normal subgroup of $G$. Then

\emph{(i)} $N/Z$ is a minimal normal subgroup of $G/Z$, $[Z,_{k}G]=1$ for some
$k$, and $[Z,N]\leq\lbrack Z,H]=1.$

\emph{(ii)} $Z\leq\mathrm{Frat}(G)$.

\emph{(iii)} If $N$ is not soluble then $N$ is quasi-semisimple with centre
$Z$ and $N/Z=S_{1}\times\cdots\times S_{n}$, where $n\geq3$ and $S_{1}%
,\ldots,S_{n}$ are isomorphic non-abelian simple groups.

\emph{(iv)} If $N$ is soluble then $N/Z$ is an elementary abelian $p$-group
for some prime $p$; also $N^{p}=1$ if $p$ is odd, $N^{2}=[N,N]$ and
$[N,N]^{2}=1$ if $p=2$.
\end{lemma}

\begin{proof}
(i) The first two statements are immediate from the definition. To show that
$[Z,H]=1$, write $Z_{i}=[Z,_{i}G]$ for $i\geq0$ (with $Z_{0}=Z$). Then
$Z_{k}=1$. Suppose we have $[Z_{i},H]=1$ for some $i$ with $k\geq i>0$. Since
$H=[H,G]$ the Three-Subgroup Lemma gives%
\[
\lbrack Z_{i-1},H]=[[Z_{i-1},H],G][Z_{i},H]=[[Z_{i-1},H],G]
\]
whence $[Z_{i-1},H]=1$ since $[Z_{i-1},H]<N$. It follows by reverse induction
that $[Z,H]=[Z_{0},H]=1$.

(ii) Suppose that $M$ is a maximal subgroup of $G$ and $M$ contains $Z_{i}$
but not $Z_{i-1},$ where $i>0$. Then $G=Z_{i-1}M$ and $H=Z_{i-1}(H\cap M)$ so
$[H,G]\leq M$, a contradiction since $Z_{i-1}\leq H=[H,G]$.

(iii) This follows from the well-known structure of minimal normal subgroups;
here $n\geq3$ because $N$ is contained in the acceptable subgroup $H$.

(iv) The first claim is standard. Since $[N,N]\leq\mathrm{Z}(N)$, the map
$x\mapsto x^{p}$ is a homomorphism of $G$-operator groups from $N$ into $Z$ if
$p$ is odd, and induces such a homomorphism from $N$ into $Z/[N,N]$ if $p=2$.
In each case the image of this homomorphism must be $1$ since $N=[N,G]$. The
final statement is easy.
\end{proof}

\bigskip

The solvability of equations like (\ref{firstsoln}) is assured by the
following results, which will be proved in later sections (the fourth one,
Proposition 11.1, is needed only for variant (B) of the Key Theorem). In each
case, $N$ denotes a QMN of $G$ and $Z=Z_{N}$.

For $\mathbf{x}=(x_{1},\ldots,x_{t}),\,\mathbf{y}=(y_{1},\ldots,y_{t})\in
G^{(t)}$ we will write%
\[
\lbrack\mathbf{x},\mathbf{y}]=\prod_{j=1}^{t}[x_{j},y_{j}].
\]
\medskip

\noindent\textbf{Proposition 7.1} \ \emph{Suppose that }$N$\emph{ is soluble
and that }$[Z,G]=1$\emph{. Put }$K=N$\emph{ if }$N$\emph{ is abelian,
}$K=N^{\prime}$\emph{ otherwise. For }$i=1,2,3$\emph{ define }$\phi
_{i}:N^{(m)}\rightarrow N$\emph{ by}%
\[
\phi_{i}(\mathbf{a})=[\mathbf{a},\mathbf{y}_{i}]
\]
\emph{where }$\mathbf{y}_{i}=(y_{i1},\ldots,y_{im})$ \emph{and the }$y_{ij}%
$\emph{ are elements of }$G$\emph{ such that }$\left\langle y_{i1}%
,\ldots,y_{im}\right\rangle K=G$\emph{ for each }$i$\emph{. Let }$\kappa\in
K$\emph{. Then there exist }$\kappa_{1},\,\kappa_{2},\,\kappa_{3}\in N$\emph{
such that }$\kappa_{1}\kappa_{2}\kappa_{3}=\kappa$\emph{ and, for each}
$i=1,2,3$,%
\[
\left|  \phi_{i}^{-1}(\kappa_{i})\right|  \geq\left|  N\right|  ^{m}|\left|
N/Z\right|  ^{-d-1}.
\]

\bigskip

The corresponding results for a non-soluble QMN involve certain constants:

\begin{description}
\item $D\geq1$ is the absolute constant specified in Theorem 1.9, and we set
$\overline{D}=4+2D$;

\item $C(q)$ and $M(q)$ are the constants specified in Theorem 1.10, and we
set\newline$z(q)=M(q)\overline{D}(q+\overline{D}).$\medskip
\end{description}

\noindent\textbf{Definition} Let $\varepsilon>0$ and $k\in\mathbb{N}$. Let
$\mathbf{y}=(y_{1},y_{2},\ldots,y_{m})\in G^{(m)}$.

\begin{description}
\item[(i)] The $m$-tuple $\mathbf{y}$ has the $(k,\varepsilon)$
\emph{fixed-point property} if in any transitive permutation action of $G$ on
a set of size $n\geq2$, at least $k$ of the elements $y_{i}$ move at least
$\varepsilon n$ points.

\item[(ii)] The $m$-tuple $\mathbf{y}$ has the $(k,\varepsilon)$
\emph{fixed-space property }if for every irreducible $\mathbb{F}_{p}G$-module
$V$ of dimension $n\geq2$, where $p$ is any prime, at least $k$ of the $y_{i}$
satisfy $\dim_{\mathbb{F}_{p}}\mathrm{C}_{V}(y_{i})\leq(1-\varepsilon)n$.
\end{description}

\bigskip

\noindent\textbf{Proposition 9.2} \ \emph{Suppose that }$N$\emph{ is
quasi-semisimple, and that }$N/Z$\emph{ is not simple. Define }$\phi
:N^{(m)}\rightarrow N$\emph{ by}%
\[
\phi(\mathbf{a})=[\mathbf{a},\mathbf{y}]
\]
\emph{where }$y_{1},\ldots,y_{m}$\emph{ are elements of }$G$\emph{ such that
}$\left\langle y_{1},\ldots,y_{m}\right\rangle N=G$\emph{. Suppose that
}$\mathbf{y}$\emph{ has the }$(k,\varepsilon)$ \emph{fixed-point property
where }$k\varepsilon\geq\overline{D}$. \emph{Then for each }$\kappa\in N$,
\[
\left|  \phi^{-1}(\kappa)\right|  \geq\left|  N\right|  ^{m}|\left|
N/Z\right|  ^{-4D}.
\]
\bigskip

\noindent\textbf{Proposition 10.1} \ \emph{Let }$q\in\mathbb{N}$.
\emph{Suppose that }$N$\emph{ is quasi-semisimple, and that its non-abelian
composition factors }$S$\emph{ satisfy} $\left|  S\right|  >C(q)$\emph{. Let
}$u_{1},\ldots,u_{m}$ $\in G$ \emph{where} $m\geq z(q)$. \emph{Then the
mapping }$\psi:N^{(m)}\rightarrow N$\emph{ defined by}%
\[
\prod_{j=1}^{m}(x_{j}u_{j})^{q}=\psi(\mathbf{x})\prod_{j=1}^{m}u_{j}^{q}%
\]
\emph{is surjective.}

\medskip

\noindent\textbf{Proposition 11.1} \ \emph{Suppose that }$N$\emph{ is
quasi-semisimple, and let } $\alpha_{1},\,\beta_{1},\ldots,$ \linebreak%
$\alpha_{D},\beta_{D}\ $\emph{be} $2D$ \emph{arbitrary automorphisms of }%
$N$\emph{. Then the mapping }$\theta:N^{(2D)}\rightarrow N$ \emph{defined by}
\[
\theta(\mathbf{a},\mathbf{b})=\prod_{j=1}^{D}T_{\alpha_{j},\beta_{j}}%
(a_{j},b_{j})
\]
\emph{is surjective.}

\subsection*{Lifting generators}

The other half of our probabilistic argument rests on the following
proposition, which will be established in \S 5. For a simple group $S$ we
define $\mu(S)$ to be the supremum of the numbers $\mu$ such that%
\[
\left|  S:M\right|  \geq\left|  S\right|  ^{\mu}%
\]
for every maximal subgroup $M$ of $S$, and for any group $N$ define
\[
\mu^{\prime}(N)=\min\left\{  \frac{1}{2},\,\mu(S)\mid S\text{ a non-abelian
composition factor of }N\right\}  .
\]
For later use, we also define%
\[
\mu(q)=\min\left\{  \frac{1}{2},\,\mu(S)\mid S\text{ simple, }\left|
S\right|  \leq C(q)\right\}  .
\]

\medskip

\noindent\textbf{Proposition 5.1} \emph{Let }$G$ \emph{be a }$d$%
\emph{-generator group and }$N$\emph{ an acceptable QMN of }$G$.\emph{ Suppose
that} $G=\left\langle y_{1},\ldots,y_{m}\right\rangle N$. \emph{Put }$Z=Z_{N}$
\emph{and let}%
\[
\mathcal{N}(\mathbf{y})=\left\{  \mathbf{a}\in N^{(m)}\mid\left\langle
y_{1}^{a_{1}},\ldots,y_{m}^{a_{m}}\right\rangle \neq G\right\}  .
\]
\emph{Let }$\varepsilon\in(0,\frac{1}{2}]$.

(i) \emph{Suppose that }$N$ \emph{is soluble and that} $\mathbf{y}$ \emph{has
the} $(k,\varepsilon)$ \emph{fixed-space property. Then}%
\[
\left|  \mathcal{N}(\mathbf{y})\right|  \leq\left|  N\right|  ^{m}\left|
N/Z\right|  ^{d-k\varepsilon}.
\]

(ii) \emph{There exists an absolute constant}$\,C_{0}$\emph{ such that if }$N$
\emph{is quasi-semisimple and }$\mathbf{y}$ \emph{has the} $(k,\varepsilon)$
\emph{fixed-point property, where} $k\varepsilon\geq\max\{2d+4,\,2C_{0}%
+2\}$\emph{, then}
\begin{align*}
\left\vert \mathcal{N}(\mathbf{y})\right\vert  &  <\left\vert N\right\vert
^{m}\text{ \ \ and}\\
\left\vert \mathcal{N}(\mathbf{y})\right\vert  &  \leq\left\vert N\right\vert
^{m}\left\vert N/Z\right\vert ^{1-s}%
\end{align*}
\emph{where}%
\[
s=\min\{\mu^{\prime}(N)(k\varepsilon/2-d-1),\,\mu^{\prime}(N)(k\varepsilon
/2-C_{0})\}.
\]

\subsection*{The proof}

Now we can prove the Key Theorem, assuming the results stated above. We will
need to know the following `derivative', obtained by direct calculation:

\begin{lemma}
\label{Xi'}Let $\mathbf{g}\in G^{(m)}$. Define $\Xi:G^{(m)}\rightarrow G$ by
$\Xi(\mathbf{v})=[\mathbf{v},\mathbf{g}]$. Then%
\[
\Xi_{\mathbf{v}}^{\prime}(\mathbf{x})=\prod_{j=1}^{m}[x_{j},g_{j}]^{\tau
_{j}(\mathbf{g},\mathbf{v})}%
\]
where%
\[
\tau_{j}(\mathbf{g},\mathbf{v})=v_{j}[g_{j-1},v_{j-1}]\ldots\lbrack
g_{1},v_{1}].
\]

\end{lemma}

Now define%
\[
k(d,q)=1+\left\lceil d\cdot\max\left\{  \,\frac{8D+2}{\mu(q)}+2d+2,\,\frac
{8D+2}{\mu(q)}+2C_{0}\right\}  \right\rceil
\]
(where $\left\lceil x\right\rceil $ denotes the least integer $\geq x$), and
let $z(q)$ be as defined above. The first claim in the next proposition gives
the Key Theorem, on putting%
\[
h_{1}(d,q)=3k(d,q).
\]

\begin{proposition}
\label{KTA}Let $G=\left\langle g_{1},\ldots,g_{d}\right\rangle $ and let $H$
be an acceptable normal subgroup of $G$. Let $m=d\cdot k(d,q)$ and define
$\mathbf{g}=(g_{1},\ldots,g_{m})$ by setting%
\[
g_{td+i}=g_{i}\qquad(0\leq t<k(d,q)).
\]
Then for each $h\in H$ there exist $\mathbf{v}(1),\,\mathbf{v}(2),\,\mathbf{v}%
(3)\in H^{(m)}$ and $\mathbf{u}\in H^{(z(q))}$ such that%
\begin{equation}
h=\prod_{i=1}^{3}[\mathbf{v}(i),\mathbf{g}]\cdot\prod_{l=1}^{z(q)}u_{l}^{q}
\label{h=Phi}%
\end{equation}
and%
\begin{equation}
\left\langle g_{1}^{\tau_{1}(\mathbf{g},\mathbf{v}(i))},\ldots,g_{m}^{\tau
_{m}(\mathbf{g},\mathbf{v}(i))}\right\rangle =G\text{ \ \ for \ }i=1,\,2,\,3.
\label{gen(i)}%
\end{equation}

\end{proposition}

The second claim, (\ref{gen(i)}), is required for the inductive proof. In
terms of the heuristic discussion above, it ensures that our solution
$(\mathbf{v}(1),\,\mathbf{v}(2),\,\mathbf{v}(3),\,\mathbf{u})$ is again
`liftable': in the guise of (\ref{genmodK(2)}) or (\ref{genmodK(3)}), it is
used directly in `Case 1', below, and in other cases enables us to quote some
of the above-stated propositions, whose hypotheses stipulate that a certain
set of elements should generate an appropriate quotient of $G$.

Let us recall that $H$ is \emph{acceptable} in $G$ if (i) $H=[H,G]$ and (ii)
no normal section of $G$ inside $H$ takes the form $S$ or $S\times S$ for a
non-abelian simple group $S$. It is clear that $H/K$ is then acceptable in
$G/K$ whenever $H\geq K$ and $K\vartriangleleft G$; we shall use this without
special mention.\bigskip

\textbf{Proof.} We will write $k=k(d,q)$ and $z=z(q)$. The result is trivial
if $H=1$; we suppose that $H>1$ and argue by induction on $\left|  H\right|
$. Since $H=[H,G]$ it follows that $H$ contains a QMN $N$ of $G$. It also
follows that $d\geq2$. Put $Z=Z_{N}$ and define a normal subgroup $K>1$ of $G$
as follows:%
\begin{equation}
K=\left\{
\begin{array}
[c]{ccc}%
\lbrack Z,G] &  & \text{if }[Z,G]>1\\
&  & \\
N &  & \text{if }[Z,G]=1\text{ and }[N,N]=1\\
&  & \\
\lbrack N,N] &  & \text{if }[Z,G]=1\text{ and }[N,N]>1
\end{array}
.\right.  \label{defn K}%
\end{equation}

Write the equation (\ref{h=Phi}) as%
\[
h=\Phi(\mathbf{v},\mathbf{u})=\Xi(\mathbf{v}(1))\cdot\Xi(\mathbf{v}%
(2))\cdot\Xi(\mathbf{v}(3))\cdot\Psi(\mathbf{u}).
\]
Inductively, we may assume that there exist $\kappa\in K$, $\mathbf{v}(i)\in
H^{(m)}$ and $\mathbf{u}\in H^{(z)}$ such that%
\[
h=\kappa\Phi(\mathbf{v},\mathbf{u})
\]
and, for \ $i=1,\,2,\,3$,%
\begin{equation}
\left\langle g_{1}^{\tau_{1}(\mathbf{v}(i))},\ldots,g_{m}^{\tau_{m}%
(\mathbf{v}(i))}\right\rangle K=G, \label{genmodK(1)}%
\end{equation}
where for brevity we write $\tau_{j}(\mathbf{x})=\tau_{j}(\mathbf{g}%
,\mathbf{x})$.

The aim is to show that there exist $\mathbf{a}(i)\in N^{(m)}$ and
$\mathbf{b}\in N^{(z)}$ such that (\ref{h=Phi}) and (\ref{gen(i)}) hold with
$\mathbf{a}(i)\cdot\mathbf{v}(i)$ replacing $\mathbf{v}(i)$ and $\mathbf{b}%
\cdot\mathbf{u}$ replacing $\mathbf{u}$. The first requirement is equivalent
to%
\begin{align}
\kappa &  =\Phi_{(\mathbf{v},\mathbf{u})}^{\prime}(\mathbf{a},\mathbf{b}%
)\nonumber\\
&  =\Xi_{\mathbf{v}(1)}^{\prime}(\mathbf{a}(1))^{\xi_{1}}\cdot\Xi
_{\mathbf{v}(2)}^{\prime}(\mathbf{a}(2))^{\xi_{2}}\cdot\Xi_{\mathbf{v}%
(3)}^{\prime}(\mathbf{a}(3))^{\xi_{3}}\cdot\Psi_{\mathbf{u}}^{\prime
}(\mathbf{b})^{\xi_{4}} \label{kappaPhi'}%
\end{align}
where $\xi_{1}=1$ and%
\[
\xi_{i}=\left(  \Xi(\mathbf{v}(1))\ldots\Xi(\mathbf{v}(i-1))\right)
^{-1}\qquad(i=2,\,3,\,4).
\]

It is convenient to reformulate the second requirement. Write%
\[
\overline{a}(i)_{j}=a(i)_{j}^{v(i)_{j}g_{j}\ldots g_{m}},\,\,\overline{g}%
_{ij}=g_{j}^{v(i)_{j}g_{j}\ldots g_{m}}.
\]

\begin{lemma}
\label{newgens}Let $\mathbf{a}(i)\in N^{(m)}$ for $i=1,\,2,\,3$. The following
are equivalent, for each $i$:%
\begin{align}
G  &  =\left\langle g_{1}^{\tau_{1}(\mathbf{a}(i)\cdot\mathbf{v}(i))}%
,\ldots,g_{m}^{\tau_{m}(\mathbf{a}(i)\cdot\mathbf{v}(i))}\right\rangle
,\label{<.gamma.>}\\
G  &  =Z\left\langle \overline{g}_{i1}^{\overline{a}(i)_{1}},\ldots
,\overline{g}_{im}^{\overline{a}(i)_{m}}\right\rangle . \label{<.barg..>}%
\end{align}

\end{lemma}

\begin{proof}
We claim that for any $m$-tuple $\mathbf{v}$,%
\begin{equation}
\left\langle g_{1}^{\tau_{1}(\mathbf{v})},\ldots,g_{m}^{\tau_{m}(\mathbf{v}%
)}\right\rangle ^{g_{1}\ldots g_{m}}=\left\langle g_{1}^{v_{1}g_{1}\ldots
g_{m}},\ldots,g_{m}^{v_{m}g_{m}}\right\rangle . \label{twogensets}%
\end{equation}
To see this, put $z_{1}=1$ and for $k>1$ set
\[
z_{k}=g_{k-1}^{v_{k-1}g_{k-2}^{-1}\ldots g_{1}^{-1}}\cdot z_{k-1}.
\]
Arguing by induction on $k$ we find that $z_{k+1}=z_{k}g_{k}^{\tau
_{k}(\mathbf{v})}$ for each $k$; this implies that%
\[
\left\langle g_{1}^{\tau_{1}(\mathbf{v})},\ldots,g_{m}^{\tau_{m}(\mathbf{v}%
)}\right\rangle =\left\langle z_{2},\ldots,z_{m+1}\right\rangle =\left\langle
g_{1}^{v_{1}},\ldots,g_{m}^{v_{m}g_{m-1}^{-1}\ldots g_{1}^{-1}}\right\rangle
\]
which is equivalent to (\ref{twogensets}).

The lemma follows on taking $\mathbf{v}=\mathbf{a}(i)\cdot\mathbf{v}(i)$, and
noting that $\overline{g}_{ij}^{\overline{a}(i)_{j}}=g_{j}^{a(i)_{j}%
v(i)_{j}g_{j}\ldots g_{m}}$ and $Z\leq\mathrm{Frat}(G)$.
\end{proof}

\medskip

Taking each $a(i)_{j}=1$ and replacing $G$ by $G/K$, we deduce that
(\ref{genmodK(1)}) implies%
\begin{equation}
G=\left\langle \overline{g}_{i1},\ldots,\overline{g}_{im}\right\rangle
K\qquad(i=1,\,2,\,3). \label{genmodK(2)}%
\end{equation}
Now write%
\[
\widetilde{g}_{ij}=g_{j}^{\tau_{j}(\mathbf{v}(i))\xi_{i}},\quad\widetilde
{a}(i)_{j}=a(i)_{j}^{\tau_{j}(\mathbf{v}(i))\xi_{i}}.
\]
Then (\ref{genmodK(1)}) is also (evidently) equivalent to%
\begin{equation}
G=\left\langle \widetilde{g}_{i1},\ldots,\widetilde{g}_{im}\right\rangle
K\qquad(i=1,\,2,\,3); \label{genmodK(3)}%
\end{equation}
and Lemma \ref{Xi'} shows that%
\begin{equation}
\Xi_{\mathbf{v}(i)}^{\prime}(\mathbf{a}(i))^{\xi_{i}}=[\widetilde{\mathbf{a}%
}(i),\widetilde{\mathbf{g}}_{i}]\qquad(i=1,\,2,\,3). \label{Xi'=[a,g]}%
\end{equation}

Thus it suffices to find $\mathbf{a}(i)$ and $\mathbf{b}$ (with entries in
$N$) such that%
\begin{equation}
\kappa=[\widetilde{\mathbf{a}}(1),\widetilde{\mathbf{g}}_{1}][\widetilde
{\mathbf{a}}(2),\widetilde{\mathbf{g}}_{2}][\widetilde{\mathbf{a}%
}(3),\widetilde{\mathbf{g}}_{3}]\Psi_{\mathbf{u}}^{\prime}(\mathbf{b}%
)^{\xi_{4}} \label{finaleq}%
\end{equation}
and such that (\ref{<.barg..>}) holds. To this end we separate several cases.

\medskip

\textbf{Case 1:} where $[Z,G]=K>1$. We think of $Z$ as a $G$-module, with $K$
acting trivially, and write it additively. From (\ref{genmodK(3)}) we have%
\[
K=Z(G-1)=\sum_{j=1}^{m}Z(\widetilde{g}_{1j}-1)=\left\{  [\mathbf{z}%
,\widetilde{\mathbf{g}}_{1}]\mid\mathbf{z}\in Z^{(m)}\right\}  .
\]
Thus there exists $\mathbf{a}(1)\in Z^{(m)}$ with $[\mathbf{a}(1),\mathbf{g}%
_{1}]=\kappa$, and we may satisfy (\ref{finaleq}) by setting $a(2)_{j}%
=a(3)_{j}=b_{j}=1$ for all $j\,$; note that $\widetilde{\mathbf{a}%
}(1)=\mathbf{a}(1)$ here since $[Z,H]=1$. As each $\overline{a}(i)_{j}$ is in
$Z$ and $K\leq Z$, in this case (\ref{<.barg..>}) follows at once from
(\ref{genmodK(2)}). $\blacksquare$

\medskip

Assume henceforth that $[Z,G]=1$. For $\kappa\in N$ and $1\leq i\leq3$ put%
\[
\mathfrak{X}_{i}(\kappa)=\left\{  \mathbf{a}(i)\in N^{(m)}\mid\lbrack
\widetilde{\mathbf{a}}(i),\widetilde{\mathbf{g}}_{i}]=\kappa\right\}  ,
\]
and let%
\[
\mathfrak{Y}_{i}=\left\{  \mathbf{a}(i)\in N^{(m)}\mid\left\langle
\overline{g}_{i1}^{\overline{a}(i)_{1}},\ldots,\overline{g}_{im}^{\overline
{a}(i)_{m}}\right\rangle Z=G\right\}  .
\]
We shall repeatedly use the following

\medskip

\noindent\textbf{Key Observation}: For each $i=1,\,2,\,3$, \emph{the}
$m$\emph{-tuple} $\overline{\mathbf{g}}_{i}$ \emph{has the} $(k,\frac{1}{d})$
\emph{fixed-space property and the} $(k,\frac{1}{d})$ \emph{fixed-point
property}.

\medskip

\noindent Indeed, since $G=\left\langle g_{1},\ldots,g_{d}\right\rangle $,
Lemma \ref{fixedpts} shows that the $d$-tuple $(g_{1},\ldots,g_{d})$ has the
$(1,\frac{1}{d})$ fixed-space property and the $(1,\frac{1}{d})$ fixed-point
property \ The claim follows because each of the generators $g_{l}$ ($1\leq
l\leq d$) is conjugate to at least $k$ of the elements $\overline{g}_{ij}$
($1\leq j\leq m$).

\medskip

\textbf{Case 2:} where $N$ is soluble, and $K=N$ if $N$ is abelian,
$K=N^{\prime}$ if not. Define $\,\phi_{i}:N^{(m)}\rightarrow N$ by%
\[
\phi_{i}(\mathbf{x})=[\mathbf{x},\widetilde{\mathbf{g}}_{i}].
\]
In view of (\ref{genmodK(3)}), we may take $y_{ij}=\widetilde{g}_{ij}$ in
Proposition 7.1 and infer that there exist $\kappa_{1},\,\kappa_{2}%
,\,\kappa_{3}\in N$ with $\kappa_{1}\kappa_{2}\kappa_{3}=\kappa$ such that%
\[
\left|  \mathfrak{X}_{i}(\kappa_{i})\right|  =\left|  \phi_{i}^{-1}(\kappa
_{i})\right|  \geq\left|  N\right|  ^{m}|\left|  N/Z\right|  ^{-d-1}%
\]
for $i=1,\,2,\,3$ (the first equality holds because $\mathbf{a}(i)\mapsto
\widetilde{\mathbf{a}}(i)$ is a bijection on $N^{(m)}$).

Let $i\in\{1,\,2,\,3\}$. With the \emph{Key Observation} and (\ref{genmodK(2)}%
), Proposition 5.1(i) shows that the number of elements $\mathbf{x}\in
N^{(m)}$ for which%
\[
\left\langle \overline{g}_{i1}^{x_{1}},\ldots,\overline{g}_{im}^{x_{m}%
}\right\rangle \neq G
\]
is at most $\left|  N\right|  ^{m}\left|  N/Z\right|  ^{d-k/d}$. Since
$\mathbf{a}(i)\mapsto\overline{\mathbf{a}}(i)$ is a bijection on $N^{(m)}$
this gives%
\[
\left|  N^{(m)}\setminus\mathfrak{Y}_{i}\right|  \leq\left|  N\right|
^{m}\left|  N/Z\right|  ^{d-k/d}.
\]
As $k>d(2d+1)$, it follows that $\left|  \mathfrak{X}_{i}(\kappa_{i})\right|
>\left|  N^{(m)}\setminus\mathfrak{Y}_{i}\right|  $.

Thus we may choose $\mathbf{a}(i)\in\mathfrak{X}_{i}(\kappa_{i})\cap
\mathfrak{Y}_{i}$, for $i=1,\,2,\,3$. Then (\ref{<.barg..>}) holds and
(\ref{finaleq}) is satisfied with $b_{l}=1$ for all $l$. $\blacksquare$

\medskip

\textbf{Case 3:} where $N=K$ is quasi-semisimple and $\left|  S\right|  \leq
C(q)$; here $S$ denotes the (unique) non-abelian composition factor of $N$.

Put $\kappa_{1}=\kappa$, $\kappa_{2}=\kappa_{3}=1$. Using Proposition 9.2 in
place of Proposition 7.1, we see just as in Case 2 that for $i=1,\,2,\,3$,%
\[
\left|  \mathfrak{X}_{i}(\kappa_{i})\right|  \geq\left|  N\right|
^{m}|\left|  N/Z\right|  ^{-4D};
\]
note that $k/d>\overline{D}$ because $D\geq1>\mu(q)$.

Now $k/d>\max\{2d+4,\,2C_{0}+2\}$, so Proposition 5.1(ii), with the \emph{Key
Observation} and (\ref{genmodK(2)}), shows that%
\[
\left\vert N^{(m)}\setminus\mathfrak{Y}_{i}\right\vert \leq\left\vert
N\right\vert ^{m}\left\vert N/Z\right\vert ^{1-s}%
\]
for each $i$, where%
\begin{align*}
s  &  =\min\{\mu(q)(k/2d-d-1),\,\mu(q)(k/2d-C_{0})\}\\
&  >4D+1.
\end{align*}
We conclude as in the preceding case that (\ref{finaleq}) and (\ref{<.barg..>}%
) can be simultaneously satisfied by a suitable choice of $\mathbf{a}%
(1),\,\mathbf{a}(2),\,\mathbf{a}(3)\in N^{(m)}$, taking each $b_{l}=1$.
$\blacksquare$

\medskip

\textbf{Case 4:} where $N=K$ is quasi-semisimple and $\left|  S\right|
>C(q)$. Applying Proposition 5.1(ii) again we infer that each of the sets
$\mathfrak{Y}_{i}$ is non-empty. Choose $\mathbf{a}(i)\in\mathfrak{Y}_{i}$ for
$i=1,\,2,\,3$. Then (\ref{<.barg..>}) holds.

Now Proposition 10.1 shows that the mapping $\Psi_{\mathbf{u}}^{\prime
}:N^{(z)}\rightarrow N$ is surjective. Hence there exists $\mathbf{b}\in
N^{(z)}$ such that%
\[
\Psi_{\mathbf{u}}^{\prime}(\mathbf{b})=\left(  \left(  [\widetilde{\mathbf{a}%
}(1),\widetilde{\mathbf{g}}_{1}][\widetilde{\mathbf{a}}(2),\widetilde
{\mathbf{g}}_{2}][\widetilde{\mathbf{a}}(3),\widetilde{\mathbf{g}}%
_{3}]\right)  ^{-1}\kappa\right)  ^{\xi_{4}^{-1}}.
\]
Then (\ref{finaleq}) is satisfied, and the proof is complete. $\blacksquare$

\medskip

\noindent\textbf{Remark.} It may be worth observing that in Case 3, the only
role played by the upper bound on $|S|$ is to provide the lower bound $\mu(q)$
for $\mu(S)$. In fact such a lower bound will obtain if we allow $S$ to range,
additionally, over groups of Lie type with bounded Lie ranks (but over finite
fields of arbitrary size); this follows from Lemma 4.8, below, for example. We
may therefore, if we prefer, restrict Case 4 to where $S$ is either
alternating of large degree or of Lie type with large Lie rank. This means
that for the Key Theorem, only the special case of Proposition 10.1 relating
to such simple groups $S$ is actually needed. This in turn depends only on the
corresponding special case of Theorem 1.10; thus (for present purposes) one
can do without the fair-sized chunk of Part II devoted to the proof of Theorem
1.10 for groups of Lie type with small Lie rank over large fields. (However,
for groups of this type we shall still need the rather easier special case of
Theorem 1.10 where $q=1$, in order to deduce Theorem 1.9.)

\subsection*{Variants (B) and (C)}

Define%
\[
k(d)=1+\left\lceil d\cdot\max\{2d+4,\,2C_{0}+2\}\right\rceil .
\]
Now modify the statement of Proposition \ref{KTA} as follows: replace $k(d,q)$
by $k(d)$, replace $z(q)$ by $2D$, and replace the formula (\ref{h=Phi}) by%
\begin{equation}
h=\prod_{i=1}^{3}[\mathbf{v}(i),\mathbf{g}]\cdot\prod_{l=1}^{D}[u_{l}%
,u_{l+D}].\label{h=Phi(B)}%
\end{equation}
This gives Key Theorem (B) if we set $h_{2}(d)=3k(d)$.

For the proof of the modified proposition, \ we set $\Psi(\mathbf{x}%
,\mathbf{y})=\prod_{j=1}^{D}[x_{j},y_{j}]$ and use

\begin{lemma}
\label{Psi'(B)}%
\[
\Psi_{(\mathbf{u},\mathbf{w})}^{\prime}(\mathbf{x},\mathbf{y})=\prod_{l=1}%
^{D}T_{\alpha_{l},\beta_{l}}(x_{l}^{\sigma_{l}},y_{l}^{\rho_{l}})
\]
where $\alpha_{l},\,\beta_{l},\,\sigma_{l}$ and $\rho_{l}$ are given by
certain group words in $u_{1},w_{1},\ldots,u_{D},w_{D}$.
\end{lemma}

This is verified by direct calculation. We now argue exactly as before, with
the following changes: omit Case 3 altogether; and in Case 4, remove the
restriction on $\left|  S\right|  $ and use Proposition 11.1 in place of
Proposition 10.1. With Lemma \ref{Psi'(B)}, this shows that the relevant
mapping $\Psi_{\mathbf{u}}^{\prime}:N^{(2D)}\rightarrow N$ is surjective.
$\blacksquare$

\bigskip

The modifications required for Key Theorem (C) are a little more drastic, so
let us state the appropriate variant of Proposition \ref{KTA}. Define%
\[
k^{\prime}(d,c)=1+\left\lceil d\cdot\max\left\{  \,\frac{8D+2}{\varepsilon
(c)}+2d+2,\,\frac{8D+2}{\varepsilon(c)}+2C_{0}\right\}  \right\rceil
\]
where $\varepsilon(c)$ is the constant appearing in Lemma \ref{GSSh} below.

\begin{proposition}
\label{KTC}Let $G$ be a $d$-generator group with $\alpha(G)=c$ and let $H$ be
an acceptable normal subgroup of $G$. Suppose that $G=H\left\langle
g_{1},\ldots,g_{r}\right\rangle $. Put $m=r\cdot k^{\prime}(d,c)$ and define
$\mathbf{g}=(g_{1},\ldots,g_{m})$ by setting%
\[
g_{tr+i}=g_{i}\qquad(0\leq t<k^{\prime}(d,c)).
\]
Then for each $h\in H$ there exist $\mathbf{v}(1),\,\mathbf{v}(2),\,\mathbf{v}%
(3)\in H^{(m)}$ such that%
\[
h=\prod_{i=1}^{3}[\mathbf{v}(i),\mathbf{g}]
\]
and%
\[
\left\langle g_{1}^{\tau_{1}(\mathbf{g},\mathbf{v}(i))},\ldots,g_{m}^{\tau
_{m}(\mathbf{g},\mathbf{v}(i))}\right\rangle =G\text{ \ \ for \ }i=1,\,2,\,3.
\]

\end{proposition}

Key Theorem (C) then follows on setting $h_{3}(d,c)=3k^{\prime}(d,c)$. For the
proof, we may no longer appeal to Lemma \ref{fixedpts}; instead we rely on

\begin{lemma}
\label{GSSh}There exists $\varepsilon=\varepsilon(c)\in(0,\frac{1}{2}]$,
depending only on $c=\alpha(G)$, such that the following hold.

\emph{(i) }If $G$ acts as a primitive permutation group on a set $\Omega$ of
size $\geq2,$ with kernel $G_{\Omega}$, then $\left|  \Omega\right|
\geq\left|  G:G_{\Omega}\right|  ^{\varepsilon}$.

\emph{(ii) }For each transitive $G$-set $\Omega$ of size $\geq2,$ there is a
proper normal subgroup $G_{0}(\Omega)$ of $G$ such that for each $x\in
G\setminus G_{0}(\Omega)$,%
\[
\left|  \mathrm{fix}_{\Omega}(x)\right|  \leq(1-\varepsilon)\left|
\Omega\right|
\]
(where $\mathrm{fix}_{\Omega}(x)$ denotes the set of fixed points of $x$ in
$\Omega$).

\emph{(iii)} For each simple $\mathbb{F}_{p}G$-module $V$ there is a proper
normal subgroup $G_{0}(V)$ of $G$ such that for each $x\in G\setminus
G_{0}(V)$,%
\[
\dim\mathrm{C}_{V}(x)\leq(1-\varepsilon)\dim V\text{.}%
\]

\end{lemma}

\begin{proof}
Gluck, Seress and Shalev prove in \cite{GSS}, Theorem 1.2 that every primitive
$G$-set $\Omega$ contains a \emph{base} $B$ of size at most $\gamma=\gamma
(c)$, a number depending only on $c$ (to say that $B$ is a base means that the
pointwise stabilizer of $B$ is equal to $G_{\Omega}$). We may suppose that
$\gamma\geq2$. This gives (i) with $\varepsilon=\gamma^{-1}$, since the action
of each element of $G$ is determined by where it moves each element of $B$.
(In fact (i) is a celebrated result of Babai, Cameron and P\'{a}lfy \cite{BCP}.)

It also implies (ii) for the case of a primitive action. To see this, let
$x\in G\setminus G_{\Omega},$ let $\omega\in\Omega$ and put $X=\{y\in
G\mid\omega x^{y}\neq\omega\}$. Then%
\[
Xg_{1}\cup\ldots\cup Xg_{\gamma}=G
\]
where $B=\{\omega g_{1},\ldots,\omega g_{\gamma}\}$ so $\left|  X\right|
\geq\gamma^{-1}\left|  G\right|  $. Therefore $\Omega\setminus\mathrm{fix}%
_{\Omega}(x)=\{\omega y^{-1}\mid y\in X\}$ has cardinality at least
$\gamma^{-1}\left|  G\right|  /\left|  G_{\omega}\right|  =\gamma^{-1}\left|
\Omega\right|  $. In this case (ii) follows with $\varepsilon=\gamma^{-1}$ and
$G_{0}(\Omega)=G_{\Omega}$.

The general case of (ii) follows on taking $G_{0}(\Omega)$ to be the kernel of
the induced action on a minimal system of imprimitivity.

Statement (iii) for a primitive $\mathbb{F}_{p}G$-module $V$ is Theorem 5.3 of
\cite{GSS}, with $G_{0}(V)=\mathrm{C}_{G}(V)$. When $V$ is imprimitive, take
$G_{0}(V)$ to be the kernel of the permutation action of $G$ on a minimal
system of imprimitivity in $V,$ and apply (ii). (A better bound for
$\gamma(c)$ is given in \cite{LS1}, Theorem 1.4.)
\end{proof}

\bigskip

The proof now proceeds as in the preceding subsection, simply omitting the
function $\Psi$. The \emph{Key Observation} is replaced by

\medskip

\noindent\textbf{Key Observation }(C). Let $k=k^{\prime}(d,c)$. \emph{For
each} $i=1,\,2,\,3,$ \emph{the image in }$(G/K)^{(m)}$ \emph{of the}
$m$\emph{-tuple} $\overline{\mathbf{g}}_{i}$ \emph{has the} $(k,\varepsilon)$
\emph{fixed-space property and the} $(k,\varepsilon)$ \emph{fixed-point
property}.

\medskip

\noindent To see this, recall (\ref{genmodK(2)}), which asserts that the
$\overline{g}_{ij}K$ \ ($j=1,\ldots,m)$ generate $G/K$. Lemma \ref{GSSh}(ii)
then implies that for any transitive $G/K$-set of size $n\geq2$, at least one
of the elements $\overline{g}_{ij}K$ must move at least $\varepsilon n$
points. Since each $\overline{g}_{ij}$ is conjugate to at least $k$ of the
$\overline{g}_{il}$ this shows that $(\overline{g}_{i1}K,\ldots,\overline
{g}_{im}K)$ has the $(k,\varepsilon)$ fixed-point property. The
$(k,\varepsilon)$ fixed-space property follows likewise from Lemma \ref{GSSh}(iii).

Now the Key Observation is applied in conjunction with Propositions 9.2 and
5.1. Both of these only really need the relevant `$(k,\varepsilon
)$-hypothesis' to be satisfied by the image of the $m$-tuple $\mathbf{y}$ in
$(G/N)^{(m)}$ (see \S 5 and \S 9). As $K\leq N$ this means that we may use Key
Observation (C) just as we used the Key Observation in the preceding subsection.

\medskip

\textbf{Cases 1, 2.} Exactly as before, replacing $1/d$ by $\varepsilon$ where necessary.

\medskip

\textbf{Case 3: }where $N=K$ is quasi-semisimple. Let $S$ denote the (unique)
non-abelian composition factor of $N$. Then Lemma \ref{GSSh}(i) shows that
$\left|  S:M\right|  \geq\left|  S\right|  ^{\varepsilon}$ for each maximal
subgroup $M$ of $S$, so we have $\mu(N)\geq\varepsilon$. The argument then
proceeds as before, with $\varepsilon$ in place of $\mu(q)$. $\blacksquare$

\section{The first inequality: lifting generators}

In this section, we fix a finite $d$-generator group $G$ and an acceptable
quasi-minimal normal subgroup $N$ of $G$. Thus $N$ contains a normal subgroup
$Z$ of $G$ with $Z\leq\mathrm{Frat}(G)$ such that $N/Z$ is a minimal normal
subgroup of $G/Z$, and if $N/Z$ is non-abelian then $N/Z$ is not the product
of fewer than $3$ simple groups.

In the latter case, the composition factors of $N/Z$ are all isomorphic to a
simple group $S$, and we have defined $\mu(S)$ to be the supremum of the
numbers $\tau$ such that%
\[
\left|  S:M\right|  \geq\left|  S\right|  ^{\tau}%
\]
for every maximal subgroup $M$ of $S$. We will write $\mu=\min\{\mu
(S),\,\frac{1}{2}\}$.

Fix positive integers $k$ and $m$ and let $\varepsilon>0$. Recall the

\medskip

\noindent\textbf{Definition} Let $\mathbf{y}=(y_{1},y_{2},\ldots,y_{m})\in
G^{(m)}$.

\begin{description}
\item[(i)] The $m$-tuple $\mathbf{y}$ has the $(k,\varepsilon)$
\emph{fixed-point property} if in any transitive permutation action of $G$ on
a set of size $n\geq2$, at least $k$ of the elements $y_{i}$ move at least
$\varepsilon n$ points.

\item[(ii)] The $m$-tuple $\mathbf{y}$ has the $(k,\varepsilon)$
\emph{fixed-space property }if for every irreducible $\mathbb{F}_{p}G$-module
$V$ of dimension $n\geq2$, where $p$ is any prime, at least $k$ of the $y_{i}$
satisfy $\dim_{\mathbb{F}_{p}}\mathrm{C}_{V}(y_{i})\leq(1-\varepsilon)n$.
\end{description}

We shall prove

\begin{proposition}
\label{prop5.1}Let $y_{1},\ldots,y_{m}\in G$ and assume that $G=\left\langle
y_{1},\ldots,y_{m}\right\rangle N$. Put%
\[
\mathcal{N}(\mathbf{y})=\left\{  \mathbf{a}\in N^{(m)}\mid\left\langle
y_{1}^{a_{1}},\ldots,y_{m}^{a_{m}}\right\rangle \neq G\right\}  .
\]
Let $\varepsilon\in(0,\frac{1}{2}]$.\newline\emph{(i)} Suppose that $N$ is
soluble and that $\mathbf{y}$ has the $(k,\varepsilon)$ fixed-space property.
Then%
\[
\left\vert \mathcal{N}(\mathbf{y})\right\vert \leq\left\vert N\right\vert
^{m}\left\vert N/Z\right\vert ^{d-k\varepsilon}.
\]
\emph{(ii)} There exists an absolute constant $C_{0}$ such that if $N$ is
quasi-semisimple and $\mathbf{y}$ has the $(k,\varepsilon)$ fixed-point
property, where $k\varepsilon\geq\max\{2d+4,\,2C_{0}+2\}$, then
\begin{align*}
\left\vert \mathcal{N}(\mathbf{y})\right\vert  &  <\left\vert N\right\vert
^{m}\text{ \ \ and}\\
\left\vert \mathcal{N}(\mathbf{y})\right\vert  &  \leq\left\vert N\right\vert
^{m}\left\vert N/Z\right\vert ^{1-s}%
\end{align*}
where%
\[
s=\min\{\mu(k\varepsilon/2-d-1),\,\mu(k\varepsilon/2-C_{0})\}.
\]

\end{proposition}

In fact, in (i) the fixed-space property of $\mathbf{y}$ will only be applied
to the action of $G$ on the elementary abelian group $N/Z$, and in (ii) the
fixed-point property of $\mathbf{y}$ will only be applied to the permutation
action of $G$ on the simple factors of $N/Z$; so in both cases it would be
enough to assume that the relevant property is possessed by the image of
$\mathbf{y}$ in $(G/N)^{(m)}.$ (This is used in the proof of Key Theorem (C).)

For $\mathbf{a}\in N^{(m)}$ write%
\[
Y(\mathbf{a})=\left\langle y_{1}^{a_{1}},\ldots,y_{m}^{a_{m}}\right\rangle ,
\]
so $\mathcal{N}(\mathbf{y})=\{\mathbf{a}\in N^{(m)}\mid Y(\mathbf{a})\neq
G\}$. Since $Z\leq\mathrm{Frat}(G)$ we have%
\[
Y(\mathbf{a})\neq G\Longleftrightarrow Y(\mathbf{a})Z\neq G,
\]
so $\mathcal{N}(\mathbf{y})$ is the union of a certain number $r$, say, of
cosets of $Z^{(m)}$. If we show that $r\leq\left|  N/Z\right|  ^{m-t}$ it will
follow that $\left|  \mathcal{N}(\mathbf{y})\right|  \leq\left|  Z\right|
^{m}\left|  N/Z\right|  ^{m-t}=\left|  N\right|  ^{m}\left|  N/Z\right|
^{-t}$. Thus we may replace $G$ by $G/Z$ and so assume henceforth that $Z=1$.

\bigskip

We now proceed with the proof. If $N$ is soluble then it is a simple
$\mathbb{F}_{p}G$-module for some prime $p$, so in case (i) at least $k$ of
the $y_{i}$ satisfy $\dim\mathrm{C}_{N}(y_{i})\leq(1-\varepsilon)n$ where
$n=\dim N$. If $N$ is not soluble then $N=S_{1}\times\cdots\times S_{n}$ where
$n\geq3$ and $G$ permutes the set $\Omega=\{S_{1},\ldots,S_{n}\}$ transitively
by conjugation; so in case (ii) at least $k$ of the $y_{i}$ move at least
$\varepsilon n$ of the factors $S_{j}$; for each such $i$ we have $\left\vert
\mathrm{C}_{N}(y_{i})\right\vert \leq\left\vert N\right\vert ^{1-\varepsilon
/2}$ (cf. Lemma \ref{subdirect} below). Thus in either case, we may relabel
the $y_{i}$ so that%
\begin{equation}
\left\vert \mathrm{C}_{N}(y_{i})\right\vert \leq\left\vert N\right\vert
^{1-\overline{\varepsilon}}\text{ \ for }1\leq i\leq k \label{k-centr}%
\end{equation}
where%
\[
\overline{\varepsilon}=\left\{
\begin{array}
[c]{ccc}%
\varepsilon &  & \text{(}N\text{ soluble)}\\
&  & \\
\varepsilon/2 &  & \text{(}N\text{ insoluble)}%
\end{array}
\right.  .
\]

Now if $\mathbf{a}\in\mathcal{N}(\mathbf{y})$ then $Y(\mathbf{a})\leq M$ for
some maximal subgroup $M$ of $G$ with $NM=G$. Write%
\[
v(M;y)=|\{a\in N\mid y^{a}\in M\}|
\]
and%
\[
v(M)=\prod_{i=1}^{m}v(M;y_{i}).
\]
Then $v(M)$ is just the number of $\mathbf{a}$ such that $Y(\mathbf{a})\leq
M$, so%
\begin{equation}
\left|  \mathcal{N}(\mathbf{y})\right|  \leq\sum_{M\in\mathcal{M}}v(M)
\label{N(y) and v(M)}%
\end{equation}
where $\mathcal{M}$ denotes the set of maximal subgroups of $G$ which
supplement $N$.

\begin{lemma}
\label{v(M,y)}Let $M\in\mathcal{M}$, $y\in G$ and put $D=M\cap N$. Then%
\[
v(M;y)=\left|  \mathrm{C}_{N}(y)\right|  \cdot\left|  \lbrack y^{b},N]\cap
D\right|
\]
for every $b\in N$ such that $y^{b}\in M$, and $v(M;y)=0$ if there is no such
$b$.
\end{lemma}

\begin{proof}
If no conjugate of $y$ lies in $M$ then $v(M;y)=0$. Otherwise, $y^{b}\in M$
for some $b\in N$; given any such $b$, for $a\in N$ we have%
\[
y^{a}\in M\Longleftrightarrow\lbrack y^{b},b^{-1}a]\in M\cap\lbrack
y^{b},N]=[y^{b},N]\cap D.
\]
The lemma follows since the fibres of the mapping $a\mapsto\lbrack
y^{b},b^{-1}a]$ are cosets of $\mathrm{C}_{N}(y)$.
\end{proof}

\bigskip

Let $\mathcal{M}_{0}$ denote the set of all $M\in\mathcal{M}$ such that $M\cap
N=1$.

\begin{lemma}
\label{M0}\emph{(i)}%
\[
\left|  \mathcal{M}_{0}\right|  \leq\left|  N\right|  ^{d}.
\]

\emph{(ii) }If $M\in\mathcal{M}_{0}$ then
\[
v(M)\leq\left\vert N\right\vert ^{m-k\overline{\varepsilon}}.
\]

\end{lemma}

\begin{proof}
(i) Follows from the well-known fact that the complements to $N$ in $G$, if
there are any, correspond bijectively to derivations from $G/N$ to $N$, and
the fact that $G$ can be generated by $d$ elements.

(ii) Since $N\cap M=1$, Lemma \ref{v(M,y)} and (\ref{k-centr}) give%
\[
v(M)=\prod_{j=1}^{m}v(M;y_{j})\leq\left\vert N\right\vert ^{k(1-\overline
{\varepsilon})}\cdot\left\vert N\right\vert ^{m-k}=\left\vert N\right\vert
^{m-k\overline{\varepsilon}}.
\]

\end{proof}

\bigskip

Part (i) of the proposition now follows: for when $N$ is abelian we have
$\mathcal{M=M}_{0}$ and so%
\[
\left|  \mathcal{N}(\mathbf{y})\right|  \leq\left|  N\right|  ^{d}\cdot\left|
N\right|  ^{m-k\varepsilon}%
\]
as required.

We assume henceforth that $N$ is \emph{non-abelian}; thus%
\[
N=S_{1}\times\cdots\times S_{n}%
\]
where $n\geq3$ and the $S_{i}$ are isomorphic simple groups. The conjugation
action of $G$ permutes the factors $S_{i}$ transitively, and we write%
\[
S_{i}^{g}=S_{i\sigma(g)}%
\]
where $\sigma(g)\in\mathrm{Sym}(n)$.

For a natural number $e$ put%
\[
\mathcal{M}(e)=\left\{  M\in\mathcal{M}\mid\left|  G:M\right|  =e\right\}  .
\]
Thus $\mathcal{M}(\left|  N\right|  )=\mathcal{M}_{0}$, and $\mathcal{M}(e)$
is non-empty only when $e\geq2$ and $e$ is a divisor of $\left|  N\right|  .$

\begin{lemma}
\label{M(e)}There is an absolute constant $C$ such that%
\[
\left|  \mathcal{M}(e)\right|  \leq e^{C}%
\]
for every proper divisor $e$ of $\left|  N\right|  $.
\end{lemma}

\begin{proof}
Let $M\in\mathcal{M}(e)$ and put $D=M\cap N$. Since $\left|  N:D\right|
=\left|  G:M\right|  =e$ we have $1<D<N$, so $D$ is not normal in $G$. As
$D\vartriangleleft M$ it follows that $M\geq\mathrm{C}_{G}(N)$. It is now
clear that $\mathrm{C}_{G}(N)$ is the core of $M$, that is, the biggest normal
subgroup of $G$ contained in $M$.

Thus $M\mapsto M/\mathrm{C}_{G}(N)$ maps $\mathcal{M}(e)$ bijectively onto the
set of core-free maximal subgroups in $G/\mathrm{C}_{G}(N)$ that supplement
but do not complement $N\mathrm{C}_{G}(N)/\mathrm{C}_{G}(N)$ and have index
$e$. It is proved by Mann and Shalev in \cite{MS} that the cardinality of this
set is bounded by $e^{C}$ where $C$ is an absolute constant: see the first
part of the proof of \cite{MS}, Corollary 2.
\end{proof}

\bigskip

\begin{lemma}
\label{subdirect}Let $V=A_{1}\times\cdots\times A_{t}$ where $t\geq2$ and the
$A_{i}$ are isomorphic finite groups. Let $g$ be an automorphism of $V$ that
permutes the subgroups $A_{i}$ and moves at least $\varepsilon t$ of them.

\emph{(i)} Let $U=B_{1}\times\cdots\times B_{t}$ where $B_{i}<A_{i}$ and
$\left|  B_{i}\right|  =\left|  B_{1}\right|  $ for each $i$. Suppose that
$U^{g}=U$. Then%
\[
\left|  \mathrm{C}_{V}(g)\right|  \cdot\left|  \lbrack g,V]\cap U\right|
\leq\left|  V\right|  \cdot\left|  V:U\right|  ^{-\varepsilon/2}.
\]

\emph{(ii)} Let $\Delta\cong A_{1}$ be a diagonal subgroup of $V$. Suppose
that $\Delta^{g}=\Delta$ and that $t\geq3$. Then%
\[
\left|  \mathrm{C}_{V}(g)\right|  \cdot\left|  \lbrack g,V]\cap\Delta\right|
\leq\left|  V\right|  \cdot\left|  V:\Delta\right|  ^{-\varepsilon/2}.
\]

\end{lemma}

\begin{proof}
Write $A=A_{1},\,B=B_{1}$. Consider a typical cycle for the permutation action
of $g,$ say $\mathcal{C}=(A_{1},\ldots,A_{l})$, and put $V_{\mathcal{C}}%
=A_{1}\times\cdots\times A_{l}$. Note that for $1\leq i\leq l$ we have%
\[
B_{i}=U\cap A_{1}^{g^{i-1}}=(U\cap A_{1})^{g^{i-1}}=B_{1}^{g^{i-1}}.
\]
A typical element of $V_{\mathcal{C}}$ takes the form%
\[
v=a_{1}\cdot a_{2}^{g}\cdot\ldots\cdot a_{l}^{g^{l-1}}%
\]
where $a_{i}\in A_{1}$ for each $i$. Then%
\begin{equation}
\lbrack g,v]=a_{l}^{-g^{l}}a_{1}\cdot(a_{1}^{-1}a_{2})^{g}\cdot\ldots
\cdot(a_{l-1}^{-1}a_{l})^{g^{l-1}} \label{[g,v]}%
\end{equation}
so $[g,v]\in U$ if and only if%
\begin{align*}
a_{i}^{-1}a_{i+1}  &  \in B_{1}\,\,(1\leq i\leq l-1)\\
\theta(a_{l})^{-1}a_{1}  &  \in B_{1}%
\end{align*}
where $\theta$ denotes the automorphism induced on $A_{1}$ by $g^{l}$. Hence
putting%
\begin{align*}
X  &  =\left\{  v\in V_{\mathcal{C}}\mid\lbrack g,v]\in U\right\}  ,\\
Y  &  =\left\{  y\in A_{1}\mid\theta(y)^{-1}y\in B_{1}\right\}
\end{align*}
we obtain a bijection%
\begin{align*}
B_{1}^{(l-1)}\times Y  &  \rightarrow X\\
(b_{1},\ldots,b_{l-1},y)  &  \mapsto(\theta(y)b_{1})\cdot(\theta(y)b_{2}%
)^{g}\cdot\ldots\cdot\left(  \theta(y)b_{l-1}\right)  ^{g^{l-2}}\cdot
y^{g^{l-1}}.
\end{align*}
Since the fibres of the map $v\mapsto\lbrack g,v]$ are cosets of
$\mathrm{C}_{V_{\mathcal{C}}}(g)$ it follows that%
\begin{align*}
\left|  \lbrack g,V_{\mathcal{C}}]\cap U\right|  \cdot\left|  \mathrm{C}%
_{V_{\mathcal{C}}}(g)\right|   &  =\left|  X\right| \\
&  =\left|  B_{1}\right|  ^{l-1}\left|  Y\right|  \leq\left|  B\right|
^{l-1}\left|  A\right|  .
\end{align*}

Now $V$ as a $\left\langle g\right\rangle $-operator group is the direct
product of the $V_{\mathcal{C}}$ over all the cycles $\mathcal{C}%
=\mathcal{C}_{1},\ldots,\mathcal{C}_{p}$ say. It follows that%
\begin{align*}
\left|  \lbrack g,V]\cap U\right|  \cdot\left|  \mathrm{C}_{V}(g)\right|   &
=\prod_{i=1}^{p}\left|  [g,V_{\mathcal{C}_{i}}]\cap U\right|  \cdot\left|
\mathrm{C}_{V_{\mathcal{C}_{i}}}(g)\right| \\
&  \leq\prod_{i=1}^{p}\left(  \left|  B\right|  ^{l_{i}-1}\left|  A\right|
\right)  =\left|  A\right|  ^{t}\left|  A:B\right|  ^{p-t}%
\end{align*}
where $l_{i}$ is the length of $\mathcal{C}_{i}$. Since at most
$(1-\varepsilon)t$ of the $l_{i}$ are equal to $1$ we have $p-t\leq
-\varepsilon t/2$. Hence%
\[
\left|  \lbrack g,V]\cap U\right|  \cdot\left|  \mathrm{C}_{V}(g)\right|
\leq\left|  A\right|  ^{t}\left|  A:B\right|  ^{-\varepsilon t/2}=\left|
V\right|  \left|  V:U\right|  ^{-\varepsilon/2}%
\]
and (i) is proved.

Taking $U=1$ in (i) we deduce that $\left|  \mathrm{C}_{V}(g)\right|
\leq\left|  A\right|  ^{p}$. Since $\left|  \Delta\right|  =\left|  A\right|
$ it follows that%
\[
\left|  \lbrack g,V]\cap\Delta\right|  \cdot\left|  \mathrm{C}_{V}(g)\right|
\leq\left|  A\right|  ^{1+p}.
\]
Suppose first that $l_{i}\geq2$ for each $i$. Then $p\leq t/2$, and as
$t\geq3$ we have%
\[
1+p\leq t-(t-1)/4.
\]
It follows that $\left|  A\right|  ^{1+p}\leq\left|  V\right|  \left|
V:\Delta\right|  ^{-1/4}$, and (ii) follows since $\varepsilon\leq\frac{1}{2}$.

Now suppose that one of the $l_{i}$ is equal to $1$, say $l_{1}=1$. Then $g$
fixes $A_{1}$. Each element of $[g,V]\cap\Delta$ is determined by its first
component, which belongs to $[g,A_{1}]$. Applying part (i) to $V^{\ast}%
=A_{2}\times\cdots\times A_{t}$, with each $B_{i}=1$, we deduce as above that
$\left|  \mathrm{C}_{V^{\ast}}(g)\right|  \leq\left|  A\right|  ^{p-1}$ and
hence that%
\[
\left|  \mathrm{C}_{V}(g)\right|  =\left|  \mathrm{C}_{A_{1}}(g)\right|
\left|  \mathrm{C}_{V^{\ast}}(g)\right|  \leq\left|  \mathrm{C}_{A_{1}%
}(g)\right|  \left|  A\right|  ^{p-1}.
\]
It follows that%
\[
\left|  \lbrack g,V]\cap\Delta\right|  \cdot\left|  \mathrm{C}_{V}(g)\right|
\leq\left|  \lbrack g,A_{1}]\right|  \left|  \mathrm{C}_{A_{1}}(g)\right|
\left|  A\right|  ^{p-1}=\left|  A\right|  ^{p}.
\]
As $p-t\leq-\varepsilon t/2<-\varepsilon(t-1)/2$ we have $\left|  A\right|
^{p}<\left|  V\right|  \left|  V:\Delta\right|  ^{-\varepsilon/2}$, again
giving (ii).
\end{proof}

\bigskip

Now fix $e$ with $2\leq e<\left|  N\right|  $ and consider $M\in
\mathcal{M}(e)$. Put $D=M\cap N$ and let $R_{i}$ denote the projection of $D$
into $S_{i}$. It is easy to see that if $g\in M$ then $R_{i}^{g}%
=R_{i\sigma(g)}$ for each $i$, so the group $\widetilde{R}=R_{1}\times
\cdots\times R_{n}$ is normalized by $M$.

Say $M$ is of \emph{type 1} if $\widetilde{R}M\neq G$. In this case
$D=\widetilde{R}$. Put $t=n$, $\ A_{i}=S_{i}$ and $B_{i}=R_{i}$. Note that%
\[
e=\left|  N:D\right|  =\left|  S_{1}:R_{1}\right|  ^{n}\geq\left|
S_{1}\right|  ^{\mu n}=\left|  N\right|  ^{\mu}.
\]

Now suppose that $\widetilde{R}M=G$. Then $\widetilde{R}=N$ and $D$ is a
subdirect product in $N=S_{1}\times\cdots\times S_{n}$. In this case, we can
re-label the $S_{i}$ so that%
\[
D=B_{1}\times\cdots\times B_{t^{\prime}}%
\]
where $t^{\prime}\mid n$ and for each $i$, $B_{i}$ is a diagonal subgroup of
$S_{r(i-1)+1}\times\cdots\times S_{r(i-1)+r}$ with $r=n/t^{\prime}\geq2$
(\cite{Cm}, Exercise 4.3). If $t^{\prime}\geq2$ say $M$ is of \emph{type 2},
and put $t=t^{\prime}$, $A_{i}=S_{r(i-1)+1}\times\cdots\times S_{r(i-1)+r}$.

If $t^{\prime}=1$ say $M$ is of \emph{type 3}, put $t=n$, take $A_{i}=S_{i}$
for each $i$ and put $\Delta=D$.

Again, we have%
\[
e=\left|  N:D\right|  =\left|  S_{1}\right|  ^{n-t^{\prime}}=\left|  N\right|
^{1-r^{-1}}\geq\left|  N\right|  ^{1/2}\geq\left|  N\right|  ^{\mu}.
\]

In each case, the action of $M$ permutes the $A_{i}$ transitively. Since
$G=NM$ it follows that $G$ also permutes the $A_{i}$ transitively. Writing
$\mathcal{J}$ to denote the set of subscripts $l\leq m$ such that $y_{l}$
moves at least $\varepsilon t$ of the $A_{i}$, we have $\left|  \mathcal{J}%
\right|  \geq k$ by the $(k,\varepsilon)$ fixed-point property of $\mathbf{y}$.

Now let $l\in\mathcal{J}$. According to Lemma \ref{v(M,y)}, if no
$N$-conjugate of $y_{l}$ lies in $M$ then $v(M;y_{l})=0$; while if $y_{l}%
^{b}\in M$ where $b\in N$ then%
\[
v(M;y_{l})=\left|  \mathrm{C}_{N}(y_{l})\right|  \cdot\left|  \lbrack
y_{l}^{b},N]\cap M\right|  .
\]
Put $g=y_{l}^{b}$. Then $g$ also moves at least $\varepsilon t$ of the $A_{i}%
$. Putting $V=N$, and $U=D$ when $M$ is of types 1 or 2, we may apply Lemma
\ref{subdirect} to deduce that%
\begin{align*}
\left|  \mathrm{C}_{N}(g)\right|  \cdot\left|  \lbrack g,N]\cap D\right|   &
\leq\left|  N\right|  \cdot\left|  N:D\right|  ^{-\varepsilon/2}\\
&  =e^{-\varepsilon/2}\left|  N\right|  .
\end{align*}

As $\left\vert \mathrm{C}_{N}(g)\right\vert =\left\vert \mathrm{C}_{N}%
(y_{l})\right\vert $ this shows that%
\[
v(M;y_{l})\leq e^{-\varepsilon/2}\left\vert N\right\vert .
\]
Hence%
\[
v(M)=\prod_{i=1}^{m}v(M;y_{i})\leq\left(  e^{-\varepsilon/2}\left\vert
N\right\vert \right)  ^{\left\vert \mathcal{J}\right\vert }\cdot\left\vert
N\right\vert ^{m-\left\vert \mathcal{J}\right\vert }\leq e^{-k\varepsilon
/2}\left\vert N\right\vert ^{m}.
\]
This holds for each $M\in\mathcal{M}(e)$. With Lemmas \ref{M0} and \ref{M(e)}
it gives%
\begin{align*}
\left\vert \mathcal{N}(\mathbf{y})\right\vert  &  \leq\sum_{M\in\mathcal{M}%
}v(M)\\
&  \leq\left\vert N\right\vert ^{m+d-k\varepsilon/2}+\sum_{e}e^{C-k\varepsilon
/2}\left\vert N\right\vert ^{m},
\end{align*}
where $e$ ranges over integers lying between $\left\vert N\right\vert ^{\mu}$
and $\left\vert N\right\vert /2.$

We can now deduce part (ii) of Proposition \ref{prop5.1}. Take $C_{0}=C+1$,
and assume that%
\[
k\varepsilon\geq\max\{2d+4,\,2C_{0}+2\}.
\]
Put $t=\min\{k\varepsilon/2-d,\,k\varepsilon/2-C\}$, write $\nu=\left\vert
N\right\vert $ and let $\zeta$ denote the Riemann zeta function. Then $t\geq
2$, so%
\[
\left\vert N\right\vert ^{-m}\left\vert \mathcal{N}(\mathbf{y})\right\vert
\leq\sum_{e\geq2}e^{-t}\leq\zeta(2)-1<1.
\]
This establishes the first claim. For the second, observe that $\zeta
(t)<2<\nu$, and so%
\[
\left\vert N\right\vert ^{-m}\left\vert \mathcal{N}(\mathbf{y})\right\vert
\leq\sum_{e\geq\nu^{\mu}}e^{-t}\leq\zeta(t)\nu^{-\mu(t-1)}\leq\nu^{1-s}%
\]
where%
\[
s=\mu(t-1)=\min\left\{  \mu(k\varepsilon/2-d-1),\,\mu\left(  k\varepsilon
/2-C_{0}\right)  \right\}  .
\]

\section{Exterior squares and quadratic maps\label{quadmaps}}

In the following section we are going to prove Proposition 7.1. This concerns
the solution of certain equations in a soluble quasi-minimal normal subgroup
$N$ of a finite group $G$. When $N$ is \emph{abelian} (`Case 1') the result is
very easy. When $N$ is \emph{non-abelian}, the problem comes down to studying
the fibres over $N^{\prime}=[N,N]$ of certain mappings $\phi_{i}$ from
$(N/Z)^{(m)}$ into $N$ (induced by commutation with certain elements of $G$);
here $Z=\mathrm{C}_{N}(G)$, $N/Z$ is a simple $\mathbb{F}_{p}G$-module for
some prime $p$, and $N^{\prime}$ is an $\mathbb{F}_{p}$-module contained in
$Z$. If $\left|  N^{\prime}\right|  =2$ (`Case 2') it turns out that the
restriction of each $\phi_{i}$ to $\phi_{i}^{-1}(N^{\prime})=V_{i}$ is a
quadratic form over $\mathbb{F}_{2}$, and the required result follows from
some elementary number theory over $\mathbb{F}_{2}$. The hardest case (`Case
3') is when $\left|  N^{\prime}\right|  >2$. The mappings $\phi_{i\left|
V_{i}\right.  }$ are still quadratic polynomial mappings over $\mathbb{F}_{p}%
$, but we may no longer suppose that their co-domain $N^{\prime}$ is
one-dimensional over $\mathbb{F}_{p}$, and higher-dimensional algebraic
geometry does not deliver the result.

To get round this difficulty, we would like to think of $N^{\prime}$ as a
one-dimensional space over a larger field. Such a structure does not arise
naturally, in general; however, $N^{\prime}$ is an epimorphic image of the
exterior square of $N/Z$, and it was shown in \cite{Sg} that the latter does
naturally have the structure of a one-dimensional space over a certain field.
This is the key to the main result of this section, Proposition \ref{phi-equn}%
, which in turn will serve to complete case 3 of the proof of Proposition 7.1.

When $p$ is odd, everything needed for the proof essentially appears in
\cite{Sg}; but the proof given in that paper for the `even' case depends
crucially on a global solubility assumption, not available to us here, and a
new approach is required. In fact we shall deal in a uniform way with the
`odd' and `even' cases, by strengthening the method used for the `odd' part in
\cite{Sg} (and the very tricky material of \S \S 8 and 9 of \cite{Sg} may now
be consigned to a historical footnote).

\bigskip We need to recall some material from \cite{Sg}, \S 4. Let
$\overline{G}$ be a group (assumed finite in \cite{Sg}, but this is not
necessary) and $R=\mathbb{Z}\overline{G}$ the group ring. Let $M$ be a finite
simple right $R$-module, so $M$ is an $\mathbb{F}_{p}\overline{G}$-module for
some prime $p$. We may consider $M$ as an $R$-bimodule via%
\[
gu=ug^{-1}\,\,\,\,(u\in M,\,g\in\overline{G}),
\]
and so define $M\otimes_{R}M$ and the exterior square%
\[
\wedge_{R}^{2}M\cong\left(  \wedge_{\mathbb{F}_{p}}^{2}M\right)  /\left(
\wedge_{\mathbb{F}_{p}}^{2}M\right)  (\overline{G}-1)
\]
(where $\overline{G}$ acts diagonally on $\wedge_{\mathbb{F}_{p}}^{2}M$). We
fix a generator $\overline{c}$ for $M$ and put%
\begin{align*}
I  &  =\mathrm{ann}_{R}(\overline{c})\\
S_{0}  &  =\left\{  r\in R\mid rI\subseteq I\right\}  .
\end{align*}
The ring $S_{0}/I$ may be identified with the finite field $\mathrm{End}%
_{R}(M)$ via $s+I\mapsto\widehat{s}$ where%
\[
\widehat{s}(\overline{c}r)=\overline{c}sr\,\,(s\in S_{0},\,r\in R).
\]
Suppose now that $M$ admits a non-zero $\overline{G}$-invariant alternating
$\mathbb{F}_{p}$-bilinear form. According to Proposition 4.4 of \cite{Sg},
there exists a subfield $k=S/I\subseteq S_{0}/I$ such that for each $s\in S,$
$a,\,b\in M,$%
\[
\widehat{s}a\otimes b=a\otimes\widehat{s}b
\]
holds in $M\otimes_{R}M$, and such that the induced action of $k$ on
$\wedge_{R}^{2}M$ makes $\wedge_{R}^{2}M$ into a $1$-dimensional vector space
over $k$. Moreover, $\dim_{k}(M)\geq2$, and the mapping%
\[
(a,b)\mapsto a\wedge b
\]
from $M\times M$ to $\wedge_{R}^{2}M$ is $k$-bilinear.

We shall consider $M$ and $\wedge_{R}^{2}M$ as left $S$-modules via
$S\rightarrow k$.\bigskip

Now assume that we are given a group $G$ and a normal subgroup $B$ such that
$G/B=\overline{G}$. Let $A/B$ be a minimal normal subgroup of $\overline{G}$
such that $A/B=M$ as a $\overline{G}$-module via conjugation. Assume also
that
\begin{align*}
\lbrack B,A]  &  =[A^{\prime},G] =1,\\
\left|  A^{\prime}:A^{\prime}\cap\lbrack B,G]\right|   &  >2
\end{align*}
and that the mapping $aB\wedge bB\mapsto\lbrack a,b]$ ($a,\,b\in A$) induces
an isomorphism%
\begin{equation}
\wedge_{R}^{2}M\rightarrow A^{\prime}. \label{1}%
\end{equation}

These hypotheses imply that $M$ does admit a non-zero $\overline{G}$-invariant
alternating $\mathbb{F}_{p}$-bilinear form: there exists an epimorphism
$\theta:A^{\prime}\rightarrow\mathbb{F}_{p}$ and then $(aB,bB)\mapsto
\theta([a,b])$ is such a form. We may therefore identify $A^{\prime}$ with the
one-dimensional $k$-space $\wedge_{R}^{2}M$ via (\ref{1}), and will use
additive and multiplicative notation interchangeably for the group operation
there. Note that%
\begin{align*}
\left|  k\right|   &  =\left|  A^{\prime}\right|  >2,\\
\left|  M\right|   &  =\left|  k\right|  ^{\dim_{k}M}\geq\left|  A^{\prime
}\right|  ^{2}.
\end{align*}
Fix $c\in A$ such that $cB=\overline{c}$, the chosen generator of $M$. Suppose
that $\overline{G}$ can be generated by $d$ elements.

\begin{proposition}
\label{quadmap}Let $x_{1},\ldots,x_{m}\in G$ satisfy $B\left\langle
x_{1},\ldots,x_{m}\right\rangle =G$. Then there exist \emph{(a)} a
$k$-subspace $U$ of $M^{(m)}$, \emph{(b)} a $k$-quadratic map $\Phi
:U\rightarrow A^{\prime}$, and \emph{(c)} for each $\mathbf{z}=(z_{1}%
,\ldots,z_{m})\in A^{(m)}$, a $k$-linear map $\alpha^{\mathbf{z}}:U\rightarrow
A^{\prime}$ such that\medskip

\emph{(i)} $\mathrm{dim}_{k}U\geq(m-d)\mathrm{dim}_{k}M\medskip$

\emph{(ii)} for each $u\in U$ there exist $a_{1},\ldots,a_{m}\in A$ with
$(a_{1}B,\ldots,a_{m}B)=u$ such that%
\[
\Phi(u)+\alpha^{\mathbf{z}}(u)=\left(  \prod_{j=1}^{m}[z_{j}a_{j}%
,x_{j}]\right)  \cdot\left(  \prod_{j=1}^{m}[z_{j},x_{j}]\right)  ^{-1}.
\]
Moreover $\alpha^{(1,\ldots,1)}=\alpha$ is surjective.
\end{proposition}

Before proving this let us deduce its primary application:

\begin{proposition}
\label{phi-equn}Let $\ast:G\rightarrow G^{\ast}$ be an epimorphism with
$\ker(\ast)\leq B$ and $[B^{\ast},G^{\ast}]=1$. For $i=1,2,3$ let
$x_{i1},\ldots,x_{im}\in G^{\ast}$ satisfy $B^{\ast}\left\langle x_{i1}%
,\ldots,x_{im}\right\rangle =G^{\ast}$, and define%
\[
\phi_{i}:A^{\ast(m)}\rightarrow A^{\ast}%
\]
by%
\[
\phi_{i}(a_{1},\ldots,a_{m})=\prod_{j=1}^{m}[a_{j},x_{ij}].
\]
Then for each $\kappa\in(A^{\ast})^{\prime}$ there exist $\kappa_{1}%
,\,\kappa_{2},\,\kappa_{3}\in A^{\ast}$ such that%
\begin{equation}
\kappa_{1}\kappa_{2}\kappa_{3}=\kappa\label{2}%
\end{equation}
and%
\begin{equation}
\phi_{i}^{-1}(\kappa_{i})\text{ contains at least }\left|  M\right|
^{m-d-1}\text{ cosets of }B^{\ast(m)}\text{ \ \ }(i=1,2,3), \label{3}%
\end{equation}
provided in case $p=2$ that $A^{\ast}=[A^{\ast},G^{\ast}]$ and $(A^{\ast}%
)^{2}=(A^{\ast})^{\prime}$.
\end{proposition}

\begin{proof}
Let $\widetilde{x}_{ij}$ denote a preimage in $G$ of $x_{ij}\in G^{\ast}$, and
let $\Phi_{i},\,\alpha_{i}:U_{i}\rightarrow A^{\prime}$ be the mappings
corresponding to $(\widetilde{x}_{i1},\ldots,\widetilde{x}_{im})$ provided in
Proposition \ref{quadmap}. Let $\widetilde{\kappa}\in A^{\prime}$ be a
preimage of $\kappa$. Write $\Phi_{i}+\alpha_{i}=f_{i}$. Note that $f_{i}$ is
not the zero map, because $\alpha_{i}$ is surjective and $\left|  k\right|
>2$, which implies that a non-zero map cannot be both linear and quadratic
over $k$; and that for each $u\in U_{i}$ there exist $a_{1},\ldots,a_{m}\in A$
such that $u=(a_{1}B,\ldots,a_{m}B)$ and%
\begin{equation}
f_{i}(u)^{\ast}=\phi_{i}(a_{1}^{\ast},\ldots,a_{m}^{\ast}). \label{4}%
\end{equation}
Since $[B^{\ast},G^{\ast}]=1$, this then holds for \emph{every} $m$-tuple
$(a_{1},\ldots,a_{m})$ with $(a_{1}B,\ldots,a_{m}B)=u$. Similarly, if
$\mathbf{z}=(z_{1},\ldots,z_{m})\in A^{(m)}$,$\,\,\alpha_{i}^{\mathbf{z}}$ are
as given in Proposition \ref{quadmap} and $f_{i}^{\mathbf{z}}=\Phi_{i}%
+\alpha_{i}^{\mathbf{z}}$, then%
\begin{equation}
f_{i}^{\mathbf{z}}(u)^{\ast}=\phi_{i}(z_{1}^{\ast}a_{1}^{\ast},\ldots
,z_{m}^{\ast}a_{m}^{\ast})\phi_{i}(z_{1}^{\ast},\ldots,z_{m}^{\ast})^{-1}
\label{4bis}%
\end{equation}
whenever $(a_{1}B,\ldots,a_{m}B)=u$.

\medskip

\textbf{Case 1:} where $p\neq2$. First we pick $\kappa_{3}$. The fibres of the
map $f_{3}:U_{3}\rightarrow A^{\prime}$ have average size at least $\left|
U_{3}\right|  /\left|  A^{\prime}\right|  >\left|  M\right|  ^{m-d-1}$, so
there exists $\widetilde{\kappa}_{3}\in A^{\prime}$ with $\left|  f_{3}%
^{-1}(\widetilde{\kappa}_{3})\right|  >\left|  M\right|  ^{m-d-1}$. Now put
$\kappa_{3}=\widetilde{\kappa}_{3}^{\ast}$. Then (\ref{4}) implies that
$\phi_{3}^{-1}(\kappa_{3})$ contains at least $\left|  M\right|  ^{m-d-1}$
cosets of $B^{\ast(m)}$.

Next, let $\widetilde{\kappa}_{4}$ be a preimage of $\kappa\kappa_{3}^{-1}$.
According to Lemma 5.1 of \cite{Sg} there exist elements $\widetilde{\kappa
}_{1},\,\widetilde{\kappa}_{2}\in A^{\prime}$ such that $\widetilde{\kappa
}_{1}+\widetilde{\kappa}_{2}=\widetilde{\kappa}_{4}$ and%
\begin{align}
\left|  f_{i}^{-1}(\widetilde{\kappa}_{i})\right|   &  \geq\left|  k\right|
^{\mathrm{dim}_{k}U_{i}-2}\nonumber\\
&  \geq\left|  M\right|  ^{m-d}\left|  k\right|  ^{-2}\geq\left|  M\right|
^{m-d-1}\,\,(i=1,2). \label{5}%
\end{align}
Now put $\kappa_{i}=\widetilde{\kappa}_{i}^{\ast}$ for $i=1,\,2$. Then
$\kappa_{1}\kappa_{2}\kappa_{3}=\kappa$, and (\ref{3}) for $i=1,\,2$ follows
from (\ref{4}) and (\ref{5}).

\medskip

\textbf{Case 2:} where $p=2$. According to Lemma 5.2 of \cite{Sg},
$f_{i}(U_{i})=P_{i}$, say, is a subgroup of index at most $2$ in $A^{\prime}$
for each $i$, and (\ref{5}) holds for each $\widetilde{\kappa}_{i}\in P_{i}%
$.\medskip

\emph{Subcase 2.1: } $P_{t}\neq P_{l}$ for some pair $t,\,l$. Then
$P_{1}+P_{2}+P_{3}=A^{\prime}$ so there exist $\widetilde{\kappa}_{i}%
=f_{i}(u_{i})\in P_{i}$ such that $\widetilde{\kappa}_{1}+\widetilde{\kappa
}_{2}+\widetilde{\kappa}_{3}=\widetilde{\kappa}$. Then both (\ref{2}) and
(\ref{3}) hold with $\kappa_{i}=\widetilde{\kappa}_{i}^{\ast}$, as in Case
1.\medskip

\emph{Subcase 2.2:} $P_{1}=P_{2}=P_{3}=P$ say, with $\left|  A^{\prime
}:P\right|  \leq2$. According to the extra hypotheses in Case 2, there exists
$a\in A^{\ast}$ such that $a^{2}\equiv\kappa\,(\operatorname{mod}P^{\ast}),$
and there exist $v_{i}\in A^{\ast(m)}$ such that $\phi_{i}(v_{i})\equiv
a\,(\operatorname{mod}A^{\ast\prime})$ (because $A^{\ast}/A^{\ast\prime}$ is a
perfect $G^{\ast}/B^{\ast}$-module and $x_{i1},\ldots,x_{im}$ generate
$G^{\ast}$ modulo $B^{\ast}$). By the pigeonhole principle, there exist $t<l$
such that $a^{-1}\phi_{t}(v_{t})\equiv a^{-1}\phi_{l}(v_{l}%
)\,(\operatorname{mod}P^{\ast});$ as $A^{\ast\prime}$ has exponent $2$ we then
have
\[
\phi_{t}(v_{t})\phi_{l}(v_{l})\equiv a^{2}\equiv\kappa\,(\operatorname{mod}%
P^{\ast}).
\]
Now put%
\begin{align*}
\kappa_{t}=\phi_{t}(v_{t}),\,\,\kappa_{l}=\phi_{l}(v_{l}),\\
\kappa_{j}=\left(  \phi_{t}(v_{t})\phi_{l}(v_{l})\right)  ^{-1}\kappa
\end{align*}
where $\{t,l,j\}=\{1,2,3\}$. Then%
\[
\kappa_{j}\in P^{\ast}=f_{j}(U_{j})^{\ast}\subseteq\phi_{j}(A^{\ast(m)}),
\]
and $\kappa_{1}\kappa_{2}\kappa_{3}=\kappa$ since $P^{\ast}$ is central.

To establish (\ref{3}), it now suffices to show that for each $i$ and each
$v\in A^{\ast(m)},$ the fibre $\phi_{i}^{-1}(\phi_{i}(v))$ contains at least
$\left|  M\right|  ^{m-d-1}$ cosets of $B^{\ast m}$. Say $v=(z_{1}^{\ast
},\ldots,z_{m}^{\ast})$. Then for each $u\in(f_{i}^{\mathbf{z}})^{-1}(0)$ we
have%
\[
\phi_{i}(z_{1}^{\ast}a_{1}^{\ast},\ldots,z_{m}^{\ast}a_{m}^{\ast})=\phi
_{i}(v)
\]
whenever $(a_{1}B,\ldots,a_{m}B)=u$, by (\ref{4bis}). Our claim now follows
from (\ref{5}) with $f_{i}^{\mathbf{z}}$ in place of $f_{i}$ and $0$ for
$\widetilde{\kappa}_{i}$.
\end{proof}

\bigskip

We turn now to the proof of Proposition \ref{quadmap}. For $a\in A$ and $g\in
G$ we shall write%
\begin{align*}
\overline{g}  &  =gB\in G/B=\overline{G}\\
\overline{a}  &  =aB\in A/B=M\\
\widetilde{a}  &  =aA^{\prime}\in A/A^{\prime}=\widetilde{A}.
\end{align*}
Since $[B,A]=[A^{\prime},G]=1$, for $a,\,d\in A$ and $g\in G$ we may set%
\begin{align*}
\lbrack\overline{a},d]  &  =[a,d],\\
\lbrack\widetilde{a},g]  &  =[a,g],\\
a^{\overline{g}}=a^{g},  &  \,\,[a,\overline{g}]=[a,g].
\end{align*}

\begin{lemma}
\label{lem1}Let $g_{1},\ldots,g_{n}\in\overline{G},$ $\varepsilon_{1}%
,\ldots,,\varepsilon_{n}\in\{1,\,-1\}$ satisfy%
\[
\sum_{j=1}^{n}\varepsilon_{j}g_{j}=0
\]
in the group ring $R=\mathbb{Z}\overline{G}$. For $a\in A$ let%
\[
\psi(a)=\prod_{j=1}^{n}a^{\varepsilon_{j}g_{j}}.
\]
Then there exist $h_{i},\,k_{i}\in\{g_{1},\ldots,g_{n}\}$ ($i=1,\ldots,l$)
such that for each $a\in A$%
\[
\psi(a)=\prod_{i=1}^{l}[a^{h_{i}},a^{k_{i}}].
\]

\end{lemma}

\begin{proof}
The hypothesis implies that $n=2t$ is even and that the sequence
$(\varepsilon_{1}g_{1},\ldots,\varepsilon_{n}g_{n})=\mathcal{S}$ is some
re-arrangement of $(y_{1},-y_{1},\ldots,y_{t},-y_{t})=\mathcal{S}^{\prime}$
where each $y_{i}$ is one of the $g_{j}$. Since $A^{\prime}$ is central in
$A$, it follows that for each $a\in A$ we have%
\[
\psi(a)=\prod_{i=1}^{t}a^{y_{i}}a^{-y_{i}}\cdot\mathfrak{x}(a)=\mathfrak{x}%
(a)
\]
where $\mathfrak{x}(a)$ is the product of certain factors of the form
$[a^{\varepsilon_{i}g_{i}},a^{\varepsilon_{j}g_{j}}]$, namely those for which
$i<j$ while $\varepsilon_{i}g_{i}$ is moved to the right of $\varepsilon
_{j}g_{j}$ when $\mathcal{S}$ is re-arranged to $\mathcal{S}^{\prime}$. The
result follows since%
\[
\lbrack a^{\varepsilon_{i}g_{i}},a^{\varepsilon_{j}g_{j}}]=\left\{
\begin{array}
[c]{cc}%
\lbrack a^{g_{i}},a^{g_{j}}] & \quad(\varepsilon_{i}\varepsilon_{j}=1)\\
& \\
\lbrack a^{g_{j}},a^{g_{i}}] & \quad(\varepsilon_{i}\varepsilon_{j}=-1)
\end{array}
\right.  .
\]

\end{proof}

\begin{corollary}
\label{psi}In the notation of Lemma \ref{lem1}, if $aB=\overline{c}\mu$ with
$\mu\in S$ then%
\[
\psi(a)=\mu^{2}\psi(c).
\]

\end{corollary}

Now fix $x_{1},\ldots,x_{m}\in\overline{G}$ with $\left\langle x_{1}%
,\ldots,x_{m}\right\rangle =\overline{G}$. Define mappings%
\begin{align*}
f:R^{m}  &  \rightarrow A\\
\mathbf{r}  &  \mapsto\prod_{i=1}^{m}[\widetilde{c}r_{i},x_{i}],
\end{align*}
and%
\begin{align*}
B:R^{m}\times R^{m}  &  \rightarrow A^{\prime}\\
(\mathbf{r},\mathbf{s})  &  \mapsto\prod_{i=1}^{m}[[\widetilde{c}r_{i}%
,x_{i}],\widetilde{c}s_{i}]\cdot\prod_{1\leq i<j\leq m}[[\widetilde{c}%
r_{i},x_{i}],[\widetilde{c}s_{j},x_{j}]]
\end{align*}
and%
\begin{align*}
\Xi:R^{m}  &  \rightarrow R(\overline{G}-1)\\
\mathbf{r}  &  \mapsto\sum_{i=1}^{m}r_{i}(x_{i}-1).
\end{align*}
(Here, $\mathbf{r}=(r_{1},\ldots,r_{m})$ etc.)

The following observations are more or less immediate; note that identifying
$A^{\prime}$ with $\wedge_{R}^{2}M$ we can equally well write%
\[
B(\mathbf{r},\mathbf{s})=\sum_{i}\overline{c}r_{i}(x_{i}-1)\wedge\overline
{c}s_{i}+\sum_{i<j}\overline{c}r_{i}(x_{i}-1)\wedge\overline{c}s_{j}%
(x_{j}-1),
\]
and that%
\[
\widetilde{f(\mathbf{r})}=\widetilde{c}\Xi(\mathbf{r})
\]
for each $\mathbf{r}\in R^{m}$.

\begin{lemma}
\label{lem2}\emph{(i)}%
\[
f(\mathbf{r}+\mathbf{s})=f(\mathbf{r})\cdot f(\mathbf{s})\cdot B(\mathbf{r}%
,\mathbf{s}).
\]

\emph{(ii)} $B$ is $S$-bilinear and $B(R^{m},I^{(m)})=B(I^{(m)},R^{m})=0$.

\emph{(iii)} $\Xi$ is a left $R$-module epimorphism.
\end{lemma}

Now since $A=B\left\langle c^{G}\right\rangle $ and $A^{\prime}\neq1$ there
exists $d\in A$ such that $[c,d]\neq1.$ We fix such a $d$.

\begin{lemma}
\label{ann(d-1)}If $s\in S$ and $s(\overline{d}-1)=0$ in $R$ then $s\in I$.
\end{lemma}

\begin{proof}
Say $s=\sum_{j=1}^{n}\varepsilon_{j}g_{j}$. Then $\sum_{j=1}^{n}%
\varepsilon_{j}g_{j}-\sum_{j=1}^{n}\varepsilon_{j}g_{j}\overline{d}=0$, and
Corollary \ref{psi} applies to the mapping $\psi$ given by%
\[
\psi(a)=[\widetilde{a}s,\overline{d}]\,\,(a\in A).
\]
Hence if $\widetilde{a}=\widetilde{c}\mu$, where $\mu\in S$, then%
\[
\psi(a)=\mu^{2}\psi(c)=\mu^{2}[\widetilde{c}s,\overline{d}]=\mu^{2}s[c,d].
\]
On the other hand, we also have%
\[
\psi(a)=[\widetilde{a}s,\overline{d}]=[\widetilde{c}\mu s,\overline{d}]=\mu
s[c,d].
\]
Since $\left|  k\right|  >2$ we may choose $\mu\in S$ so that $\mu-\mu
^{2}\not \equiv 0\,(\operatorname{mod}I),$ and deduce that $s[c,d]=0$. The
result follows since $I$ is the annihilator of each non-zero element in
$A^{\prime}$.
\end{proof}

\bigskip

Put%
\[
V=\Xi^{-1}(S(\overline{d}-1)).
\]
Thus $V$ is a left $S$-submodule of $R^{m}$, and $\Xi$ maps $V$ onto
$S(\overline{d}-1).$ Moreover,%
\[
f(V)\subseteq A^{\prime},
\]
since if $\Xi(v)=s(\overline{d}-1)$ then
\[
\widetilde{f(v)}=\widetilde{c}s(\overline{d}-1)=0.
\]

In view of Lemma \ref{ann(d-1)}, there is a well-defined mapping%
\[
\alpha:V\rightarrow A^{\prime}%
\]
such that%
\[
\alpha(v)=[\widetilde{c}s,d]=s[c,d]
\]
when $\Xi(v)=s(\overline{d}-1)$, $\,s\in S$. Evidently $\alpha$ is a left
$S$-module epimorphism. Define%
\[
\Phi:V\rightarrow A^{\prime}%
\]
by%
\[
\Phi(v)=f(v)-\alpha(v).
\]

\begin{lemma}
\label{quad}For each $v\in V$ and $\mu\in S$ we have%
\[
\Phi(\mu v)=\mu^{2}\Phi(v).
\]

\end{lemma}

\begin{proof}
Say $v=(r_{1},\ldots,r_{m})\in V$ and $\Xi(v)=s(\overline{d}-1)$ with $s\in
S$. For $a\in A$ put%
\[
\psi(a)=\prod_{i=1}^{m}[\widetilde{a}r_{i},x_{i}]\cdot\lbrack\widetilde
{a}s,d]^{-1}.
\]
Since $\sum r_{i}(x_{i}-1)-s(\overline{d}-1)=0,$ we may apply Corollary
\ref{psi} to deduce that if $\widetilde{a}=\widetilde{c}\lambda$ where
$\lambda\in S$ then $\psi(a)=\lambda^{2}\psi(c)$. But%
\[
\prod_{i=1}^{m}[\widetilde{a}r_{i},x_{i}]=\prod_{i=1}^{m}[\widetilde{c}\lambda
r_{i},x_{i}]=f(\lambda v)
\]
and%
\[
\lbrack\widetilde{a}s,d]=[\widetilde{c}\lambda s,d]=\alpha(\lambda v),
\]
so $\psi(a)=\Phi(\lambda v)$. Thus in particular $\Phi(v)=\psi(c)$ and
$\Phi(\mu v)=\psi(a)$ where $\widetilde{a}=\widetilde{c}\mu$, and the lemma follows.
\end{proof}

\bigskip

For $w\in R^{m}$ and $v\in V$ put%
\[
\alpha^{w}(v)=\alpha(v)+B(w,v).
\]
Then $\alpha^{w}:V\rightarrow A^{\prime}$ is a left $S$-module homomorphism,
and Lemma \ref{lem2} shows that%
\begin{align*}
\Phi(v)+\alpha^{w}(v)  &  =f(v)+B(w,v)\\
&  =f(w+v)f(w)^{-1}%
\end{align*}
for each $v\in V$.

It follows from Lemma \ref{lem2} that for $u,\,v\in V,$%
\[
\Phi(u+v)=\Phi(u)+\Phi(v)+B(u,v).
\]
With Lemma \ref{quad} this implies that $\Phi$ factors through $V\rightarrow
V/IV$, and that the mapping $\overline{\Phi}:V/IV\rightarrow A^{\prime}$
induced by $\Phi$ is quadratic as a map of $k$-vector spaces; similarly,
$\alpha^{w}$ factors through $V\rightarrow V/IV$ and induces a $k$-linear map
$\overline{\alpha}^{w}:V/IV\rightarrow A^{\prime}$.

\begin{lemma}
\label{factoring}Each of the maps $\Phi$ and $\alpha^{w}$ factors through
$V\rightarrow(V+I^{(m)})/I^{(m)}$; and $B$ factors through $R^{m}\times
R^{m}\rightarrow R^{m}/I^{(m)}\times R^{m}/I^{(m)}$.
\end{lemma}

\begin{proof}
The claim regarding $B$ is immediate from Lemma \ref{lem2}(ii). For the rest,
we separate two cases.

\medskip

\textbf{Case 1: }$p\neq2$. Let $v\in V$. Then%
\[
4\Phi(v)=\Phi(2v)=2\Phi(v)+B(v,v)
\]
so $\Phi(v)=\frac{1}{2}B(v,v)$ which depends only on the coset of $v$ modulo
$I^{(m)}$. In particular, $\Phi(I^{(m)}\cap V)=0,$ so $f(I^{(m)}\cap
V)=\alpha(I^{(m)}\cap V)$ is a $k$-subspace of $A^{\prime}$. But if
$\mathbf{s}\in I^{(m)}\cap V$ then%
\[
f(\mathbf{s})=\prod[\widetilde{c}s_{i},x_{i}]\in\lbrack B,G]\cap A^{\prime}.
\]
As $[B,G]\cap A^{\prime}<A^{\prime}$ and $A^{\prime}$ is a $1$-dimensional
$k$-space it follows that $\alpha(I^{(m)}\cap V)=f(I^{(m)}\cap V)=0$. Since
$B(w,I^{(m)})=0$ we also have $\alpha^{w}(I^{(m)}\cap V)=0$.

Thus both $\Phi$ and $\alpha^{w}$ factor through $V\rightarrow(V+I^{(m)}%
)/I^{(m)}$.

\medskip

\textbf{Case 2: }$p=2$. Consider the $k$-vector space $W=(I^{(m)}\cap V)/IV$,
and write%
\[
\overline{f}=\left(  \overline{\Phi}+\overline{\alpha}\right)  _{\left|
W\right.  }:W\rightarrow A^{\prime}.
\]
Suppose that $\overline{f}(W)\neq\{0\}$. According to Lemma 5.2 of \cite{Sg},
$\overline{f}(W)$ is then an additive subgroup of index at most $2$ in
$A^{\prime}$. However, $\overline{f}(W)=f(I^{(m)}\cap V)\subseteq\lbrack
B,G]\cap A^{\prime}$, as observed above, so $\left|  A^{\prime}:[B,G]\cap
A^{\prime}\right|  \leq2.$ This contradicts our original hypothesis; it
follows that
\[
f(I^{(m)}\cap V)=\overline{f}(W)=\{0\}
\]
and hence that $\overline{\alpha}_{\left|  W\right.  }=\overline{\Phi
}_{\left|  W\right.  }$. Since $\overline{\alpha}$ is linear, $\overline{\Phi
}$ is quadratic and $\left|  k\right|  >2$ it follows that $\overline{\alpha
}_{\left|  W\right.  }=\overline{\Phi}_{\left|  W\right.  }=0$.

Hence $\alpha(I^{(m)}\cap V)=\Phi(I^{(m)}\cap V)=0$, and the proof is now
completed as in Case 1.
\end{proof}

\begin{lemma}
\label{dim}$\dim_{k}((V+I^{(m)})/I^{(m)})\geq(m-d)\dim_{k}(M).$
\end{lemma}

\begin{proof}
Put $h=\dim_{k}(M)$. Then $R/I$ is generated as a left $S$-module by $h$
elements, one of which may be taken to be $1_{R}$; as $I\subseteq S$ it
follows that $R$ is an $h$-generator left $S$-module. Since $R(\overline
{G}-1)$ is a $d$-generator left $R$-module, it is a $dh$-generator left
$S$-module. As $V\supseteq\ker\Xi$ it follows that $R^{m}/V$ is a
$dh$-generator left $S$-module, and hence that $R^{m}/(V+I^{(m)})$ is a
$k$-vector space of dimension at most $dh$. On the other hand, $R^{m}%
/I^{(m)}\cong M^{(m)}$ is a $k$-vector space of dimension $mh$. The lemma
follows since%
\[
\dim_{k}((V+I^{(m)})/I^{(m)})=\dim_{k}(R^{m}/I^{(m)})-\dim_{k}(R^{m}%
/(V+I^{(m)})).
\]

\end{proof}

\bigskip

We can now complete the proof of Proposition \ref{quadmap}. Fix $\mathbf{z}\in
A^{(m)}$. Since $\widetilde{A}=\widetilde{B}+\widetilde{c}R$ we can write
$\widetilde{z}_{i}=\widetilde{b}_{i}+\widetilde{c}t_{i}$ with $b_{i}\in B$ and
$t_{i}\in R,$ and we put%
\[
w=(t_{1},\ldots,t_{m})\in R^{m}.
\]

The map%
\[
(r_{1},\ldots,r_{m})\mapsto(\overline{c}r_{1},\ldots,\overline{c}r_{m})
\]
induces a left $S$-module isomorphism $\theta:R^{m}/I^{(m)}\rightarrow
M^{(m)}$. Put%
\[
U=\theta((V+I^{(m)})/I^{(m)})\leq M^{(m)}.
\]
Then $U$ is a $k$-subspace of $M^{(m)}$ and Lemma \ref{dim} shows that
$\mathrm{dim}_{k}U\geq(m-d)\mathrm{dim}_{k}M$. According to Lemma
\ref{factoring}, the maps $\Phi$ and $\alpha^{w}$ induce maps $\widetilde
{\Phi}$ and $\widetilde{\alpha}^{w}$ from $(V+I^{(m)})/I^{(m)}$ to $A^{\prime
}$. Then%
\[
\Phi_{0}=\widetilde{\Phi}\circ\theta^{-1}:U\rightarrow A^{\prime}%
\]
is quadratic over $k$ and%
\[
\alpha_{0}^{w}=\widetilde{\alpha}^{w}\circ\theta^{-1}:U\rightarrow A^{\prime}%
\]
is linear over $k$; also $\alpha_{0}^{\mathbf{0}}$ is surjective; these all
follow from the corresponding properties of $\overline{\Phi},\,\overline
{\alpha}^{w}:V/IV\rightarrow A^{\prime}$ established above.

Let $u=\theta(\mathbf{r})\in U$, and let $a_{i}\in A$ be such that
$\widetilde{a}_{i}=\widetilde{c}r_{i}$. Then $(a_{1}B,\ldots,a_{m}B)=u.$ Now%
\begin{align}
\prod_{i=1}^{m}[z_{i}a_{i},x_{i}]  &  =\prod_{i=1}^{m}[\widetilde{b}%
_{i}+\widetilde{c}(t_{i}+r_{i}),x_{i}]\nonumber\\
&  =\prod_{i=1}^{m}[\widetilde{c}(t_{i}+r_{i}),x_{i}]\cdot\prod_{i=1}%
^{m}[\widetilde{b}_{i},x_{i}] \label{factor1}%
\end{align}
since $[B,A]=1,$ and similarly%
\begin{equation}
\prod_{i=1}^{m}[z_{i},x_{i}]=\prod_{i=1}^{m}[\widetilde{c}t_{i},x_{i}%
]\cdot\prod_{i=1}^{m}[\widetilde{b}_{i},x_{i}]. \label{factor2}%
\end{equation}
On the other hand,%
\begin{align*}
\Phi_{0}(u)+\alpha_{0}^{w}(u)  &  =\Phi(\mathbf{r})+\alpha^{w}(\mathbf{r})\\
&  =f(\mathbf{t}+\mathbf{r})f(\mathbf{t})^{-1}\\
&  =\prod_{i=1}^{m}[\widetilde{c}(t_{i}+r_{i}),x_{i}]\cdot\left(  \prod
_{i=1}^{m}[\widetilde{c}t_{i},x_{i}]\right)  ^{-1}.
\end{align*}
With (\ref{factor1}) and (\ref{factor2}) this shows that%
\[
\Phi_{0}(u)+\alpha_{0}^{w}(u)=\left(  \prod_{i=1}^{m}[z_{i}a_{i}%
,x_{i}]\right)  \cdot\left(  \prod_{i=1}^{m}[z_{i},x_{i}]\right)  ^{-1}.
\]
This is precisely claim (ii) of Proposition \ref{quadmap}, if we write $\Phi$
for $\Phi_{0}$ and $\alpha^{\mathbf{z}}$ for $\alpha_{0}^{w}$. In the special
case $\mathbf{z}=(1,\ldots,1)$ we can take $\mathbf{t}=\mathbf{0}$ to ensure
that $\alpha_{0}^{w}$ is surjective. This completes the proof.

\section{The second inequality, soluble case}

We are now ready to establish one of the main steps in the proof of the Key
Theorem, concerning the case where $N$ is a soluble quasi-minimal normal
subgroup of $G$. The following notation and hypotheses are in force throughout
this section.

$G$ is a finite $d$-generator group, $N$ is a soluble quasi-minimal normal
subgroup of $G$, $Z=Z_{N}$ is the maximal normal subgroup of $G$ properly
contained in $N$, and we write $M=N/Z$. We assume in addition that $[Z,G]=1$.

Recall (Lemma 4.2) that $M$ is a simple $\mathbb{F}_{p}G$-module for some
prime $p$ and that $[N,G]=N$,%
\begin{align*}
N^{p}  &  =1\,\ \ \text{if }p\neq2\\
N^{2}  &  =N^{\prime}\text{ \ and \ }N^{\prime2}=1\text{ \ if }p=2.
\end{align*}
Note that if $N^{\prime}\neq1$ then $N/Z$ cannot be cyclic so $\left|
M\right|  \geq p^{2}$.

Set%
\[
K=\left\{
\begin{array}
[c]{ccc}%
N &  & \text{if \ }N^{\prime}=1\\
&  & \\
N^{\prime} &  & \text{if }N^{\prime}>1
\end{array}
\right.  \text{.}%
\]
For $i=1,2,3$ let $x_{i1},\ldots,x_{im}$ satisfy $K\left\langle x_{i1}%
,\ldots,x_{im}\right\rangle =G$. Define%
\begin{align*}
\phi_{i}:N^{(m)}  &  \rightarrow N\\
(a_{1},\ldots,a_{m})  &  \mapsto\prod_{j=1}^{m}[a_{j},x_{ij}].
\end{align*}

\begin{proposition}
\label{sol-eqns} Let $\kappa\in K$. Then there exist $\kappa_{1},\,\kappa
_{2},\,\kappa_{3}\in N$ such that%
\[
\kappa_{1}\kappa_{2}\kappa_{3}=\kappa
\]
and for $i=1,2,3$%
\begin{equation}
\left|  \phi_{i}^{-1}(\kappa_{i})\right|  \geq\left|  N\right|  ^{m}\left|
M\right|  ^{-d-1}. \label{|fibre|}%
\end{equation}

\end{proposition}

\begin{proof}
Note that (\ref{|fibre|}) holds if (and only if) $\phi_{i}^{-1}(\kappa_{i})$
contains at least $\left|  M\right|  ^{m-d-1}$ cosets of $Z^{(m)}$. We
separate cases.\medskip

\textbf{Case 1:} where $N=K$ is abelian. Write $N$ additively, and suppose
that $G=\left\langle g_{1},\ldots,g_{d}\right\rangle $. The mapping
$(a_{1},\ldots,a_{d})\mapsto\sum a_{i}(g_{i}-1)$ induces an epimorphism from
$M^{(d)}$ to $[N,G]=N$, so $\left|  N\right|  \leq\left|  M\right|  ^{d}$.
Similarly, $\phi_{i}$ induces an epimorphism $\overline{\phi}_{i}%
:M^{(m)}\rightarrow N$. Take $\kappa_{1}=\kappa$ and $\kappa_{2}=\kappa_{3}%
=1$. Now $\phi_{i}^{-1}(\kappa_{i})$ consists of $\left|  \overline{\phi}%
_{i}^{-1}(\kappa_{i})\right|  $ cosets of $Z^{(m)}$, and the result follows
since%
\[
\left|  \overline{\phi}_{i}^{-1}(\kappa_{i})\right|  =\left|  \ker
\overline{\phi}_{i}\right|  =\left|  M^{m}\right|  /\left|  N\right|
\geq\left|  M\right|  ^{m-d}.
\]
\medskip

\textbf{Case 2: }where $\left|  N^{\prime}\right|  =2$, $K=N^{\prime}$. Let
$\overline{\phi}_{i}:M^{(m)}\rightarrow N$ and $\widetilde{\phi}_{i}%
:M^{(m)}\rightarrow N/N^{\prime}$ denote the maps naturally induced by
$\phi_{i}$. As above, each $\widetilde{\phi}_{i}$ is an epimorphism, and each
fibre of $\widetilde{\phi}_{i}$ has size at least $\left|  M\right|  ^{m-d}$.
There exists $c\in N$ with $c^{2}\neq1,$ and then $N^{\prime}=\{1,\,c^{2}\}$.
For each $i$ we now have%
\[
\overline{\phi}_{i}^{-1}(c)\cup\overline{\phi}_{i}^{-1}(c^{3})=\widetilde
{\phi}_{i}^{-1}(cN^{\prime})
\]
so $\left|  \overline{\phi}_{i}^{-1}(c^{\varepsilon(i)})\right|  \geq\frac
{1}{2}\left|  M\right|  ^{m-d}\geq\left|  M\right|  ^{m-d-1}$ where
$\varepsilon(i)$ is $1$ or $3$. One of these two values must occur at least
twice as $i$ ranges over $\{1,2,3\}$; say $\varepsilon(s)=\varepsilon
(t)=\varepsilon$. Now if $\kappa=c^{2}$ put $\kappa_{s}=\kappa_{t}%
=c^{\varepsilon},\,\kappa_{u}=1$ where $\{1,2,3\}=\{s,t,u\}$; if $\kappa=1$
put $\kappa_{1}=\kappa_{2}=\kappa_{3}=1$. In either case we then have
$\kappa_{1}\kappa_{2}\kappa_{3}=\kappa$ (since $c^{6}=c^{2}$).

Now $\phi_{i}^{-1}(\kappa_{i})$ is the union of $\left|  \overline{\phi}%
_{i}^{-1}(\kappa_{i})\right|  $ cosets of $Z^{(m)}$; so to complete the proof
in this case it remains to show that $\left|  \overline{\phi}_{i}%
^{-1}(1)\right|  \geq\left|  M\right|  ^{m-d-1}$. Put $V=\overline{\phi}%
_{i}^{-1}(N^{\prime})=\ker\widetilde{\phi}_{i}$, so $\left|  V\right|
\geq\left|  M\right|  ^{m-d}$. We claim that $\overline{\phi}_{i\left|
V\right.  }:V\rightarrow N^{\prime}$ is a quadratic form over $\mathbb{F}_{2}%
$, if $N^{\prime}$ is identified with $\mathbb{F}_{2}$. To see this, define
$B:V\times V\rightarrow N^{\prime}$ by%
\[
B(\mathbf{u},\mathbf{v})=\overline{\phi}_{i}(\mathbf{u}+\mathbf{v}%
)-\overline{\phi}_{i}(\mathbf{u})-\overline{\phi}_{i}(\mathbf{v});
\]
one readily verifies that if $\mathbf{u}=(u_{1}Z,\ldots,u_{m}Z),\,\mathbf{v}%
=(v_{1}Z,\ldots,v_{m}Z)$ then
\[
B(\mathbf{u},\mathbf{v})=\sum_{j=1}^{m}[[u_{j},x_{ij}],v_{j}]+\sum
_{j<l}[[u_{j},x_{ij}],[v_{l},x_{il}]],
\]
and hence that $B$ is bilinear as a map of $\mathbb{F}_{2}$-spaces. This
establishes the claim, which then implies that each fibre of $\overline{\phi
}_{i\left|  V\right.  }$ has size at least%
\[
2^{\dim_{\mathbb{F}_{2}}(V)-2}=\frac{1}{4}\left|  V\right|  \geq\left|
M\right|  ^{m-d-1}%
\]
(cf. Lemma 5.2 of \cite{Sg}). The result follows.

\medskip

\textbf{Case 3}: where $\left|  N^{\prime}\right|  >2$, $K=N^{\prime}$. Let
$F$ be a free group and $\pi:F\rightarrow G$ an epimorphism. Set $A=\pi
^{-1}(N)$ and $B=\pi^{-1}(Z)$. Then $A$ is free, and it is well known that the
mapping $(a,b)\mapsto\lbrack a,b]$ induces an isomorphism%
\[
\theta_{1}:A/A^{\prime}\wedge A/A^{\prime}\rightarrow A^{\prime}/[A^{\prime
},A].
\]
Write $M_{1}=A/B$. Noting that $A^{\prime}A^{p}\leq B$, one verifies easily
that $\theta_{1}$ induces an isomorphism%
\[
\theta_{2}:M_{1}\wedge M_{1}\rightarrow A^{\prime}/[B,A].
\]
The group $F/B$ acts by conjugation on $A/[B,A]$; and $\theta_{2}$ becomes an
isomorphism of $R=\mathbb{Z}(F/B)$-modules when $F/B$ is made to act
diagonally on $A/B\wedge A/B$, so $\theta_{2}$ induces an isomorphism%
\[
\theta_{3}:\wedge_{R}^{2}M_{1}=\frac{M_{1}\wedge M_{1}}{[M_{1}\wedge M_{1}%
,F]}\rightarrow\frac{A^{\prime}}{[A^{\prime},F][B,A]}.
\]

Now let $^{-}:F\rightarrow F/[A^{\prime},F][B,A]$ denote the quotient map.
Since $[N^{\prime},G][Z,N]\leq\lbrack Z,G]=1$, the map $\pi$ induces an
epimorphism $\ast:\overline{F}\rightarrow G$. Evidently%
\begin{align*}
\lbrack\overline{B},\overline{A}]  &  =[\overline{A}^{\prime},\overline
{F}]=1,\\
\ker(\ast)  &  \leq\overline{B},\,\,\,\,\overline{A}^{\ast}%
=N,\,\,\,\,\overline{B}^{\ast}=Z,
\end{align*}
and%
\[
\left|  \overline{A}^{\prime}:\overline{A}^{\prime}\cap\lbrack\overline
{B},\overline{F}]\right|  \geq\left|  N^{\prime}\right|  >2.
\]

Thus all the hypotheses of Section 6 are satisfied if we take $\overline{F}$
for $G$, $\overline{A}$ for $A$ and $\overline{B}$ for $B$; Proposition
\ref{sol-eqns} thus reduces in the present case to an application of
Proposition 6.2, with $G$ taking the role of $G^{\ast}$.

This completes the proof.
\end{proof}

\section{Word combinatorics}

In the next three sections we examine the solution of equations in a direct
product of quasisimple groups. This preparatory section is devoted to some
observations on the shape of abstract group words, generalizing Lemma 1 of
\cite{N2}: these will help us to keep track of the equations as the unknowns
are successively eliminated.

The material here is rather abstract, and won't make much sense until it is
applied. However, it seems inevitable, given the nature of our main theorems,
that at some stage we will have to get to grips with the detailed rewriting of
words in a group; by separating off in this section some of the most technical
steps, we hope to make the complicated arguments of the later sections a
little less opaque.

\bigskip

Let $\Gamma$ be a group and $Y$ a non-empty set. The \emph{free }$\Gamma
$\emph{-group on} $Y$ is the free group on the alphabet $Y^{\Gamma}%
=\{y^{g}\mid y\in Y,\,g\in\Gamma\}$, on which $\Gamma$ acts by permuting the
basis in the obvious way. We denote it by%
\[
F_{\Gamma}(Y);
\]
it may be identified with the normal closure of the free group $F(Y)$ in the
free product $F(Y)\ast\Gamma$.

A subset $Z$ of $F_{\Gamma}(Y)$ will be called \emph{independent} if every map
from $Z$ into an arbitrary $\Gamma$-group $S$ can be extended to a $\Gamma
$-equivariant homomorphism from $F$ to $S$ (thus for example every subset of
$Y$ is independent). The following `invariance' and `exchange' principles are
more or less self-evident: (i) if $Z$ is independent and $g(y)\in\Gamma$ for
each $y\in Z$ then $\{y^{g(y)}\mid y\in Z\}$ is independent; (ii) if
$Z\cup\{x\}$ is independent and $P,\,Q\in\left\langle Z^{\Gamma}\right\rangle
$ then $Z\cup\{PxQ\}$ is independent. A family of elements $\{z_{1}%
,z_{2},\ldots\}$ is called independent if its terms are all distinct and form
an independent set.

As a matter of notation, we will usually write $y$ for $y^{1}$ and $y^{-g}$ in
place of $\left(  y^{g}\right)  ^{-1}$ ($y\in Y,\,g\in\Gamma$).\bigskip

Now we fix two disjoint sets, a set $X$ of \emph{variables} and a set $P$ of
\emph{parameters}$,$ and consider the free $\Gamma$-group
\[
F=F_{\Gamma}(X\cup P).
\]
Let $\mathcal{M}$ denote the the free monoid on the set $\{y^{\pm g}\mid y\in
X\cup P,\,g\in\Gamma\}$; this is the set of `unreduced' group words on the
alphabet $X^{\Gamma}\cup P^{\Gamma}$. Let $W\subseteq\mathcal{M}$ denote the
free monoid on $X\cup X^{-1}$. There is a natural map $\overline{\phantom{m}}:
\mathcal{M}\rightarrow F$ (evaluation), and we define a mapping $\widehat
{\phantom{m}}:\mathcal{M}\rightarrow W$ as follows: for $U\in\mathcal{M}$, let
$\widehat{U}\in W$ denote the word obtained from $U$ by deleting all terms
belonging to $P^{\Gamma}\cup P^{-\Gamma}$ and replacing each term $x^{\pm g}$
with $x^{\pm1}$ ($x\in X$, $g\in\Gamma$).

For $U,\,V\in\mathcal{M}$ we write%
\[
U=_{F}V\Longleftrightarrow\overline{U}=\overline{V}%
\]
(the notation $U=V$ for $U,\,V\in\mathcal{M}$ will always mean that $U$ and
$V$ are identical as words).

We write%
\[
\left|  x\right|  =\left\{
\begin{array}
[c]{ccc}%
x &  & (x\in X)\\
&  & \\
x^{-1} &  & (x\in X^{-1})
\end{array}
\right.  ,
\]
and for $w\in W$ put
\[
\sup(w)=\{\left|  x\right|  \mid x\ \text{occurs in }w\}.
\]
We call $w\in W$ \emph{balanced} if each element of $\sup(w)\cup\sup(w)^{-1}$
occurs exactly once in $w$.

\begin{lemma}
\label{nontriv}Suppose that $w\in W$ is balanced and $w\neq_{F}1$. Then%
\begin{equation}
w=Ax^{-1}By^{-1}CxDyE \label{w=wd(1)}%
\end{equation}
for some $x,\,y\in X\cup X^{-1}$ with $\left|  x\right|  \neq\left|  y\right|
$ and $\{\left|  x\right|  ,\,\left|  y\right|  \}\cap\sup(ABCDE)=\emptyset$.
\end{lemma}

\begin{proof}
The hypotheses imply that $w=u_{1}y^{-1}vyu_{2}$ where $y\in X\cup X^{-1}$ and
$v\neq\emptyset$. Choose such an expression with $v$ as short as possible. Say
$x$ occurs in $v$ where $x\in X\cup X^{-1}$. Then $x^{-1}$ must occur in
$u_{1}$ or in $u_{2}$; in the first case we have (\ref{w=wd(1)}), in the
second case we get (\ref{w=wd(1)}) on replacing $x$ by $x^{-1}$ and then
interchanging $x$ and $y$. The final claim is clear since $w$ is balanced.
\end{proof}

\begin{proposition}
\label{1twistedcom}Let $V\in\mathcal{M}$. Suppose that $\widehat{V}$ is
balanced and $\widehat{V}\neq_{F}1$. Then%
\[
V=_{F}T_{a,b}(\xi,\eta)\cdot V_{1}%
\]
for some $a,\,b\in\Gamma$ and $\xi,\eta,\,V_{1}\in\mathcal{M}$ such that
\emph{(i) }the family $\{\overline{\xi},\,\overline{\eta}\}\cup\sup
(\widehat{V_{1}})\cup P$ is independent, and \emph{(ii)} ignoring exponents
from $\Gamma$, each term from $P\cup P^{-1}$ occurs with the same multiplicity
in $V_{1}$ as it has in $V$.
\end{proposition}

\noindent Recall that
\[
T_{a,b}(\xi,\eta)=\xi^{-1}\eta^{-1}\xi^{a}\eta^{b}.
\]

\begin{proof}
Lemma \ref{nontriv} shows that%
\[
\widehat{V}=Ax^{-1}By^{-1}CxDyE,
\]
say, for suitable words $A,\,B$ etc. in $W$ and $x,\,y\in X\cup X^{-1}$ with
$\left|  y\right|  \neq\left|  x\right|  $, $\{\left|  x\right|  ,\,\left|
y\right|  \}\cap\sup(ABCDE)=\emptyset.$ It follows that%
\[
V=A^{\prime}x^{-e}B^{\prime}y^{-f}C^{\prime}x^{ea}D^{\prime}y^{fb}E^{\prime}%
\]
where $a,\,b,\,e,\,f\in\Gamma$ and $A^{\prime},\,B^{\prime},\ldots
\in\mathcal{M}$ satisfy $\widehat{A^{\prime}}=A,\,\,\widehat{B^{\prime}}=B$
etc. Now put%
\begin{align*}
U_{1}  &  =A^{\prime ab^{-1}}D^{\prime b^{-1}},\,U_{2}=U_{1}^{a^{-1}}C^{\prime
a^{-1}},\\
\xi &  =U_{2}x^{e}A^{\prime-1},\\
\eta &  =U_{1}y^{f}B^{\prime-1}U_{2}^{-1}.
\end{align*}
A direct calculation shows that%
\[
V=_{F}T_{a,b}(\xi,\eta)\cdot V_{1}%
\]
where%
\[
V_{1}=A^{\prime ab^{-1}a^{-1}b}D^{\prime b^{-1}a^{-1}b}C^{\prime a^{-1}%
b}B^{\prime b}E^{\prime}.
\]
Note that $\widehat{V_{1}}=ADCBE$. The claim (i) follows from the invariance
and exchange principles, and the claim (ii) is clear.
\end{proof}

\bigskip

For the next proposition we need some further notation. Fix a mapping
$\chi:X\rightarrow\{1,\ldots,m\},$ and for each $x\in X$ define $\chi
(x^{-1})=-\chi(x)$. We call $\chi(x)$ the \emph{colour} of $x$. For
$w=y_{1}y_{2}\ldots y_{k}\in W$ (with each $y_{i}\in X\cup X^{-1}$) define
$\chi(w)$ to be the sequence%
\[
\chi(w)=\left(  \chi(y_{1}),\ldots,\chi(y_{k})\right)  .
\]
A new sequence $\tau(w)$, the \emph{colour type of }$w$, is now defined as
follows: first, wherever a segment consisting of consecutive equal negative
terms occurs in $\chi(w)$, delete all but one of them (so a maximal segment
$(-r,-r,\ldots,-r)$ is contracted to$\ -r$); then replace each term by its
absolute value. For example $(1,1,-2,2,-2,-2,-2,-3)\mapsto
(1,1,-2,2,-2,-3)\mapsto(1,1,2,2,2,3)$.

For sequences $S$ and $T$ we write%
\[
S\leq T
\]
to indicate that $S$ is a subsequence of $T$. Put%
\[
L_{n}=(1,2,\ldots,m,1,2,\ldots,m,\ldots,1,2,\ldots,m)
\]
with $n$ repetitions of $(1,2,\ldots,m)$.

\begin{lemma}
\label{balance}Let $w\in W$ be balanced. Put $Y=\sup(w)$, and suppose that
$\tau(w)\leq L_{n},$ where $1\leq n<\left|  Y\right|  $. Then there exist
$x,\,y\in Y\cup Y^{-1}$ with $\left|  y\right|  \neq\left|  x\right|  $ such
that%
\begin{equation}
w=Ax^{-1}By^{-1}CxDyE \label{w=word}%
\end{equation}
where $w_{0}=ADCBE$ is balanced and $\tau(w_{0})\leq L_{n}$.
\end{lemma}

\begin{proof}
We claim that $w\neq_{F}1$. The proof is by induction on $n$. Suppose that
$w=_{F}1$. Then $w=uxx^{-1}v$ where $x\in X\cup X^{-1}$. Suppose that $x\in X$
and $x$ has colour $i$. Then $\tau(w)=(\tau(u),i,i,S)$ where $\tau(v)$ is
either $S$ or $(i,S)$. In either case, $\tau(uv)\leq(\tau(u),i,S)\leq
L_{n-1}.$ One sees similarly that $\tau(uv)\leq L_{n-1}$ if $x\in X^{-1}$. As
$uv$ is balanced and $\left|  \sup(uv)\right|  =\left|  \sup(w)\right|
-1>n-1$ the inductive hypothesis gives $uv\neq_{F}1,$ a contradiction.

Applying Lemma \ref{nontriv} we obtain the expression (\ref{w=word}).

It is clear that $w_{0}$ is again balanced. To establish the final claim,
suppose for example that $\chi(x)=i>0$ and $\chi(y)=j>0$. Let $A_{0}%
,\,B_{0},\,C_{0}$ be the words obtained from $A,\,B,\,C$ respectively by
removing all terms coloured $-i$ from the end of $A$ and the beginning of $B$,
and all terms coloured $-j$ from the end of $B$ and the beginning of $C.$
Unless $B_{0}=\emptyset$ and $i=j$ we then have%
\begin{align*}
\tau(w)=(\tau(A_{0}),i,\tau(B_{0}),j,\tau(C_{0}),i,\tau(D),j,\tau(E)),\\
\tau(w_{0})\leq(\tau(A_{0}),i,\tau(D),j,\tau(C_{0}),i,\tau(B_{0}),j,\tau(E)).
\end{align*}
It is easy to see that if the first sequence is a subsequence of $L_{n}$, then
so is the second. The other cases are dealt with similarly.
\end{proof}

\begin{proposition}
\label{L_nProp}Let $V\in\mathcal{M}$ and $k\in\mathbb{N}$, and put
$w=\widehat{V}$. Suppose \emph{(a)} $w$ is balanced and \emph{(b)}
$\tau(w)\leq L_{n}$ for some $n\geq1$ with%
\[
\left|  \sup(w)\right|  \geq n+2k.
\]
Then there exist $V_{k},$\thinspace$\xi_{1},\eta_{1},\ldots,\xi_{k},\eta
_{k}\in\mathcal{M}$ such that $\left|  \sup(\widehat{V_{k}})\right|  =\left|
\sup(w)\right|  -2k$ and\newline\smallskip\emph{(i)} $\{\overline{\xi_{1}%
},\overline{\eta_{1}},\ldots,\overline{\xi_{k}},\overline{\eta_{k}}\}\cup
\sup(\widehat{V_{k}})\cup P$ is an independent family;\newline\smallskip
\emph{(ii)}%
\[
V=_{F}T_{a_{1},b_{1}}(\xi_{1},\eta_{1})\cdot\ldots\cdot T_{a_{k},b_{k}}%
(\xi_{k},\eta_{k})\cdot V_{k}%
\]
for some $a_{i},\,b_{i}\in\Gamma$.
\end{proposition}

\begin{proof}
Put $Y=\sup(w)$. Lemma \ref{balance} shows that%
\[
w=Ax^{-1}By^{-1}CxDyE,
\]
say, for suitable words $A,\,B,\ldots$ \ in $W$ and $x,\,y\in Y\cup Y^{-1}$
with $\left|  y\right|  \neq\left|  x\right|  $. We may now define $\xi
_{1}=\xi$, $\eta_{1}=\eta$, $\ a_{1}=a,\,b_{1}=b$ and $V_{1}$ as in the proof
of Proposition \ref{1twistedcom} above, to obtain%
\[
V=_{F}T_{a_{1},b_{1}}(\xi_{1},\eta_{1})\cdot V_{1},
\]
where $\widehat{V_{1}}=ADCBE$ and the family $\{\overline{\xi_{1}}%
,\overline{\eta_{1}}\}\cup\sup(\widehat{V_{1}})\cup P$ is independent.
Evidently%
\[
\left|  \sup(\widehat{V_{1}})\right|  =\left|  \sup(w)\right|  -2.
\]
If $k=1$ we are done.

Suppose that $k>1$. Put $w_{1}=\widehat{V_{1}}$. According to Lemma
\ref{balance} the word $w_{1}$ is balanced and satisfies $\tau(w_{1})\leq
L_{n}$. Also%
\begin{align*}
\left|  \sup(w_{1})\right|   &  =\left|  \sup(w)\setminus\{\left|  x\right|
,\left|  y\right|  \}\right| \\
&  \geq n+2(k-1).
\end{align*}
Arguing by induction on $k$, we may therefore suppose that%
\[
V_{1}=_{F}T_{a_{2},b_{2}}(\xi_{2},\eta_{2})\cdot\ldots\cdot T_{a_{k},b_{k}%
}(\xi_{k},\eta_{k})\cdot V_{k}%
\]
is of the required form, and the result follows.
\end{proof}

\section{Equations in semisimple groups, 1: the second inequality}

Let $N$ be a quasi-semisimple group with centre $Z$ and let $g_{1}%
,\ldots,g_{m}$ be automorphisms of $N$. We assume that $N/Z$ is the direct
product of $n\geq2$ simple groups, and that the group generated by
$g_{1},\ldots,g_{m}$ permutes these transitively.

For each $i$ let $c(g_{i})$ denote the number of cycles in this permutation
representation of $g_{i}$. We shall establish the following, where
$D\in\mathbb{N}$ is the absolute constant appearing in Theorem 1.9:

\begin{proposition}
\label{ssB}Suppose that%
\begin{equation}
\sum_{i=1}^{m}c(g_{i})\leq(m-2)n-2D. \label{c(g)-cond}%
\end{equation}
Then for each $\kappa\in N$ the number of solutions $\mathbf{u}=(u_{1}%
,...,u_{m})\in N^{(m)}$ to the equation
\begin{equation}
\kappa=[\mathbf{u},\mathbf{g}]:=\prod_{i=1}^{m}[u_{i},g_{i}] \label{sseq}%
\end{equation}
is at least $|N|^{m}|N/Z|^{-4D}$.
\end{proposition}

Before proving this, let us deduce the version used in \S 4 for the proof of
the Key Theorem:

\begin{proposition}
\label{prop92} Let $G$ be a finite group and $N$ a quasi-semisimple
quasi-minimal normal subgroup of $G$, such that $N/Z_{N}$ is not simple.
Suppose that $G=\left\langle y_{1},\ldots,y_{m}\right\rangle N$, and that the
$m$-tuple $\mathbf{y}$ has the $(k,\varepsilon)$ fixed-point property where
$k\varepsilon\geq2D+4$. Define $\phi:N^{(m)}\rightarrow N$ by%
\[
\phi(\mathbf{a})=\prod_{i=1}^{m}[a_{i},y_{i}].
\]
Then for each $\kappa\in N$ we have%
\[
\left|  \phi^{-1}(\kappa)\right|  \geq\left|  N\right|  ^{m}\left|
N/Z_{N}\right|  ^{-4D}.
\]

\end{proposition}

\noindent The various terms used in this statement were introduced in Section
4. Rather than repeating the definitions wholesale, we recall those
consequences that are relevant here: these may be taken as the hypotheses for
Proposition \ref{prop92}.

\begin{itemize}
\item The normal subgroup $N$ satisfies $N=[N,N]>Z_{N}=\mathrm{Z}(N),$ and
$N/Z_{N}=T_{1}\times\cdots\times T_{n}$ where $n\geq2$ and the $T_{i}$ are
isomorphic simple groups;

\item The conjugation action of $G$ permutes the set $\{T_{1},\ldots,T_{n}\}$
transitively, and at least $k$ of the $y_{j}$ move at least $\varepsilon n$ of
the $T_{i}$.
\end{itemize}

\begin{proof}
Apply Proposition \ref{ssB}, taking $g_{i}$ to be the image of $y_{i}$ in
$\mathrm{Aut}(N)$. It is only necessary to verify the condition
(\ref{c(g)-cond}). Let $\left|  \mathrm{fix}(g_{j})\right|  $ denote the
number of fixed points of $g_{j}$ in the set $\{T_{1},\ldots,T_{n}\}$. Then%
\[
n\geq\left|  \mathrm{fix}(g_{j})\right|  +2(c(g_{j})-\left|  \mathrm{fix}%
(g_{j})\right|  ),
\]
and for at least $k$ values of $j$ we have $\left|  \mathrm{fix}%
(g_{j})\right|  \leq n-\varepsilon n$. Therefore%
\[
2\sum_{i=1}^{m}c(g_{i})\leq\sum_{i=1}^{m}\left(  n+\left|  \mathrm{fix}%
(g_{i})\right|  \right)  \leq k(2n-\varepsilon n)+(m-k)\cdot2n
\]
and (\ref{c(g)-cond}) follows since $k\varepsilon\geq2D+4$ and $n\geq2$.
\end{proof}

\bigskip

We proceed to prove Proposition \ref{ssB}, and from now on write
$G=\left\langle g_{1},\ldots,g_{m}\right\rangle $. The universal cover of $N$
is a direct product $\widetilde{N}=S_{1}\times\cdots\times S_{n}$ where each
$S_{i}$ is a quasisimple group of universal type and $n\geq2$. The action of
$G$ on $N$ lifts to an action on $\widetilde{N}$, and $G$ permutes
$\{S_{1},\ldots,S_{n}\}$ the same way it permutes the simple factors of $N/Z$.

Now $N=\widetilde{N}/A$ for some $A\leq\widetilde{Z}=\mathrm{Z}(\widetilde
{N})$. If the proposition holds with $\widetilde{N}$ in place of $N$, and
$\widetilde{\kappa}$ is a preimage of $\kappa$, then $[\widetilde{\mathbf{u}%
},\mathbf{g}]=\widetilde{\kappa}$ holds for at least $|\widetilde{N}%
|^{m}|\widetilde{N}/\widetilde{Z}|^{-4D}$ values of $\widetilde{\mathbf{u}}%
\in\widetilde{N}^{(m)}$. These project to at least%
\[
\frac{|\widetilde{N}|^{m}|\widetilde{N}/\widetilde{Z}|^{-4D}}{\left|
A\right|  ^{m}}=|N|^{m}|N/Z|^{-4D}%
\]
solutions $\mathbf{u}$ of (\ref{sseq}) in $N^{(m)}$. Thus we may, and shall,
assume henceforth that $N=\widetilde{N}$.

By way of notation we shall write%
\[
S_{i}^{g^{-1}}=S_{^{g}i}\qquad\,(g\in G,\,1\leq i\leq n)
\]
(so $i\mapsto\,^{g}i$ gives the left action of $G$ on $\{1,\ldots,n\}$
corresponding to its right action on $\{S_{1},\ldots,S_{n}\}$). Since the
action is transitive, the groups $S_{i}$ are all isomorphic; we fix an
identification of each $S_{i}$ with a fixed quasisimple group $S$. Thus
elements of $N$ will be written in the form%
\[
x=(x(1),x(2),\ldots,x(n))
\]
with each $x(i)\in S$, and the action of $G$ takes the form%
\[
x^{g}=(x(^{g}1)^{g(1)},x(^{g}2)^{g(2)},\ldots,x(^{g}n)^{g(n)});
\]
here $g(i)\in\mathrm{Aut}(S)$ is induced by $g_{\left|  S_{\left(
^{g}i\right)  }\right.  }:S_{\left(  ^{g}i\right)  }\rightarrow S_{i}$ (when
each $S_{j}$ is identified with $S$).

For a subset $\Delta$ of $\{1,\ldots,n\}$,%
\[
\pi_{\Delta}:N\rightarrow\prod_{i\in\Delta}S_{i}%
\]
will denote the natural projection.

\bigskip

We are going to think of (\ref{sseq}) as the equation%
\begin{equation}
\kappa=x_{1}x_{2}\ldots x_{m}, \label{kapp=x}%
\end{equation}
to be solved for $x_{1},\ldots,x_{n}\in N$ subject to the conditions%
\begin{equation}
x_{i}\in\lbrack N,g_{i}] \label{x in [N,gi]}%
\end{equation}
for each $i.$ The equation (\ref{kapp=x}) is equivalent to the system of
equations $\mathcal{F}=(F_{1},\ldots,F_{n}):$%
\begin{equation}
\kappa(s)=x_{1}(s)x_{2}(s)\ldots x_{m}(s). \tag{$F_s$}\label{Fs}%
\end{equation}

To analyse the condition (\ref{x in [N,gi]}), let $\Omega_{i}$ denote the set
of orbits of $g_{i}$ on the set $\{1,\ldots,n\}$. Then $x\in\lbrack N,g_{i}]$
if and only if $\pi_{\Delta}(x)\in\lbrack\pi_{\Delta}(N),g_{i}]$ for each
orbit $\Delta\in\Omega_{i}$. For each such $\Delta$ let $k_{\Delta}$ be the
first member of $\Delta$ and put $n(\Delta)=\left|  \Delta\right|  $. Then
$g_{i}^{n(\Delta)}$ maps $S_{k_{\Delta}}$ to itself, inducing the automorphism
$\beta_{i}(\Delta)=g^{n(\Delta)}(k_{\Delta})$ of $S.$

We claim that $\pi_{\Delta}(x_{i})\in\lbrack\pi_{\Delta}(N),g_{i}]$ if and
only if there exists $u_{i}(\Delta)\in S$ such that%
\begin{equation}
\prod_{j=0}^{n(\Delta)-1}x_{i}(^{g_{i}^{j}}k_{\Delta})^{g_{i}^{j}(k_{\Delta}%
)}=u_{i}(\Delta)^{-1}u_{i}(\Delta)^{\beta_{i}(\Delta)}. \tag{$H_i,_\Delta
$}\label{Hidelta}%
\end{equation}
Indeed, dropping the subscript $i$ for the moment and putting $k=k_{\Delta}$,
$n=n(\Delta)$, if $\pi_{\Delta}(x)=[\pi_{\Delta}(u),g_{i}]$ then%
\begin{align}
x(k)  &  =u(k)^{-1}u(^{g}k)^{g(k)}\nonumber\\
x(^{g}k)  &  =u(^{g}k)^{-1}u(^{g^{2}}k)^{g(^{g}k)}\nonumber\\
&  \vdots\label{steps}\\
x(^{g^{n-1}}k)  &  =u(^{g^{n-1}}k)^{-1}u(^{g^{n}}k)^{g(^{g^{n-1}}k)}\nonumber
\end{align}
and (\ref{Hidelta}) follows with $u_{i}(\Delta)=u(k)$ (note that%
\[
g(^{g^{r-1}}k)\ldots g(^{g}k)g(k)=g^{r}(k)
\]
for each $r$). Conversely, if (\ref{Hidelta}) holds then putting
$u(k)=u_{i}(\Delta)$ we can solve (\ref{steps}) for $u(^{g}k),\,u(^{g^{2}%
}k),\ldots,u(^{g^{n-1}}k)$ in turn and so determine $\pi_{\Delta}(u)$ with
$\pi_{\Delta}(x)=[\pi_{\Delta}(u),g_{i}]=\pi_{\Delta}(x)$. This establishes
the claim.

\bigskip

Put%
\begin{align*}
X  &  =\{x_{i}(s)\mid1\leq i\leq m,\ 1\leq s\leq n\},\\
\mathcal{U}  &  =\{u_{i}(\Delta)\mid1\leq i\leq m,\ \Delta\in\Omega_{i}\},\\
\mathcal{K}  &  =\{\kappa(1),\kappa(2),\ldots,\kappa(n)\},\\
P  &  =\mathcal{U}\cup\mathcal{K}.
\end{align*}
Note that%
\[
\left|  \mathcal{U}\right|  =\sum_{i=1}^{m}\left|  \Omega_{i}\right|
=\sum_{i=1}^{m}c(g_{i}).
\]
We start by considering these as sets of abstract symbols, and call $P$ the
set of \emph{parameters} and $X$ the set of \emph{variables}. Each term
$x_{i}(s)$ is assigned the \emph{colour} $i$. We shall apply the results of
Section 8, taking $\Gamma=\mathrm{Aut}(S)$ and $F=F_{\Gamma}(X\cup P)$.

We are going to reduce the system $\mathcal{F}$ subject to the conditions
(\ref{x in [N,gi]}) to a single equation. First of all, for each $i$ and each
$\Delta\in\Omega_{i}$ we solve equation (\ref{Hidelta}) for $x_{i}(k_{\Delta
})$ and substitute the resulting expression in equation ($F_{k_{\Delta}}$).
That is, replace $x_{i}(k_{\Delta})$ by%
\[
u_{i}(\Delta)^{-1}u_{i}(\Delta)^{\beta_{i}(\Delta)}\cdot\left(  \prod
_{j=1}^{n(\Delta)-1}x_{i}(^{g_{i}^{j}}k_{\Delta})^{g_{i}^{j}(k_{\Delta}%
)}\right)  ^{-1}.
\]
At this stage, the conditions (\ref{x in [N,gi]}) and all the variables
$x_{i}(k_{\Delta})$ have been eliminated, at the cost of introducing some
parameters from $\mathcal{U}$. Call the resulting system of equations
$\mathcal{F}^{\prime}=(F_{1}^{\prime},\ldots,F_{n}^{\prime})$, and let $U_{s}$
be the word on $X^{\Gamma}\cup P^{\Gamma}$ on the right-hand side of
$F_{s}^{\prime}$.

Together, the words $\widehat{U_{1}},\ldots,\widehat{U_{n}}$ contain the
variables%
\[
x_{i}(s),x_{i}(s)^{-1}\,\,\,\,(s\neq k_{\Delta}\text{ for }\Delta\in\Omega
_{i}),
\]
that is,%
\[
mn-\sum_{i=1}^{m}\left|  \Omega_{i}\right|  =mn-\sum_{i=1}^{m}c(g_{i})
\]
matching pairs $x,\,x^{-1}$.

Recalling the definition of \emph{colour type} from the previous section,
observe also that the colour type of each $\widehat{U_{s}}$ satisfies%
\[
\tau(\widehat{U_{s}})\leq L_{1}=(1,\ldots,m).
\]

Next, we successively reduce the number of equations by a process of
substitution of variables. Suppose that $x\in X\cup X^{-1}$ occurs in $U_{1}$
but $x^{-1}$ does not; then $x^{-1}$ appears in $U_{l}$ for some $l\neq1$.
Solve $F_{l}^{\prime}$ for $x$ and substitute the resulting expression in
$U_{1}$. We call this a substitution $(l\rightarrow1)$. Each such operation
reduces by one both the number of equations in $\mathcal{F}^{\prime}$ and the
total number of variables. We claim now that it is possible to apply $n-1$
substitutions and thus reach an equivalent system consisting of the single
equation
\begin{equation}
\kappa(1)=U, \label{one2}%
\end{equation}
where $U$ is a certain word on $X^{\Gamma}\cup P^{\Gamma}$.

To establish the claim, let us call two equations $F_{s}^{\prime}%
,\,F_{t}^{\prime}$ \emph{linked} if they share a variable from $X$ (which then
must appear with positive exponent in one of them and negative exponent in the
other), and let $\mathcal{R}$ be the equivalence relation on $\mathcal{F}%
^{\prime}$ generated by the linked pairs. Now $F_{s}^{\prime}$ and
$\,F_{t}^{\prime}$ are linked precisely when $s$ and $t$ lie in the same orbit
of $g_{i}$ for some $i$. As $G=\left\langle g_{1},\ldots,g_{m}\right\rangle $
acts transitively on $\{1,\ldots,n\}$ it follows that $\mathcal{F}^{\prime}$
consists of one equivalence class under $\mathcal{R}$. Now a substitution
$(l\rightarrow1)$ eliminates only the variable used to link $F_{l}^{\prime}$
with $F_{1}^{\prime}$, and simultaneously eliminates the equation
$F_{l}^{\prime}$; so the resulting system of equations $\mathcal{F}%
^{^{\prime\prime}}$ still consists of a single $\mathcal{R}$-equivalence
class. If $\left|  \mathcal{F}^{^{\prime\prime}}\right|  >1$ there exists
$l^{\prime}\neq1$ such that the new equation $F_{1}^{^{\prime\prime}}$ is
linked to $F_{l^{\prime}}^{\prime}\in\mathcal{F}^{^{\prime\prime}}$, and we
can perform a substitution $(l^{\prime}\rightarrow1)$. Evidently the process
may be repeated as long as more than one equation remains in the system, and
the claim is now clear.

Since each substitution eliminates precisely one pair $x,\,x^{-1}$ ($x\in X$),
the word $\widehat{U}$ is balanced and%
\[
\left\vert \sup(\widehat{U})\right\vert =mn-\sum_{i=1}^{m}c(g_{i})-(n-1)\geq
n+2D+1.
\]
Moreover, we claim that $\tau(\widehat{U})\leq L_{n}.$ To see this, suppose
that after $f-1$ substitutions $U_{1}$ has been transformed into $V_{f}$ where
$\tau(\widehat{V_{f}})\leq L_{f}$. The next substitution $(l\rightarrow1)$ has
one of the following effects on $\widehat{V_{f}}$:%
\begin{align*}
\widehat{V_{f}}=A_{0}\alpha x^{-1}\beta B_{0} &  \mapsto\widehat{V_{f+1}%
}=A_{0}\alpha\cdot DC\cdot\beta B_{0}\\
\widehat{V_{f}}=A_{0}\alpha x\beta B_{0} &  \mapsto\widehat{V_{f+1}}%
=A_{0}\alpha\cdot DC\cdot\beta B_{0}%
\end{align*}
where $\tau(\alpha x^{-1}\beta)\leq L_{1}$ and $\ \widehat{U_{l}}=CxD$ in the
first case, $\tau(\alpha x\beta)\leq L_{1}$ and$\ \ \widehat{U_{l}}=Cx^{-1}D$
in the second case. Say $\chi(x)=i$. Since $\tau(\widehat{U_{l}})\leq L_{1}$,
in the first case we have
\begin{align*}
\tau(C) &  \leq(1,\ldots,i-1),\,\,\quad\tau(D)\leq(i+1,\ldots,m),\\
\tau(\alpha) &  \leq(1,\ldots,i),\,\,\quad\tau(\beta)\leq(i,\ldots,m),
\end{align*}
while in the second case
\begin{align*}
\tau(C) &  \leq(1,\ldots,i),\,\,\quad\tau(D)\leq(i,\ldots,m),\\
\tau(\alpha) &  \leq(1,\ldots,i-1),\,\,\quad\tau(\beta)\leq(i+1,\ldots,m).
\end{align*}
In either case, $\tau(\alpha\cdot DC\cdot\beta)\leq L_{2}$. Therefore
$\tau(\widehat{V_{f+1}})\leq L_{f+1}$, and the claim follows by induction.

We may now apply Proposition 8.4, which shows that%
\[
U=_{F}T_{a_{1},b_{1}}(x_{1},y_{1})\cdot\ldots\cdot T_{a_{D},b_{D}}(x_{D}%
,y_{D})\cdot U_{0}%
\]
where $a_{i},\,b_{i}\in\Gamma$ and
\begin{equation}
\left\{  \overline{x_{1}},\,\overline{y_{1}},\ldots,\overline{x_{D}%
},\,\overline{y_{D}}\right\}  \cup\sup(\widehat{U_{0}})\cup\mathcal{U}%
\cup\mathcal{K}\label{indep}%
\end{equation}
is an independent family. Moreover, $|\sup(\widehat{U_{0}})|=|\sup(\widehat
{U})|-2D$ so putting $\mathcal{X}_{0}=\sup(\widehat{U_{0}})\cup\mathcal{U}$ we
have%
\[
\left\vert \mathcal{X}_{0}\right\vert =mn-\sum_{i=1}^{m}c(g_{i}%
)-(n-1)-2D+\left\vert \mathcal{U}\right\vert =mn-(n-1)-2D.
\]

Define $\psi:\mathcal{K}\rightarrow S$ by $\kappa(j)\mapsto\kappa(j)$, and
extend $\psi$ arbitrarily to $\mathcal{K}\cup\mathcal{X}_{0}$. Let $\mu\in S$
be the value of $U_{0}$ determined by $\psi$. According to Theorem 1.9, there
exist $\xi_{1},\,\eta_{1},\ldots,\xi_{D},\,\eta_{D}\in S$ such that%
\[
T_{a_{1},b_{1}}(\xi_{1},\eta_{1})\cdot\ldots\cdot T_{a_{D},b_{D}}(\xi_{D}%
,\eta_{D})=\kappa(1)\mu^{-1}.
\]
Since (\ref{indep}) is an independent family, we can extend $\psi$ to a
$\Gamma$-equivariant homomorphism $F\rightarrow S$ sending $x_{i}$ to $\xi
_{i}$ and $y_{i}$ to $\eta_{i}$ for each $i,$ and then $\psi(U)=\kappa(1)$.

Each such mapping $\psi$ thus gives rise to a solution of the original
equation (\ref{sseq}). Distinct mappings give distinct solutions, because the
values of all the variables $x_{i}(s)$ are determined by the values of the
$u_{i}(s)$ via (\ref{steps}). The number of solutions is therefore at least
equal to the number of possible maps $\psi,$ which is at least%
\[
\left|  S\right|  ^{\left|  \mathcal{X}_{0}\right|  }\geq\left|  S\right|
^{mn-n-2D+1}=\frac{\left|  N\right|  ^{m}}{\left|  S\right|  ^{n+2D-1}}.
\]
Now $\left|  S\right|  <\left|  S/\mathrm{Z}(S)\right|  ^{2}$ (\cite{GLS},
\S 6.1) so
\[
\left|  S\right|  ^{n+2D-1}<\left|  S\right|  ^{2nD}<\left|  S/\mathrm{Z}%
(S)\right|  ^{4nD}=\left|  N/Z\right|  ^{4D},
\]

and the proposition follows.

\section{Equations in semisimple groups, 2: powers}

Fix a positive integer $q$. The constants $D$, $C=C(q)$ and $M=M(q)$ are those
appearing in Theorems 1.9 and 1.10, and we put%
\[
\overline{D}=4+2D,\,\,\,z(q)=M\overline{D}(q+\overline{D}).
\]
In this section we establish

\begin{proposition}
\label{Dpowers}Let $N$ be a quasi-semisimple normal subgroup of a group $G$
and $h_{1},\ldots,h_{m}\in G$. Assume that $m\geq z(q)$ and that $\left|
T\right|  >C$ for each non-abelian composition factor $T$ of $N$. Then the
mapping $\psi:N^{(m)}\rightarrow N$ given by%
\[
\prod_{i=1}^{m}(a_{i}h_{i})^{q}=\psi(a_{1},\ldots,a_{m})\cdot\prod_{i=1}%
^{m}h_{i}^{q}.
\]
is surjective.
\end{proposition}

Let $H=\left\langle h_{1},\ldots,h_{m}\right\rangle $. It is clear that $\psi$
depends only on the action of the $h_{i}$ on $N$. The action of $H$ on $N$
lifts to an action on the universal cover $\widetilde{N}$, and it will suffice
to prove the result for the case where $N=\widetilde{N}$. Thus we shall assume
that $N=S_{1}\times\cdots\times S_{r}$ where each $S_{i}$ is a quasisimple
group; the action of $H$ then permutes the $S_{i}$. If $N=N_{1}\times
\cdots\times N_{t}$ where each $N_{i}$ is $H$-invariant, then it is easy to
see that $\psi=\psi_{\mid N_{1}}\times\cdots\times\psi_{\mid N_{t}}$; so we
may assume in addition that this permutation action is transitive. It follows
that $S_{i}\cong S$ for each $i,$ where $S$ is quasisimple with $\left|
S/\mathrm{Z}(S)\right|  >C$.

The explicit expression for $\psi$ is thoroughly unpleasant. Instead of
confronting it directly we proceed as follows. For $\mathbf{x},\,\mathbf{b}\in
N^{(m)}$ and $1\leq i\leq m$ put%
\[
a_{i}(\mathbf{x},\mathbf{b})=x_{i}^{b_{i}}[b_{i},h_{i}^{-1}],
\]
so $a_{i}(\mathbf{x},\mathbf{b})h_{i}=(x_{i}h_{i})^{b_{i}}$. Then%
\begin{align}
\psi(a_{1}(\mathbf{x},\mathbf{b}),\ldots,a_{m}(\mathbf{x},\mathbf{b}))  &
=\prod_{i=1}^{m}\left(  (x_{i}h_{i})^{q}\right)  ^{b_{i}}\cdot\left(
\prod_{i=1}^{m}h_{i}^{q}\right)  ^{-1}\nonumber\\
&  =\prod_{i=1}^{m}\left(  (x_{i}h_{i})^{q}\right)  ^{b_{i}}\cdot\left(
\prod_{i=1}^{m}(x_{i}h_{i})^{q}\right)  ^{-1}\cdot\psi(\mathbf{x})\nonumber\\
&  =\prod_{i=1}^{m}[b_{i},(x_{i}h_{i})^{-q}]^{\tau_{i}(\mathbf{xh})}\cdot
\psi(\mathbf{x}) \label{Psi(x,b)}%
\end{align}
where%
\begin{align*}
\tau_{i}(\mathbf{xh})  &  =(x_{i-1}h_{i-1})^{-q}\ldots(x_{1}h_{1})^{-q}\\
&  =\xi_{i}(x_{1},\ldots,x_{i-1})\tau_{i}(\mathbf{h}),
\end{align*}
say. We shall prove

\begin{proposition}
\label{q-comm}Let $N=S^{(r)}$ where $S$ is a quasisimple group with $\left|
S/\mathrm{Z}(S)\right|  >C$, and let $H=\left\langle k_{1},\ldots
,k_{m}\right\rangle \leq\mathrm{Aut}(N)$ act transitively on the set of simple
factors of $N$. Suppose that $m\geq z(q)$. Then there exist $y_{1}%
,\ldots,y_{m}\in N$ such that
\[
N=\prod_{j=1}^{m}[N,(y_{j}k_{j})^{q}].
\]

\end{proposition}

\noindent(Here and later, we do not distinguish between an element of $N$ and
the inner automorphism it induces).

This suffices to complete the proof of Proposition \ref{Dpowers}. Indeed,
suppose we want to solve the equation $\psi(\mathbf{a})=\kappa$. In view of
(\ref{Psi(x,b)}), it will suffice to find $\mathbf{x}$ and $\mathbf{b}$ such
that%
\begin{equation}
\prod_{i=1}^{m}[b_{i},(x_{i}h_{i})^{-q}]^{\tau_{i}(\mathbf{xh})}=\kappa
\psi(\mathbf{x})^{-1}. \label{x-h-eqn}%
\end{equation}
For each $i$ let $k_{i}\in\mathrm{Aut}(N)$ be induced by $h_{i}^{-\tau
_{i}(\mathbf{h})}$; it is easy to see that then $\left\langle k_{1}%
,\ldots,k_{m}\right\rangle $ is the group of automorphisms induced by
$\left\langle h_{1},\ldots,h_{m}\right\rangle $, hence acts transitively on
the $S_{j}$. Let $y_{1},\ldots,y_{m}\in N$ be as specified in the proposition,
and define $x_{1},\ldots,x_{m}$ recursively by%
\[
x_{i}^{h_{i}}=[h_{i},\xi_{i}^{-1}]\cdot y_{i}^{-\tau_{i}^{-1}\xi_{i}^{-1}}%
\]
where $\xi_{1}=\tau_{1}=1$ and $\xi_{i}=\xi_{i}(x_{1},\ldots,x_{i-1}%
),\,\tau_{i}=\tau_{i}(\mathbf{h})$ for $i>1$. According to the proposition,
there exists $\mathbf{a}\in N^{(m)}$ such that $\prod_{i=1}^{m}[a_{i}%
,(y_{i}k_{i})^{q}]=\kappa\psi(\mathbf{x})^{-1}$. Since $y_{i}k_{i}$ acts as
$(x_{i}h_{i})^{-\tau_{i}(\mathbf{xh})}$ on $N$, we may now solve
(\ref{x-h-eqn}) by setting $b_{i}=a_{i}^{\tau_{i}(\mathbf{xh})^{-1}}$.

\bigskip

The rest of this section is devoted to the proof of Proposition \ref{q-comm}.
As in \S 9, we write $N=S_{1}\times\cdots\times S_{r}$ and fix an
identification of each $S_{i}$ with $S$. For $h\in H$ and $x=(x(1),\ldots
,x(r))\in N$ we write%
\[
x^{h}=\left(  x(^{h}\!1)^{h(1)},\ldots,x(^{h}\!r)^{h(r)}\right)  ,
\]
where $i\mapsto\,^{h}\!i$ is the permutation $\sigma(h^{-1})$ of
$\{1,\ldots,r\}$ induced by the action of $h^{-1}$ on $\{S_{1},\ldots,S_{r}\}$
and $h(i)\in\mathrm{Aut}(S)$ is induced by $h_{\mid S_{(^{h}i)}}:S_{(^{h}%
i)}\rightarrow S_{i}$. For $\Delta\subseteq\{1,\ldots,r\}$ the projection
$N\rightarrow\prod_{i\in\Delta}N_{i}$ is denoted $\pi_{\Delta}$.

The set of fixed points of $\sigma(h)$ is denoted $\mathrm{fix}(h)$, and we
write%
\begin{align*}
\mathrm{fix}^{\ast}(i)  &  =\{j\in\{1,\ldots,m\}\mid i\in\mathrm{fix}%
(k_{j}^{q})\},\\
\lambda(\Delta)  &  =m\left|  \Delta\right|  -\sum_{i\in\Delta}\left|
\mathrm{fix}^{\ast}(i)\right|  .
\end{align*}
Thus $\lambda(\Delta)$ is the number of pairs $(i,j)$ with $i\in\Delta$ such
that $k_{j}^{q}$ moves $S_{i}$.

Put $G_{1}=\left\langle k_{1}^{q},\ldots,k_{m}^{q}\right\rangle $ and let
$\Omega$ be an orbit of $\sigma(G_{1})$ on $\{1,\ldots,r\}$. We say that
$\Omega$ is of \emph{type I }if $\lambda(\Omega)<\overline{D}\left|
\Omega\right|  $, of \emph{type II} otherwise.

When $\Omega$ is of type I there exists at least one $i\in\Omega$ for which
$\left|  \mathrm{fix}^{\ast}(i)\right|  >m-\overline{D}$; we choose such a
value of $i$ and denote it $i_{\Omega}$. Put%
\[
\mathcal{S}=\bigcup\left\{  (\Omega,j)\mid j\in\mathrm{fix}^{\ast}(i_{\Omega
})\right\}
\]
where $\Omega$ ranges over all the $G_{1}$-orbits of type I. Two pairs
$(\Omega,j)$ and $(\Omega^{\prime},j^{\prime})$ will be called
\emph{independent} if either $j\neq j^{\prime}$ or $j=j^{\prime}$ and
$i_{\Omega}$ and $i_{\Omega^{\prime}}$ lie in distinct orbits of $k_{j}$; a
subset of $\mathcal{S}$ is independent if its members are pairwise independent.

\begin{lemma}
\label{indeplemma}Suppose that $m\geq z(q).$ Then for each $G_{1}$-orbit
$\Omega$ of type I there exist an interval $J_{\Omega}\subseteq\mathrm{fix}%
^{\ast}(i_{\Omega})$ and a subset $I_{\Omega}\subseteq J_{\Omega}$ such that

\emph{(i)} $\left|  I_{\Omega}\right|  =M$,

\emph{(ii)} the set%
\[
\mathcal{T}=\left\{  (\Omega,j)\in\mathcal{S}\mid j\in I_{\Omega}\right\}
\]
is independent.
\end{lemma}

This lemma will be proved below. Now let $\Omega$ be a $G_{1}$-orbit of type
I. If $j\in I_{\Omega}$ then $\sigma(k_{j})^{q}$ fixes $i_{\Omega}$, so the
$k_{j}$-cycle $C(\Omega,j)$ of $i_{\Omega}$ has length $e(\Omega,j)$, say,
dividing $q$. Put $q_{\Omega j}=q/e(\Omega,j)$ and let $\beta_{\Omega j}$
denote the automorphism of $S_{i_{\Omega}}$ induced by the action of
$k_{j}^{e(\Omega,j)}$. According to Theorem 1.10, we may choose elements
$x_{\Omega j}\in S_{i_{\Omega}}$ so that%
\begin{equation}
S_{i_{\Omega}}=\prod_{j\in I_{\Omega}}[S_{i_{\Omega}},(x_{\Omega j}%
\beta_{\Omega j})^{q_{\Omega j}}]. \label{S-comm-eq}%
\end{equation}

In this way we obtain a family of elements $x_{\Omega j}\in S_{i_{\Omega}}$ as
$(\Omega,j)$ ranges over $\mathcal{T}$. Now for each $j\in\{1,\ldots,m\}$ let%
\[
y_{j}=\prod_{I_{\Omega}\backepsilon j}x_{\Omega j}\in\prod_{I_{\Omega
}\backepsilon j}S_{i_{\Omega}}%
\]
($y_{j}=1$ if the range of the product is empty). The independence of
$\mathcal{T}$ ensures that if $j\in I_{\Omega}$ then $\pi_{C(\Omega,j)}%
(y_{j})=x_{\Omega j},$ and hence that $(y_{j}k_{j})^{e(\Omega,j)}$ acts on
$S_{i_{\Omega}}$ as $x_{\Omega j}\beta_{\Omega j}$. Thus writing $g_{j}%
=(y_{j}k_{j})^{q}$ for each $j$, we have%
\begin{equation}
S_{i_{\Omega}}=\prod_{j\in I_{\Omega}}[S_{i_{\Omega}},g_{j}]=\prod_{j\in
J_{\Omega}}[S_{i_{\Omega}},g_{j}] \label{S-omega}%
\end{equation}
for each $G_{1}$-orbit $\Omega$ of type I.

\bigskip

Now put $G=\left\langle g_{1},\ldots,g_{m}\right\rangle $, and note that
$\sigma(g_{j})=\sigma(k_{j}^{q})$ for each $j,$ $\sigma(G)=\sigma(G_{1})$. For
each $G$-orbit $\Omega$ let $N_{\Omega}=\prod_{i\in\Omega}S_{i}$. Then $N$ is
the direct product of the $N_{\Omega},$ each of which is invariant under $G$.
Thus to prove Proposition \ref{q-comm} it will suffice to show that%
\begin{equation}
N_{\Omega}=\prod_{j=1}^{m}[N_{\Omega},g_{j}] \label{N-omega}%
\end{equation}
for each $G$-orbit $\Omega$.

\bigskip

\textbf{Case 1}: where $\Omega$ is of type II. Assume for ease of notation
that $\Omega=\{1,\ldots,n\}$. Then%
\[
\lambda(\Omega)=mn-\sum_{i=1}^{n}\left|  \mathrm{fix}^{\ast}(i)\right|  \geq
n\overline{D},
\]
where $\mathrm{fix}^{\ast}(i)=\{j\mid\,^{g_{j}}\!i=i\}$. Note that this
entails $n\geq2$. Let $c(g_{j})$ denote the number of cycles of $\sigma
(g_{j})$. Then%
\begin{align*}
\sum_{i=1}^{m}c(g_{j})  &  \leq\frac{1}{2}\sum_{j=1}^{m}(n+\left|
\mathrm{fix}_{\Omega}(g_{j})\right|  )\\
&  =\frac{1}{2}(mn+\sum_{i=1}^{n}\left|  \mathrm{fix}^{\ast}(i)\right|  )\\
&  =mn-\frac{1}{2}\lambda(\Omega)\leq mn-\frac{1}{2}\overline{D}n=(m-2)n-nD.
\end{align*}
Since $n\geq2$, the identity (\ref{N-omega}) now follows from Proposition 9.1.

\bigskip

\textbf{Case 2}: where $\Omega$ is of type I. Say $\Omega=\{1,\ldots,n\}$, and
that $i_{\Omega}=1$. Suppose that the interval $J_{\Omega}$ is $\{l,l+1,\ldots
,p\}$, so $\sigma(g_{j})$ fixes $1$ for $l\leq j\leq p$.

Given $\kappa\in N_{\Omega}$ we have to solve the equation%
\begin{equation}
\kappa=x_{1}\ldots x_{m} \label{kappa=x}%
\end{equation}
subject to the conditions%
\begin{equation}
x_{i}\in\lbrack N_{\Omega},g_{i}], \label{x(i)}%
\end{equation}
$i=1,\ldots,m.$ Let $\Omega_{i}$ denote the set of orbits of $\sigma(g_{i})$
in $\Omega$; for $\Delta\in\Omega_{i}$ write $k_{\Delta}$ for its first member
and put $n(\Delta)=\left|  \Delta\right|  $. As shown in the proof of
Proposition 9.1, the condition (\ref{x(i)}) is satisfied if and only if%
\[
\prod_{j=0}^{n(\Delta)-1}x_{i}(^{g_{i}^{j}}k_{\Delta})^{g_{i}^{j}(k_{\Delta}%
)}\in\lbrack S,\beta_{i}(\Delta)]
\]
for each $\Delta\in\Omega_{i}$, where $\beta_{i}(\Delta)=g_{i}^{n(\Delta
)}(k_{\Delta})$ is the automorphism of $S$ induced by the action of
$g_{i}^{n(\Delta)}$ on $S_{k(\Delta)}$ (see formula ($H_{i,\Delta}$) in \S 9).

Write (\ref{kappa=x}) as the system of equations $\mathcal{F}$:%
\begin{equation}
\kappa(s)=x_{1}(s)\ldots x_{m}(s), \tag{$F_s$}\label{Fsnew}%
\end{equation}
$s=1,\ldots,n$. For each $s$, let $F_{s}^{\prime}$ be the equation obtained
from $F_{s}$ as follows: for each pair $(i,\Delta)$ with $\Delta\in\Omega_{i}$
and $k_{\Delta}=s$, replace $x_{i}(s)$ by the expression%
\begin{equation}
V_{i}(\Delta)\cdot\left(  \prod_{j=1}^{n(\Delta)-1}x_{i}(^{g_{i}^{j}}%
k_{\Delta})^{g_{i}^{j}(k_{\Delta})}\right)  ^{-1}, \label{subst}%
\end{equation}
where $V_{i}(\Delta)$ is a new symbol. Note that for $i\in J_{\Omega}$ we have
an orbit $\Delta=\{1\}\in\Omega_{i}$, so the first equation becomes%
\begin{equation}
\kappa(1)=\overline{x}_{1}(1)\ldots\overline{x}_{l-1}(1)\cdot\prod_{i=l}%
^{p}V_{i}(\{1\})\cdot\overline{x}_{p+1}(1)\ldots\overline{x}_{m}(1),
\tag{$F_1^\prime$}\label{F'1}%
\end{equation}
where $\overline{x}_{i}(1)$ stands for the expression (\ref{subst}) with
$\Delta$ the $g_{i}$-orbit of $1$.

The resulting system $\mathcal{F}^{\prime}$ of equations contains the unknowns
$V_{i}(\Delta)$ for $\Delta\in\Omega_{i}$ and the $x_{i}(s)$ for every $s$ not
of the form $k_{\Delta}$, $\Delta\in\Omega_{i}$; each such $x_{i}(s)$ now
occurs exactly once with its inverse. We are required to solve $\mathcal{F}%
^{\prime}$ with each $x_{i}(s)\in S$ and each $V_{i}(\Delta)\in\lbrack
S,\beta_{i}(\Delta)]$.

Next, we reduce $\mathcal{F}^{\prime}$ to a single equation using the
procedure described in the proof of Proposition 9.1. That is, if a term
$x=x_{i}(j)^{\pm1}$ appears in $F_{1}^{\prime}$ but its inverse does not, then
$x^{-1}$ appears in some $F_{l}^{\prime}$, $l\neq1$. Solve $F_{l}^{\prime}$
for $x$ and substitute the resulting expression in $F_{1}^{\prime}$; cross out
the equation $F_{l}^{\prime}$, and iterate. As we saw in the preceding
section, the transitivity of $\sigma(G)$ ensures that after $n-1$ such steps
the equations $F_{2}^{\prime},\ldots,F_{n}^{\prime}$ will have been eliminated.

Since the `middle part' of $F_{1}^{\prime}$ is unaffected by this process, the
resulting equation takes the form%
\[
\kappa(1)=A\cdot\prod_{i=l}^{p}V_{i}(\{1\})\cdot B
\]
where $A\cdot B$ is the product, in some order, of certain terms
$x_{i}(j)^{\varepsilon}$, all the$\,V_{i}(\Delta)$ with $\Delta\neq\{1\}$ when
$l\leq i\leq p$, and $\kappa(2)^{-1},\ldots,\kappa(n)^{-1}$, possibly with an
automorphism attached. Setting each such $x_{i}(j)$ and each such
$V_{i}(\Delta)$ equal to $1,$ we are reduced to solving%
\[
\prod_{i=l}^{p}V_{i}(\{1\})=\kappa^{\ast}%
\]
for a certain $\kappa^{\ast}\in S$, subject to the conditions $V_{i}%
(\{1\})\in\lbrack S,\beta_{i}(\{1\})]$ for $l\leq i\leq p$. But $\beta
_{i}(\{1\})$ is just the automorphism induced by the action of $g_{i}$ on
$S_{1}$; the solubility of this equation is therefore assured by
(\ref{S-omega}).

This completes the proof.

\bigskip

It remains to give the\medskip

\textbf{Proof of Lemma \ref{indeplemma}} \ Let $\mathcal{O}$ denote the set of
all $G_{1}$-orbits of type I$.$ For each $\Omega\in\mathcal{O}$ we are given
$i_{\Omega}\in\Omega$ such that $\sigma(k_{j})^{q}$ fixes $i_{\Omega}$ for all
but at most $\overline{D}-1$ values of $j$; thus the set $\mathrm{fix}^{\ast
}(i_{\Omega})$ is the union of at most $\overline{D}$ intervals.

We make the following\medskip

\noindent\emph{Claim}: Let $X\subseteq\{1,\ldots,m\}$ be a subset with
$\left|  X\right|  \geq q+\overline{D}$. Then there exists a mapping
$j:\mathcal{O}\rightarrow X$ such that%
\[
\mathcal{S}_{X}=\left\{  \left(  \Omega,j(\Omega)\right)  \mid\Omega
\in\mathcal{O}\right\}
\]
is an independent subset of $\mathcal{S}$.\medskip

Accepting the claim for now, partition the sequence $\{1,\ldots,m\}$ into
$M\overline{D}$ intervals $X(1),\ldots,X(M\overline{D})$ of length at least
$q+\overline{D}$, and put%
\[
\widetilde{\mathcal{T}}=\bigcup_{i=1}^{M\overline{D}}\mathcal{S}_{X(i)}.
\]
This is evidently an independent set. Now fix $\Omega\in\mathcal{O}$ and
consider the set $\widetilde{\mathcal{T}}_{\Omega}=\{j\mid(\Omega
,j)\in\widetilde{\mathcal{T}}\}$. This meets each $X(i)$, so has cardinality
at least $M\overline{D}$. Since $\widetilde{\mathcal{T}}_{\Omega}%
\subseteq\mathrm{fix}^{\ast}(i_{\Omega})$ it follows that $\left|
\widetilde{\mathcal{T}}_{\Omega}\cap J\right|  \geq M$ for at least one of the
(at most) $\overline{D}$ intervals $J$ that make up $\mathrm{fix}^{\ast
}(i_{\Omega})$. Put $J_{\Omega}=J$ and let $I_{\Omega}$ be any subset of
$\widetilde{\mathcal{T}}_{\Omega}\cap J$ of size $M$. These then satisfy all
the requirements of the lemma.

To prove the \emph{Claim}, we will apply Hall's `marriage theorem' (see e.g.
\cite{PB}, Chapter 22). The `men' are pairs $(\Delta,j)$ where $j\in X$ and
$\Delta$ is an orbit of $\sigma(k_{j})$ with $\left|  \Delta\right|  \mid q$.
The set of `women' is just $\mathcal{O},$ and we say that $\Omega$ `knows'
$(\Delta,j)$ (and vice versa) precisely when $i_{\Omega}\in\Delta$.
\ Evidently each man knows at most $q$ women; while each woman $\Omega$ knows
at least $q$ men, namely the%
\[
(C(\Omega,j),j),\qquad j\in X\cap\mathrm{fix}^{\ast}(i_{\Omega})
\]
where $C(\Omega,j)$ is the $\sigma(k_{j})$-cycle containing $i_{\Omega}$. It
follows (counting possible `couples' in two ways) that for every $n$, each set
of $n$ women collectively knows at least $n$ men.\ Hall's theorem now ensures
that each woman $\Omega$ can find a husband $(\Delta(\Omega),j(\Omega))$ with
$i_{\Omega}\in\Delta(\Omega)$. The monogamy rule means that if $\Omega
\neq\Omega^{\prime}$ then $(\Delta(\Omega),j(\Omega))\neq(\Delta
(\Omega^{\prime}),j(\Omega^{\prime}))$; this is precisely the statement that
the pairs $(\Omega,j(\Omega))$ and $(\Omega^{\prime},j(\Omega^{\prime}))$ are independent.

\section{Equations in semisimple groups, 3: twisted commutators}

Theorem 1.9, stated in the Introduction, asserts that every element of any
finite quasisimple group can be written as a product of boundedly many twisted
commutators. Here we generalize this result. Recall the notation%
\[
T_{\alpha,\beta}(x,y)=x^{-1}y^{-1}x^{\alpha}y^{\beta},
\]
and let $D$ be the absolute constant given in Theorem 1.9.

\begin{proposition}
\label{qss-twisted}Let $N$ be a quasi-semisimple group and $\alpha_{l}%
,\beta_{l}\ \ $($l=1,2,...,D$) arbitrary automorphisms of $N$. Then
\begin{equation}
\prod_{i=1}^{D}T_{\alpha_{i},\beta_{i}}(N,N)=N. \label{T-equn}%
\end{equation}

\end{proposition}

The universal cover of $N$ is a direct product $\widetilde{N}=S_{1}%
\times\cdots\times S_{n}$ where each $S_{i}$ is a quasisimple group of
universal type. Each automorphism of $N$ lifts to one of $\widetilde{N}$, so
it will suffice to prove the result in the case $N=\widetilde{N}$, which we
assume henceforth.

Let $G=\left\langle \alpha_{l},\beta_{l}\mid1\leq l\leq D\right\rangle $ be
the subgroup of $\mathrm{Aut}(N)$ generated by the given automorphisms. Then
$G$ permutes the factors $S_{1},\ldots,S_{n}$, with orbits $\Lambda_{i}$ say.
Now $N$ is the direct product of the subgroups $N(i)=\prod_{j\in\Lambda_{i}%
}S_{j}$, on each of which $G$ acts by restriction, and it will suffice to
prove (\ref{T-equn}) with $N(i)$ in place of $N$, for each $i$. Thus we may,
and shall, assume that the permutation action of $G$ on $\{S_{1},\ldots
,S_{n}\}$ is transitive. As in the preceding sections, we shall write%
\[
S_{i}^{\alpha^{-1}}=S_{^{\alpha}i}\,(\alpha\in G,\,1\leq i\leq n).
\]
Since this action of $G$ is transitive, the groups $S_{i}$ are all isomorphic;
we fix an identification of each $S_{i}$ with a fixed quasisimple group $S$.
Thus elements of $N$ will be written in the form%
\[
x=(x(1),x(2),\ldots,x(n))
\]
with each $x(i)\in S$, and the action of $G$ takes the form%
\[
x^{\alpha}=(x(^{\alpha}1)^{\alpha(1)},x(^{\alpha}2)^{\alpha(2)},\ldots
,x(^{\alpha}n)^{\alpha(n)});
\]
here $\alpha(1),\ldots,\alpha(n)\in\mathrm{Aut}(S)$ depend on $\alpha\in G$
(and the fixed identifications $S_{i}\rightarrow S$). We put%
\[
\Gamma=\left\langle \alpha_{i}(s),\,\beta_{i}(s)\mid1\leq i\leq D,\,1\leq
s\leq n\right\rangle \leq\mathrm{Aut}(S).
\]

Let $\kappa\in N$. We have to show that there exist $\kappa_{1},\ldots
,\kappa_{D}\in N$ such that%
\begin{equation}
\kappa=\kappa_{1}\ldots\kappa_{D} \label{com}%
\end{equation}
and such that for each $i$ there exist $x_{i},\,y_{i}\in N$ with%
\begin{equation}
\kappa_{i}=x_{i}{}^{-1}y_{i}{}^{-1}x_{i}{}^{\alpha_{i}}y_{i}{}^{\beta_{i}}.
\label{kappa(i)}%
\end{equation}

To begin with, we fix $i$ and analyse the equation (\ref{kappa(i)}). This
equation is solvable in $N$ if and only if there exist elements $x_{i}%
(s),\,y_{i}(s)\in S$ ($s=1,\ldots,n$) such that (writing $\alpha=\alpha
_{i},\,\beta=\beta_{i}$)
\begin{equation}
\kappa_{i}(s)=x_{i}(s)^{-1}y_{i}(s)^{-1}x_{i}(^{\alpha}s)^{\alpha(s)}%
y_{i}(^{\beta}s)^{\beta(s)} \tag{$E_s$}\label{Es}%
\end{equation}
holds for $s=1,\ldots,n$. We consider $\mathcal{E}=(E_{1},\ldots,E_{n})$ as a
system of simultaneous equations in the unknowns $x_{i}(s),\,y_{i}(s)$.

Put%
\[
G_{i}=\left\langle \alpha_{i},\beta_{i}\right\rangle
\]
and let $\Omega=\Omega_{i}$ denote the set of orbits of $G_{i}$ on
$\{1,\ldots,n\}$. The system $\mathcal{E}$ breaks up into $\left|
\Omega\right|  $ independent systems of equations, one for each orbit
$\Delta\in\Omega:$%
\[
\mathcal{E}_{\Delta}=(E_{s})_{s\in\Delta}.
\]
We fix an orbit $\Delta$ of size $n_{\Delta}$, and introduce the alphabet
$X\cup P$ where
\begin{align*}
X  &  =X_{\Delta}=\{x_{i}(s),\,y_{i}(s)\mid s\in\Delta\}\\
P  &  =P_{\Delta}=\{\kappa_{i}(s)\mid s\in\Delta\};
\end{align*}
for now the elements of $X\cup P$ are considered as abstract symbols. Here $P$
is the set of \emph{parameters }and\emph{\ }$X$ is the set of \emph{variables}%
. Let $F_{\Delta}=F_{\Gamma}(X_{\Delta}\cup P_{\Delta})$ be the free $\Gamma
$-group, defined in Section 8. For $x\in X^{\pm1}\cup P^{\pm1}$ we will write
$x^{\ast}$ to denote an arbitrary element of the form $x^{\gamma}$, $\gamma
\in\Gamma$.

We consider the right-hand sides of the equations $E_{s}$ in $\mathcal{E}%
_{\Delta}$ as words on the alphabet $X^{\Gamma}\cup P^{\Gamma}$. Note that
each variable occurs exactly once, as does its inverse, in the system
$\mathcal{E}_{\Delta}$.

Let $k=k_{\Delta}$ be the first symbol in $\Delta$. We shall modify
$\mathcal{E}_{\Delta}$ by the familiar process of eliminating variables.
Suppose that $x^{\ast}$ occurs in the first equation $E_{k}$, where $x\in
X\cup X^{-1}$; then a term $x^{-\ast}$ occurs in some equation $E_{l}$. If
$l\neq k$, solve $E_{l}$ for $x$ and substitute the resulting value of $x$ in
$E_{k}$. Let us call this process a \emph{substitution }$(l\rightarrow k)$.
Each substitution reduces by one both the number of variables and the number
of equations in the system $\mathcal{E}_{\Delta}$.

We claim that it is possible to apply $n_{\Delta}-1$ substitutions and thus
reach an equivalent system consisting of the single equation
\begin{equation}
\kappa_{i}(k_{\Delta})=U_{\Delta}, \label{one}%
\end{equation}
where $U_{\Delta}$ is a certain word on $X^{\Gamma}\cup P^{\Gamma}$. This
follows just as in \S 9 from the fact that $G_{i}$ acts transitively on
$\Delta$.

As in \S 8, let $\mathcal{M}=\mathcal{M}_{\Delta}$ denote the the free monoid
on $X^{\pm\Gamma}\cup P^{\pm\Gamma}$. Recall that for words $U$,$\,U^{\prime
}\in\mathcal{M}$ the expression $U=_{F}U^{\prime}$means that $U$ and
$U^{\prime}$ take the same value $\overline{U}$ in the group $F$, and that
$\widehat{U}$ denotes the word in $W$, the free monoid on $X\cup X^{-1}$,
obtained from $U$ when all symbols from $P^{\pm\Gamma}$ are deleted and
$x^{\gamma}$ is replaced by $x$ for each $x\in X\cup X^{-1}$,\thinspace
\ $\gamma\in\Gamma$.

\begin{lemma}
\label{U_delta}There exist $x=x_{\Delta},\,y=y_{\Delta},\,V=V_{\Delta}%
\in\mathcal{M}$ and $a=a_{\Delta},\,b=b_{\Delta}\in\Gamma$ such that%
\[
U_{\Delta}=_{F}T_{a,b}(x,y)\cdot V,
\]
the family $\{\overline{x},\,\overline{y}\}\cup\sup(\widehat{V})\cup P$ is
independent, each $\kappa_{i}(s)$ for $s\in\Delta\setminus\{k_{\Delta}\}$
occurs exactly once in $V$ with exponent $-\gamma$ for some $\gamma\in\Gamma$,
and $\kappa_{i}(k_{\Delta})$ does not occur in $V$.
\end{lemma}

\begin{proof}
Let $U_{1}$ be the word on the right-hand side of $E_{k}$, and let $U_{f}$
denote the word obtained from $U_{1}$ after $f-1$ substitutions have been
carried out. A substitution $(l\rightarrow k)$ has one of the following
effects (we drop the subscript $i$ from $x(s),\,y(s),$ $\alpha$ and $\beta$):%
\[%
\begin{tabular}
[c]{lllll}%
\underline{$l$} &  & \underline{$U_{f}$} &  & \underline{$U_{f+1}$}\\
$s$ &  & $A\cdot x(s)^{\ast}\cdot B$ &  & $A\cdot y(l)^{-\ast}x(^{\alpha
}l)^{\ast}y(^{\beta}l)^{\ast}\kappa_{i}(l)^{-\ast}\cdot B$\\
$\!\!\!\!\!\!^{\alpha^{-1}}\!\!s$ &  & $A\cdot x(s)^{-\ast}\cdot B$ &  &
$A\cdot y(^{\beta}l)^{\ast}\kappa_{i}(l)^{-\ast}x(l)^{-\ast}y(l)^{-\ast}\cdot
B$\\
$s$ &  & $A\cdot y(s)^{\ast}\cdot B$ &  & $A\cdot x(^{\alpha}l)^{\ast
}y(^{\beta}l)^{\ast}\kappa_{i}(l)^{-\ast}x(l)^{-\ast}\cdot B$\\
$\!\!\!\!\!\!^{\beta^{-1}}\!\!s$ &  & $A\cdot y(s)^{-\ast}\cdot B$ &  &
$A\cdot\kappa_{i}(l)^{-\ast}x(l)^{-\ast}y(l)^{-\ast}x(^{\alpha}l)^{\ast}\cdot
B$%
\end{tabular}
\]

Since each variable occurs exactly once with its inverse in the system
$\mathcal{E}_{\Delta},$ it is easy to see that the same holds for the final
word $U_{\Delta}=U_{n_{\Delta}}$, except for the $n_{\Delta}-1$ matching pairs
$x,\,x^{-1}$ that have been eliminated. Thus the word $\widehat{U_{\Delta}}$
is balanced.

We claim also that $\widehat{U_{\Delta}}\neq_{F}1$. To see this, let
$\Phi=\left\langle \xi,\,\eta\right\rangle $ be a the free nilpotent group of
class two on two free generators, and define a (monoid) homomorphism
$\theta:W\rightarrow\Phi$ by%
\[
x(s)^{\varepsilon}\mapsto\xi^{\varepsilon},\,y(s)^{\varepsilon}\mapsto
\eta^{\varepsilon}%
\]
($\varepsilon=\pm1,$ $s\in\Delta$). Now we can write $\theta(\widehat{U_{f}%
})=\theta(\widehat{A})\theta(z)\theta(\widehat{B})$ where $z$ is one of
$x(s)^{\pm1},\,y(s)^{\pm1}$, and we see that in the four cases listed we get,
respectively,%
\begin{align*}
\theta(\widehat{U_{f+1}})  &  =\theta(\widehat{A})\theta(z)\cdot\lbrack
\xi,\eta]\cdot\theta(\widehat{B})\\
\theta(\widehat{U_{f+1}})  &  =\theta(\widehat{A})\theta(z)\cdot\lbrack
\xi^{-1},\eta^{-1}]\cdot\theta(\widehat{B})\\
\theta(\widehat{U_{f+1}})  &  =\theta(\widehat{A})\theta(z)\cdot\lbrack
\eta,\xi^{-1}]\cdot\theta(\widehat{B})\\
\theta(\widehat{U_{f+1}})  &  =\theta(\widehat{A})\theta(z)\cdot\lbrack
\eta^{-1},\xi]\cdot\theta(\widehat{B})
\end{align*}
each of which is equal to $[\xi,\eta]\theta(\widehat{U_{f}})$. As
$\theta(\widehat{U_{1}})=[\xi,\eta]$ it follows that $\theta(\widehat
{U_{\Delta}})=[\xi,\eta]^{n_{\Delta}}\neq1$. Since $\theta$ factors through
$F$ this establishes the claim.

Thus $U_{\Delta}$ satisfies the conditions of Proposition 8.2. This now gives
the result, provided only that the multiplicities of the $\kappa_{i}(s)$ in
$U_{\Delta}$ are as described in the statement. But this is clear, since each
substitution $(l\rightarrow k)$ as above introduces the term $\kappa
_{i}(l)^{-\ast}$ (and no other terms from $P\cup P^{-1}$), and the label $l$
runs over the set $\Delta\setminus\{k_{\Delta}\}$ as $f$ goes from $1$ to
$n_{\Delta}-1$.
\end{proof}

\bigskip

The preceding reduction now shows that the equation (\ref{kappa(i)}) is
solvable in $N$ if and only if for each orbit $\Delta\in\Omega_{i}$ there
exists a $\Gamma$-homomorphism%
\[
\phi_{\Delta}:F_{\Delta}\rightarrow S
\]
sending each symbol $\kappa_{i}(s)$ to the element with the same name in $S$
and satisfying%
\[
\phi_{\Delta}(T_{a_{\Delta},b_{\Delta}}(x_{\Delta},y_{\Delta})V_{\Delta
})=\kappa_{i}(k_{\Delta}).
\]

Now put $Z_{\Delta}=\sup(\widehat{V_{\Delta}})$, and consider a new alphabet
$Y\cup P\cup\mathcal{K}$ where%
\begin{align*}
Y  &  =\bigcup_{i=1}^{D}\bigcup_{\Delta\in\Omega_{i}}Z_{\Delta}\cup
\{x_{\Delta},y_{\Delta}\}\\
P  &  =\bigcup_{i=1}^{D}\bigcup_{\Delta\in\Omega_{i}}P_{\Delta}=\{\kappa
_{i}(s)\mid1\leq i\leq D,\,1\leq s\leq n\}\\
\mathcal{K}  &  =\{\kappa(1),\ldots,\kappa(n)\}.
\end{align*}

The equation (\ref{com}) is equivalent to the system of equations
$\mathcal{F}=(F_{1},\ldots,F_{n}):$%
\begin{equation}
\kappa(s)=\kappa_{1}(s)\ldots\kappa_{D}(s). \tag{$F_s$}%
\end{equation}
For each pair $(i,\Delta)$ with $\Delta\in\Omega_{i}$ we substitute the
expression $T_{a_{\Delta},b_{\Delta}}(x_{\Delta},y_{\Delta})V_{\Delta}$ for
$\kappa_{i}(k_{\Delta})$ in the equation $F_{k_{\Delta}}$, to obtain a system
$\mathcal{F}^{\prime}=(F_{1}^{\prime},\ldots,F_{n}^{\prime})$:%
\begin{equation}
\kappa(s)=W_{s} \tag{$F_s^\prime$}\label{F's}%
\end{equation}
where $W_{s}$ is a certain word on the alphabet $Y^{\Gamma}\cup P^{\Gamma}$.
Now recall that $V_{\Delta}$ contains $\kappa_{i}(s)^{-\ast}$ exactly once for
each $s\in\Delta\setminus\{k_{\Delta}\},$ and no other terms from
$P^{\pm\Gamma}$; it follows that $W_{1}W_{2}\ldots W_{n}$ contains the terms
$\kappa_{i}(s),\,\kappa_{i}(s)^{-\ast}$ once each whenever $s\notin
\{k_{\Delta}\mid\Delta\in\Omega_{i}\},$ and no other terms from $P^{\pm\Gamma
}$.

We now repeat the elimination procedure used above. Suppose that $\mu\in P\cup
P^{-1}$ and $\mu^{\ast}$ occurs in $F_{1}^{\prime}$ while $\mu^{-\ast}$ occurs
in $F_{l}^{\prime}$ for some $l\neq1$. Solve $F_{l}^{\prime}$ for $\mu$ and
substitute the resulting expression into $F_{1}^{\prime}.$ It is easy to see
that two equations $F_{p}^{\prime},\,F_{q}^{\prime}$ are `linked', in the
sense that they share a parameter from $P$, if and only if there exists $i$
such that $p$ and $q$ lie in the same orbit of $G_{i}$. Since the $G_{i}$
generate $G$ which is transitive on $\{1,\ldots,n\},$ it follows as before
that we can perform $n-1$ such substitutions and obtain an equivalent system
consisting of one equation%
\begin{equation}
\kappa(1)=V. \label{k(1)=V}%
\end{equation}
Each substitution $(l\rightarrow1)$ eliminates a pair $\kappa_{i}%
(l),\,\kappa_{i}(l)^{-1}$ and introduces into the right-hand member of
$F_{1}^{\prime}$ both a term $\kappa(l)^{-1}$ and all the terms $\ T_{\Delta
}=T_{a_{\Delta},b_{\Delta}}(x_{\Delta},y_{\Delta})$ that appear in
$F_{l}^{\prime}$ (ignoring exponents from $\Gamma$). It follows that $V$
contains each of the terms $\kappa(2)^{-1},\ldots,\kappa(n)^{-1}$ exactly
once, and each of the terms $T_{\Delta}$ ($\Delta\in\Omega_{i},\,1\leq i\leq
D$) exactly once. The other factors of $V$ (still ignoring exponents from
$\Gamma$) all belong to $P^{\pm1}\cup\bigcup_{\Delta}Z_{\Delta}^{\pm1}$.

Let $\mathcal{X}=\{x_{i}(s)\mid1\leq i\leq D,\,1\leq s\leq n\}$. Recall now
(Lemma \ref{U_delta}) that each of the families $\{\overline{x_{\Delta}%
},\,\overline{y_{\Delta}}\}\cup Z_{\Delta}\cup P_{\Delta}$ is independent in
the free $\Gamma$-group $F_{\Delta}$. This implies that the family%
\[
\bigcup_{\Delta\in\Omega_{i},\,1\leq i\leq D}\left(  \left\{  \overline
{x_{\Delta}},\overline{y_{\Delta}}\right\}  \cup Z_{\Delta}\right)  \cup
P\cup\mathcal{K}%
\]
is independent in the free $\Gamma$-group $F$ on $\mathcal{X}\cup
P\cup\mathcal{K}$. Hence for any choice of elements $\xi_{\Delta}%
,\,\eta_{\Delta}\in S$ there is a $\Gamma$-equivariant homomorphism $\phi
_{\xi,\eta}:F\rightarrow S$ sending $x_{\Delta}$ to $\xi_{\Delta}$,
$y_{\Delta}$ to $\eta_{\Delta}$, each symbol $\kappa(i)$ to the given element
$\kappa(i)$ of $S$, and each term of $P\cup\bigcup Z_{\Delta}$ that appears in
$V$ to $1$. Then%
\[
\phi_{\xi,\eta}(\kappa(1))=\kappa(1),
\]
while
\[
\phi_{\xi,\eta}(V)=h_{0}\prod_{\Delta\in\Omega_{i},\,1\leq i\leq
D}T_{a_{\Delta},b_{\Delta}}(\xi_{\Delta},\eta_{\Delta})^{\gamma_{\Delta
}h_{\Delta}}%
\]
(in some order) where the $\gamma_{\Delta}\in\Gamma$ and $h_{0},\,h_{\Delta
}\in S$ do not depend on $\xi,\,\eta$.

Using the identity%
\[
T_{a,b}(x,y)^{\gamma}=T_{a^{\gamma},b^{\gamma}}(x^{\gamma},y^{\gamma})
\]
we rewrite the above as%
\begin{equation}
h_{1}^{-1}\phi_{\xi,\eta}(V)=\prod_{\Delta\in\Omega_{i},\,1\leq i\leq
D}T_{a_{\Delta}^{\prime},b_{\Delta}^{\prime}}(\xi_{\Delta}^{\prime}%
,\eta_{\Delta}^{\prime}) \label{T-product}%
\end{equation}
where $a_{\Delta}^{\prime},b_{\Delta}^{\prime}\in\mathrm{Aut}(S)$ and
$\xi_{\Delta}^{\prime},\eta_{\Delta}^{\prime}$ are the images of $\xi_{\Delta
},\eta_{\Delta}$ under certain fixed automorphisms of $S$.

Now Theorem 1.9 asserts that $S=\prod T_{a_{\Delta}^{\prime},b_{\Delta
}^{\prime}}(S,S)$ provided there are at least $D$ factors in the product.
Hence we can choose values for $\xi_{\Delta},\eta_{\Delta}$ in $S$ so that the
product on the right of (\ref{T-product}) takes the value $h_{1}^{-1}%
\kappa(1)$.

The original equation (\ref{com}) is now solved subject to the conditions
(\ref{kappa(i)}) by giving each unknown $x_{i}(s)$ the value $\phi_{\xi,\eta
}(x_{i}(s))$. This completes the proof of Proposition \ref{qss-twisted}.

\texttt{\bigskip}

\texttt{Nikolay Nikolov}

\texttt{New College}

\texttt{Oxford OX1 3BN}

\texttt{UK}.

\bigskip

\texttt{Dan Segal}

\texttt{All Souls College}

\texttt{Oxford OX1 4AL}

\texttt{UK.}

\end{document}